\documentclass[11pt,a4paper]{article}
\usepackage{amsmath}
\usepackage{amssymb}
\usepackage{color}

\usepackage{graphicx}
\usepackage{amsfonts,amsmath, amssymb}

\numberwithin{equation}{section}

\usepackage{graphics}
\usepackage{epsfig}
\newtheorem{claim}{\bf \t}[part]

\newtheorem{theorem}{Theorem}[section]
\newtheorem{corollary}[theorem]{Corollary}
\newtheorem{lemma}[theorem]{Lemma}

\newtheorem{remark}[theorem]{Remark}

\def\v{\varepsilon}

\def\t{\theta}
\def\k{\kappa}

\def\a{\alpha}
\def\b{\beta}

\def\d{\delta}
\def\l{\lambda}

\def\r{\rho}

\def\z{\zeta}

\def\di{\displaystyle}
\def\i{\infty}
\def\f{\frac}
\begin{document}

\title{ Justification of Diffusion limit for the Boltzmann Equation with a non-trivial Profile}

\author{ {Feimin Huang$^{a,b}$, \ Yi Wang$^{a,b}$, \
 Yong Wang$^{a}$\footnote{Corresponding author.  \newline \indent
Email addresses: fhuang@amt.ac.cn(Feimin Huang),
wangyi@amss.ac.cn(Yi Wang), yongwang@amss.ac.cn(Yong Wang), matyang@cityu.edu.hk(Tong Yang)},  \ and   Tong Yang$^{c}$}
\\
\ \\
   {\small \it $^a$Institute of Applied Mathematics, AMSS, CAS, Beijing 100190, China}\\
\  \\
   {\small \it $^b$Beijing Center of Mathematics and Information Sciences, Beijing
100048, P.R.China } \\
 {\small and} \\
  {\small \it $^c$Department of Mathematics, City University of
Hong Kong, Hong Kong}\\
}

\date{ }
\maketitle

\begin{abstract}

Under the diffusion scaling and a scaling assumption on the
microscopic component, a non-classical fluid dynamic system was derived in
\cite{BGLY} that is related to the system of ghost effect derived in \cite{Sone-2} in different settings.
This paper aims to justify this limit system for a non-trivial background
profile with slab symmetry. The result reveals
not only the  diffusion phenomena in the temperature and density,
but also the flow of higher order in Knudsen number due to the
gradient of the   temperature.  Precisely,
we show that the solution to the Boltzmann equation  converges to a diffusion wave with  decay rates in both Knudsen number and time.

\

Keywords: Boltzmann equation,  Knudsen number, diffusive scaling, diffusion wave

\

AMS: 35Q35, 35B65, 76N10
\end{abstract}

\tableofcontents

\section{Introduction}

Consider the Boltzmann equation with slab symmetry under the diffusive scaling
\begin{equation}\label{1.1}
\v \partial_tf^\v+\xi_1 f^\v_x
	=\frac {1}{\v}Q(f^\v,f^\v),\quad (t,x,\xi)\in{\mathbb R}_+
	\times{\mathbb R}\times{\mathbb R}^3.
\end{equation}
Here $f^\v(t,x,\xi)\ge 0$ is the distribution density of particles at $(t,x)$ with velocity $\xi$,
$Q(f,f)$ is  the collision operator which
is a non-local bilinear operator in the velocity variable with a
kernel determined by the physics of particle
interaction. For monatomic
gas, the rotational invariance of the particle leads to the
collision operator $Q(f, f)$ as a bilinear
collision operator in the form of, cf. \cite{Boltzmann}:
\begin{equation*}\label{B1.2}
Q(f,g)(\xi) \equiv \frac{1}{2}
\int_{\mathbb{R}^3}\int_{{\mathbb S}_+^2}
\Big(f(\xi')g(\xi_*')+ f(\xi_*')g(\xi')-f(\xi) g(\xi_*) - f(\xi_*)g(\xi)
\Big) B(|\xi-\xi_*|, \theta)d \xi_* d \Omega,
\end{equation*}
with $\theta$ being the angle between the relative velocity and the
unit vector $\Omega$. Here ${\mathbb S}^2_+=\{\Omega\in {\mathbb
	S}^2:\ (\xi-\xi_*)\cdot \Omega\geq 0\}$. The conservation of
momentum and energy gives the following relation between velocities
before and after collision:
$$
\left\{
\begin{array}{l}
\xi'= \xi -[(\xi-\xi_*)\cdot \Omega] \; \Omega, \\[3mm]
\xi_*'= \xi_* + [(\xi-\xi_*)\cdot \Omega] \; \Omega.
\end{array}
\right.
$$

In this paper, we will consider the two
basic models, i.e., the hard sphere model and the hard potential
with angular cut-off, for which the collision kernel
$B(|\xi-\xi_*|,\theta)$ takes the form of
$$
B(|\xi-\xi_*|,\theta)=|(\xi-\xi_*, \Omega)|,
$$
and
$$
B(|\xi-\xi_*|,\theta)=|\xi-\xi_*|^{\frac{p-5}{p-1}}b(\theta),\quad
b(\theta)\in L^1([0, \pi]),~p\geq5,
$$
respectively. Here, $p$ is the index in potential of
the inverse power law
 that is  proportional to $r^{1-p}$ with $r$ being the distance
between two particles.

Motivated by \cite{Sone-2},  the following macroscopic and microscopic decomposition with scalings was introduced in \cite{BGLY}:
\begin{equation}\label{1.3}
f^\v={M}_{[\rho^\v,\v u^\v, \theta^\v]} +\v G^\v.
\end{equation}
Here $M_{[\rho^\v,\v u^\v, \theta^\v]}$ is the local Maxwellian and
$G^\v$ is the microscopic component. Moreover,
the  local Maxwellian $M_{[\rho^\v,\v u^\v, \theta^\v]}$ is defined by the
five conserved quantities, that is,
the mass density $\rho^\v(t,x)$, momentum density
$m^\v(t,x)=\v\rho^\v(t,x) u^\v(t,x)$ and energy
density $e^\v(t,x)+\frac{1}{2}|\v u^\v(t,x)|^2$
given by
\begin{equation}\label{B1.4}
\left\{
\begin{array}{l}
\rho^\v(t,x) \equiv {\displaystyle\int_{{ \mathbb R}^3}} f^\v(t,x,\xi) d \xi, \\[3mm]
m_i^\v(t,x) \equiv {\displaystyle\int_{{\mathbb R}^3}} \psi_i(\xi)
f^\v(t,x,\xi)d \xi \  {\textrm {for} }\  i=1,2,3, \\[3mm]
\left[\rho^\v\left({ e^\v}+\frac{\v^2}{2}  |u^\v|^2 \right)\right](t,x) \equiv
{\displaystyle\int_{{\mathbb R}^3}} \psi_4(\xi) f^\v(t,x,\xi) d \xi,
\end{array}
\right.
\end{equation}
as
\begin{equation}\label{B1.4-1}
{ M}\equiv { M}_{[\rho^\v,\v u^\v,\theta^\v]} (t,x,\xi) \equiv
\frac{\rho^\v(t,x)}{\sqrt{ (2 \pi R \theta^\v(t,x))^3}}
\exp\left(-\frac{|\xi -\v u^\v(t,x)|^2}{2R\theta^\v(t,x)}\right).
\end{equation}
Here   $\psi_\alpha(\xi)$ are the collision invariants:
$$
\left\{
\begin{array}{l}
\psi_0(\xi) \equiv 1, \\[2mm]
\psi_i(\xi) \equiv \xi_i \ \ {\textrm {for} }\ \  i=1,2,3, \\[2mm]
\psi_4(\xi) \equiv \frac{1}{2} |\xi|^2,
\end{array}
\right.
$$
satisfying
$$
\int_{{\mathbb{R}}^3} \psi_j(\xi) Q(h,g) d \xi =0,\quad {\textrm {for} }
\ \  j=0,1,2,3,4.
$$
Here, $\theta^\v$ is the temperature  related to the internal energy
$e^\v$ by  $e^\v=\frac{3}{2} R\theta^\v$
with $R$ being the gas constant, and $\v u^\v$ is the bulk velocity.  Note
that even though $u^\v$ is of higher order, it
is the scaled velocity that appears in the equations for
the macroscopic variables $\rho^\v$ and $\theta^\v$.

The Boltzmann equation is a  fundamental equation in statistical physics for
 rarefied gas which  describes
the time evolution of  particle distribution.
There has been tremendous progress  on the mathematical theories
for the Boltzmann equation with $\v$ being a fixed constant, such as the global existence
of weak (renormalized) solution for large data in \cite{Diperna-Lions}
and classical
solutions as  small perturbations of  equilibrium states (Maxwellian) in
\cite{Guo,Liu-Yang-Yu, Ukai-1}  and the references therein, etc.

On the other hand, the study on
the hydrodynamic limit of Boltzmann equation is important
and challenging. For this,  it is well known that the classical works of Hilbert, Chapman-Enskog reveal the  relation of the
Boltzmann equation to the classical systems of fluid dynamics through
asymptotic
expansions with respect to the  Knudsen number.
For the hydrodynamic limit of Boltzmann equation to the compressible Euler system,  we refer  \cite{BGL, BGLY} for the formal derivation. If the Euler system is assumed to have smooth solution,  this hydrodynamic limit is proved rigorously in  \cite{Ukai-3,Caflisch} with and without initial layer respectively.

However, it is well known that the compressible Euler system  develops singularity in finite time even for sufficiently smooth initial data. The Riemann problem is the  basic problem to the compressible Euler system, and its solution turns out to be  fundamental in the
theory of hyperbolic conservation laws because it not only captures the local and
global behavior of solutions but also reveals the effect of nonlinearity in
the structure of the solutions. There are three basic wave patterns for the Euler system, that is, shock wave, rarefaction wave,
and contact discontinuity. For  the hydrodynamic limit of the Boltzmann equation
in the setting of Riemann solutions, we refer  \cite{ Huang-Wang-Wang-Yang, Huang-Wang-Yang,  Huang-Wang-Yang-2, Xin-Zeng, Yu}.

Under the diffusive scaling, usually, the  density function $f^\v(t,x,\xi)$ is
set as a perturbation of a global Maxwellian $M_{[1,0,1]}$, i.e.
\begin{align}\label{expansion}
f^\v(t,x,\xi)=M_{[1,0,1]}+  M_{[1,0,1]}\Big(\v f_1(t,x,\xi)+\cdots+\v^n f_n^\v(t,x,\xi)\Big). 
\end{align}
There has been extensive study on the hydrodynamic limit $\v\rightarrow0$ of the Boltzmann equation  to the incompressible Navier-Stokes-Fourier system,
for example, to justify the DiPerna-Lions' renormalized solution in \cite{Diperna-Lions} of the Boltzmann equation to the Leray-Hopf weak solutions of the incompressible Navier-Stokes-Fourier system. For this,
Bardos-Golse-Levermore \cite{BGL} first studied this problem under certain a priori assumption.  Recently, a breakthrough was achieved by Golse-Raymond in \cite{GS1}   which  established a proof of such  limit for certain class of collision kernels. After that, some progress was made for more general collision kernels,
cf. \cite{Levermore}. In fact, there are also  a lot of
   important contributions on this problem over the years, see \cite{GL,GS,L-M,L-M1,Masmoudi,Masmoudi-S,Saint} and the references therein.

In the framework of classical solutions to the incompressible Navier-Stokes-Fourier system, it was proved in \cite{Esposito1} that one can find a  Boltzmann solution $f^\v(t,x,\xi)$ such that $f^\v_2$ is of order $\v^2$, but it is not clear about the amplitude $f^\v_2$ at the initial time. Later,  the Navier-Stokes-Fourier limit was proved for $f^\v(0,x,\xi)$
with small data in \cite{BU}. Recently, Guo in \cite{Guo1} justified the diffusive expansion \eqref{expansion} when $f_1(0,x,\xi)$ has small amplitude while $f_i^\v(0,x,\xi)$ can have arbitrarily large amplitude for $i\geq2$ in
a torus. This work was later generalized to some
other settings, cf. \cite{Liu-Zhao,Jang,Xiong-Jiang}. Moreover, based on the $L^2-L^\infty$ estimate,  Esposito-Guo-Kim-Marra \cite{Guo2} proved the hydrodynamic limit of the rescaled Boltzmann equation  to the incompressible Navier-Stokes-Fourier system in a bounded domain if the initial data is small.

Notice that all the results under the diffusive scaling mentioned above are either about large perturbation of
vacuum or small perturbation of a  global Maxwellian.
A natural question to ask is how about the perturbation of a non-trivial
profile. The purpose of this paper is to study this problem in
the setting of \eqref{1.3}.

In fact, under the assumption \eqref{1.3},
when $\v\rightarrow 0$,  formally we have
\begin{align}\label{diffusive}
f^\v={M}_{[\rho^\v,\v u^\v, \theta^\v]} +\v G^\v\rightarrow{M}_{[\rho, 0, \theta]},
\end{align}
which shows that in the macroscopic level,  only
the unknown limit functions $\r,\t$
survive because the macroscopic velocity is zero.
However, as shown in \cite{BGLY} and will be recalled
in the next section, the equations of $\r$ and $\t$  are actually closely related to the scaled velocity $u$. Indeed, this diffusive scaling  induces diffusion phenomenon for both the temperature $\t$ and density $\r$, and
the non-zero gradient of temperature  induces a non-trivial flow in the higher order along the same direction.
In this paper, we will construct such diffusion wave and study the hydrodynamic limit of the rescaled Boltzmann equation to such a diffusion wave global in time.

The rest of the paper will be organized  as follows.  The construction of the diffusion wave and the main theorem will be given in the next section.
 We will reformulate the problem and derive some a priori estimates in Section 3. Based on the a priori estimates, the main theorem will be proved in Section 4.

\

\noindent\textbf{Notations:} Throughout this paper, the positive
generic constants that are independent of $\v$ are denoted by
$c,C,C_i(i=1,2,3,\cdots)$. And we will use $\|\cdot\|$ to denote the standard
$L_2(\mathbb{R};dz)$ norm, and $\|\cdot\|_{H^i}~(i=1,2,3,\cdots)$ to
denote the standard Sobolev $H^i(\mathbb{R};dz)$ norm with $z=x$ or $y$. Sometimes, we also use
$O(1)$ to denote a uniform bounded constant independent of
$\v$.

\section{Construction of Profile and the Main Result}

We will drop the superscript $\v$  in the case of no confusion for simple notation. The inner product of $h,\ g$ in  $L^2_{\xi}({\mathbb
R}^3)$ with respect to a given Maxwellian $\tilde{ M}$ is defined
by:
$$
 \langle h,g\rangle_{\tilde{ M}}\equiv \int_{{\mathbb R}^3}
 \frac{1}{\tilde{ M}}h(\xi)g(\xi)d \xi,
$$
when the  integral is well defined. If $\tilde{M}$ is the local
Maxwellian $M$, with respect to this inner product,
the macroscopic space is spanned by the following five pairwise
orthogonal functions
$$
\left\{
\begin{array}{l}
 \chi_0(\xi) \equiv {\displaystyle\frac1{\sqrt{\rho}}M}, \\[2mm]
 \chi_i(\xi) \equiv {\displaystyle\frac{\xi_i-\v u_i}{\sqrt{R\theta\rho}}M} \ \ {\textrm {for} }\ \  i=1,2,3, \\[2mm]
 \chi_4(\xi) \equiv
 {\displaystyle\frac{1}{\sqrt{6\rho}}(\frac{|\xi-\v u|^2}{R\theta}-3)M},\\
\langle\chi_i,\chi_j\rangle=\delta_{ij}, ~i,j=0,1,2,3,4.
 \end{array}
\right.
$$
Using these functions, we define the macroscopic
projection $P_0$ and microscopic projection $P_1$ as follows:
$$
\left\{
\begin{array}{l}
 P_0h \equiv {\displaystyle\sum_{j=0}^4 \langle h,\chi_j\rangle\chi_j,} \\[2mm]
 P_1h\equiv h-P_0h.
 \end{array}
\right.
$$
The projections $P_0$ and $P_1$ are orthogonal:
$$
P_0P_0=P_0, P_1P_1=P_1, P_0P_1=P_1P_0=0.
$$
A function $h(\xi)$ is called microscopic or non-fluid if
$$
\int h(\xi)\psi_j(\xi)d\xi=0,~j=0,1,2,3,4.
$$
Under this decomposition, the solution $f(t,x,\xi)$ of the
Boltzmann equations satisfies
\begin{equation*}\label{B1.5}
P_0f=M,~P_1f=\v G,
\end{equation*}
and the Boltzmann equation becomes
\begin{equation*}\label{B1.6}
(\v M+\v^2G)_t+\xi_1(M+\v G)_x=2 Q(M,G)+\v Q(G,G),
\end{equation*}
which is equivalent to the following fluid-type system for the
fluid components (see \cite{Liu-Yang-Yu} and
\cite{Liu-Yu} for details):
\begin{equation}\label{B1.7}
\begin{cases}
\di \v\rho_{t}+(\v\rho u_1)_x=0,\\
\di \v(\v\rho u_1)_t+(\v^2\rho
u_1^{2}+{p})_{x}=-\v\int\xi_1^2G_x d\xi,  \\
\di \v(\v\rho u_i)_t+(\v^2\rho u_1u_i)_{x}=-\v\int\xi_1\xi_iG_xd\xi,~ i=2,3,\\
\di \v \bigl[\rho(e+\frac{|\v u|^2}{2})\bigr]_{t}+ \bigl[\v\rho
u_1(e+\frac{|\v u|^2}{2})+\v pu_1\bigr]_x=-\v\int\frac12\xi_1|\xi|^2G_xd\xi,
\end{cases}
\end{equation}
or more precisely,
\begin{equation}\label{B1.8}
\begin{cases}
\v\rho_{t}+(\v\rho u_1)_x=0,\\
\di \v(\v\rho u_1)_t+(\v^2\rho
u_1^{2}+{p})_{x}=\frac43\v(\mu(\theta)
\v u_{1x})_x-\v\int\xi_1^2\Theta_xd\xi,  \\
\di \v(\v\rho u_i)_t+(\v^2\rho u_1u_i)_{x}=\v(\mu(\theta)
\v u_{ix})_x-\v\int\xi_1\xi_i\Theta_xd\xi,~ i=2,3,\\
\v\bigl[\rho(e+\frac{|\v u|^2}{2})\bigr]_{x}+ \bigl[\v\rho
\di u_1(e+\frac{|\v u|^2}{2})+\v pu_1\bigr]_x=\v(\k(\theta)\theta_x)_x\\
\di \quad~~+\frac43\v(\v^2\mu(\theta)u_1u_{1x})_x+\sum_{i=2}^3\v(\v^2\mu(\theta)u_iu_{ix})_x
-\v\int\frac12\xi_1|\xi|^2\Theta_xd\xi,
\end{cases}
\end{equation}
together with an equation for the non-fluid component ${ G}$:
\begin{equation}\label{B1.9}
 \v^2G_t+P_1(\xi_1M_x)+\v P_1(\xi_1G_x)=L_MG+\v Q(G, G),
\end{equation}
where
\begin{equation*}\label{B1.10}
G=L_M^{-1}(P_1(\xi_1 M_x)) +\Theta,
\end{equation*}  and
\begin{equation*}\label{B1.11}
\Theta=L_M^{-1}(\v^2G_t+\v P_1(\xi_1G_x)-\v Q(G, G)).
\end{equation*}
Here $L_M$ is the linearized operator of the collision operator
with respect to the local Maxwellian $M$:
$$
L_M h=Q(M, h)+ Q(h, M),
$$
and the null space $N$ of $L_M$ is spanned by the macroscopic
variables:
$$
\chi_j, ~j=0,1,2,3,4.
$$
Furthermore, there exists a positive constant
$\sigma_0(\rho,u,\theta)>0$ such that for any function $h(\xi)\in
N^\bot$, see \cite{Grad},
$$
\langle h,L_Mh \rangle\le -\sigma_0 \langle\nu(|\xi|)h,h \rangle,
$$
where $\nu(|\xi|)$ is the collision frequency. For the hard sphere
and the hard potential with angular cut-off, the collision
frequency $\nu(|\xi|)$ has the following property
$$
0<\nu_0<\nu(|\xi|)<c(1+|\xi|)^{\beta},
$$
for some positive constants $\nu_0, c$ and $0<\beta\le 1$.

In the above presentation, we normalize the gas constant $R$ to be
$\frac 23$ for simplicity so that $e=\f32R\theta=\t$ and
$p=R\r\t=\frac23\rho\theta$. Notice also that the viscosity coefficient
$\mu(\theta)>0$ and the heat conductivity coefficient
$\k(\theta)>0$ are smooth functions of the temperature
$\theta$. And the following relation holds between these two
functions, \cite{CC,Grad},
\begin{equation}\label{B1.12}
\k(\theta)=\frac{15}{4} R \mu(\theta) =\frac 52 \mu(\theta),
\end{equation}
after taking $R=\frac 23$. It should be pointed out that \eqref{B1.12} is crucially used in the following analysis. In fact, in our analysis, it is required that $$
\inf_{\theta}\k(\theta)>\frac {5}{4}
\sup_{\theta}\mu(\theta)$$ for all $\theta$ under consideration. By
\eqref{B1.12}, it is known that the above condition holds provided that the variation of the
temperature is suitably small.

Now we are in a position to derive the limit equations for $(\r,u,\t)$ in the diffusive limit \eqref{diffusive} formally.
As \cite{BGLY}, we assume that
\begin{equation}\label{pressure}
p^\v=\mbox{const}+O(1)\v^2,
\end{equation}
then, as $\v\rightarrow0$,  $\eqref{B1.8}_1, \eqref{B1.8}_2$ and $\eqref{B1.8}_4$ yields formally that
\begin{equation}\label{B1.8-1}
\begin{cases}
\di p=\mbox{const},\\
\di \rho_{t}+(\rho u_1)_x=0,\\
\di (\rho u_1)_t+(\rho
u_1^{2})_x+P^{\ast}_{x}=\frac43(\mu(\theta)
u_{1x})_x,\\
\di(\rho \t)_{t}+ (\rho
 u_1\t+pu_1)_x=(\k(\theta)\theta_x)_x,\\
\end{cases}
\end{equation}
where $P^\ast$ is unknown function. The equation \eqref{B1.8-1} reveals how  the zero order function $\r,\t$ depend  on the scaled velocity even though the macroscopic velocity tends to zero.

\

With slab symmetry, in the
macroscopic level, it is more convenient to rewrite the system
by using the {\it Lagrangian} coordinates as in the study of conservation laws.
That is, consider the coordinate transformation:
\begin{equation}\nonumber
(x,t)\rightarrow\Big(\int_{(0,0)}^{(x,t)}\r(y,s)dy-(\r u_1)(y,s)ds, s\Big),
\end{equation}
which is  still denoted  as $(x,t)$ without confusion. Denote that $v=\f{1}{\r}$,
the system \eqref{1.1}  and \eqref{B1.7} in the
Lagrangian coordinates become
\begin{equation}\label{B1.13}
\v f_t-\frac{\v u_1}{v}f_x+\frac{\xi_1}{v}f_x=\f{1}{\v}Q(f,f),
\end{equation}
and
\begin{eqnarray}\label{B1.14}
\begin{cases}
\v v_{t}-\v u_{1x}=0,\\
\di \v^2 u_{1t}+p_{x}=-\v\int\xi_1^2G_{x}d\xi,\\
\di \v^2 u_{it}=-\v\int\xi_1\xi_i G_{x}d\xi, ~i=2,3,\\
\di \v \bigl(e+\frac{|\v u|^{2}}{2}\bigr)_{t}+ (\v pu_1)_{x}=-\v\int\frac12\xi_1|\xi|^2G_{x}d\xi,
\end{cases}
\end{eqnarray}
respectively.
Moreover, \eqref{B1.8} and \eqref{B1.9} take the form
\begin{eqnarray}\label{B1.15}
\begin{cases}
\v v_{t}-\v u_{1x}=0,\\[1mm]
\di \v^2 u_{1t}+p_{x}=\frac43\v^2(\frac {\mu(\theta)}vu_{1x})_{x}-\v\int\xi_1^2\Theta_{1x}d\xi,\\[1mm]
\di \v^2 u_{it}=\v^2(\frac{\mu(\theta)}{v}u_{ix})_x-\v\int\xi_1\xi_i\Theta_{1x}d\xi,~i=2,3,\\[1mm]
\di \v \bigl(e+\frac{|\v u|^{2}}{2}\bigr)_{t}+
(\v pu_1)_{x}=\v(\frac{\k(\theta)}{v}\theta_x)_x+\frac43\v^3(\frac{\mu(\theta)}{v}u_1u_{1x})_x\\[1mm]
\di ~~~~+\sum_{i=2}^3\v^3(\frac{\mu(\theta)}{v}
u_iu_{ix})_x-\v\int\frac12\xi_1|\xi|^2\Theta_{1x}d\xi,
\end{cases}
\end{eqnarray}
and
\begin{equation}\label{B1.16}
\v^2 G_t-\frac{\v^2 u_1}{v}G_x+P_1(\frac{\xi_1}{v}M_x)+\v P_1(\frac{\xi_1}{v}G_x)=L_MG+\v Q(G,G),
\end{equation}
with
\begin{equation*}\label{B1.17}
G=L^{-1}_M( P_1(\frac{\xi_1}{v} M_x))+\Theta_1,
\end{equation*}
and
\begin{equation}\label{B1.18}
\Theta_1=L_M^{-1}\Big(\v^2G_t-\frac{\v^2 u_1}vG_x+\frac{\v}{v} P_1(\xi_1G_x)-\v Q(G,G)\Big).
\end{equation}
The limiting equation \eqref{B1.8-1} becomes
\begin{equation}\label{B1.8-2}
\begin{cases}
\di p=\mbox{const},\\
\di v_{t}-u_{1x}=0,\\
\di u_{1t}+P^{\ast}_{x}=\frac43(\f{\mu(\theta)}{v}
u_{1x})_x,\\
\di\t_t+ pu_{1x}=(\f{\k(\theta)}{\t}\theta_x)_x.
\end{cases}
\end{equation}

\subsection{Construction of  profile}

We will construct a background solution to \eqref{B1.8-2} in this subsection. Without loss of generality, set
\begin{align}
p=\f{2\t}{3v}=\f23,
\end{align}
that is
\begin{align}\label{1.22}
v=\t.
\end{align}
Assume the boundary conditions at the far fields given by
\begin{equation}\label{B1.18-1}
\lim_{x\rightarrow \pm\infty}(v,\t)(x,t)=( v_\pm, \t_\pm),~~\mbox{and}~~\f{\t_+}{v_+}=\f{\t_-}{v_-}=1,~~\mbox{with}~\t_-\neq\t_+.
\end{equation}
Note that if $\t_-=\t_+$, then $v=\t=1,u_1=0$ is a trivial solution to \eqref{B1.8-2}, and the diffusive limit of the rescaled Boltzmann equation to the incompressible Navier-Stokes-Fourier system is well studied as mentioned in the introduction.

Noting  \eqref{1.22}, the equation $\eqref{B1.8-2}_4$ is rewritten as
\begin{eqnarray}\label{B1.30i}
\theta_t+\frac23u_{1x}=(\frac{\k(\t)}{\t}\theta_x)_x.
\end{eqnarray}
Substituting $\eqref{B1.8-2}_2$ into \eqref{B1.30i} and noting \eqref{1.22}, we have the following scalar nonlinear diffusion equation
\begin{eqnarray}\label{B1.30-1}
\theta_t=(a(\t)\theta_x)_x, ~a(\t)=\frac{3\k(\t)}{5\t},~~\mbox{with}~~\lim_{x\rightarrow\pm\infty}\t(x,t)=\t_{\pm}.
\end{eqnarray}
From \cite{Atkinson-Peletier} and \cite{Duyn-Peletier}, it is known that
the nonlinear diffusion equation \eqref{B1.30-1} admits a self-similar solution $\hat\t(\eta)$ with $\eta=\f{x}{\sqrt{1+t}}$
satisfying the  boundary conditions $
\hat\t(\pm\i,t)=\t_\pm$. Furthermore, $\hat{\t}(\eta)$ is a monotonic function. Let $\delta=|\t_+-\t_-|$, then
$\hat\t(t,x)$ has the property that
\begin{equation}\label{N-1.4}
\hat\t_x(t,x)=\f{O(1)\delta}{\sqrt{1+t}}e^{-\f{x^2}{4a(\t_\pm)(1+t)}},\qquad {\rm as}~~x\rightarrow\pm\i.
\end{equation}
Define
\begin{equation}\label{tilde}
(\tilde{v},\tilde{u}_1, \tilde{\t})\doteq (\hat\t, a(\hat\t)\hat\t_x ,\hat\t)(x,t),
\end{equation}
then it is easy to check that $(\tilde{v},\tilde{u}_1, \tilde{\t})$ satisfying \eqref{B1.8-2} as
\begin{equation}\label{1.23}
\begin{cases}
\di \tilde p=\f{2\tilde\t}{3\tilde v}=\f23,\\
\di \tilde{v}_{t}-\tilde{u}_{1x}=0,\\
\di \tilde{u}_{1t}+P^{\ast}_{x}=\frac43(\f{\mu(\tilde\theta)}{\tilde v}
\tilde{u}_{1x})_x,\\
\di\tilde\t_t+ \tilde{p}\tilde{u}_{1x}=(\f{\k(\tilde\theta)}{\tilde\t}\tilde\theta_x)_x,
\end{cases}
\end{equation}
where $P^\ast=-a(\tilde\t)\tilde\t_t+\f{4\mu(\tilde\t)}{\tilde\t}(a(\tilde\t)\tilde\t_x)_x$.

\begin{remark}\label{rmk2.1}
By \eqref{tilde} and \eqref{1.23}, we actually
construct  a diffusion wave to the  limit system. On the other hand,
if $\t_-<\t_+$, then $\tilde{u}_1=a(\hat\t)\hat\t_x>0$, that is, the variation of temperature along the x-axis induces a nontrivial scaled flow along the same direction, see  Figure 1. The case $\t_->\t_+$ is similar, see Figure 2.
\end{remark}

\begin{center}
	\includegraphics[scale=0.8]{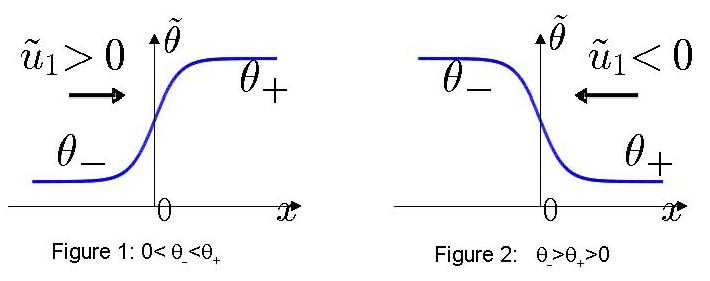}
\end{center}

\begin{remark}\label{rmk2.2}
The construction of the profile  $(\tilde{v},\tilde{u}_1,
\tilde{\t})$ is motivated by  the viscous contact wave of
compressible Navier-Stokes equations,  see
\cite{Huang-Matsumura-Xin}, \cite{Huang-Xin-Yang} and
\cite{Huang-Yang}. The viscous contact wave is used to approximate
the contact discontinuity for compressible Euler equation and its
pressure keeps constant.
\end{remark}

In order to justify the hydrodynamic limit of the rescaled Boltzmann equation to the limit system \eqref{1.23},  if we  use the profile
$(\tilde{v},\tilde{u}_1, \tilde{\t})$, then some non-integrable
error terms with respect to time coming from the non-fluid
component for the system about perturbation.
Therefore, one needs to construct another profile $(\bar{v},\v \bar u, \bar\t)$
for the rescaled Boltzmann equation, based on  $(\tilde{v},\tilde{u}_1, \tilde{\t})$.
For this, we require
that the approximate pressure $p$ satisfies
\begin{equation}\label{B1.18-2}
\bar p=\f{2\bar\t}{3\bar v}=\f23+O(1)\v^2\doteq p_++O(1)\v^2.
\end{equation}
Motivating by \cite{Huang-Yang},  we
first notice that the main part of the non-fluid component in
the solution $G$ and part of $\Theta_1$ defined in \eqref{B1.18}, are
given by
\begin{equation*}\label{B1.19}
w=\frac1vL_M^{-1}(P_1(\xi_1M_x))=\frac1{Rv\theta}L^{-1}_M\{P_1[\xi_1(\frac{|\xi-\v u|^2}{2\theta}
\theta_x+\xi\cdot \v u_{x})M]\},
\end{equation*}
 and
\begin{equation*}\label{B1.20}
\hat{\Theta}_1=L_M^{-1}(\frac{\v}{v}P_1(\xi_1w_x)-\v Q(w,w)),
\end{equation*}
respectively. To distinguish the leading term coming from the non-fluid component, we rewrite the Boltzmann equation \eqref{B1.15} as
\begin{eqnarray}\label{B1.21}
\begin{cases}
{ \v v_{t}-\v u_{1x}=0,}\\[2mm]
\di { \v^2 u_{1t}+p_{x}=\frac43\v^2 (\frac {\mu(\theta)}vu_{1x})_{x}-\sum_{j=1}^2\v \int\xi_1^2\Theta^j_{1x}d\xi,}\\[2mm]
\di {
\v^2 u_{it}=\v^2 (\frac{\mu(\theta)}{v}u_{ix})_x-\sum_{j=1}^2\v \int\xi_1\xi_i\Theta^j_{1x}d\xi,~
i=2,3,}\\[2mm]
\di {\v \bigl(e+\frac{|\v u|^{2}}{2}\bigr)_{t}+
(\v pu_1)_{x}=\v (\frac{\k(\theta)}{v}\theta_x)_x-\sum_{j=1}^2\v \int\frac12\xi_1|\xi|^2\Theta^j_{1x}d\xi
+H_{x},}
\end{cases}
\end{eqnarray}
with
\begin{eqnarray}\label{B1.22}
&&\v^2\tilde{G}_t-L_M\tilde{G}=-\frac1{Rv\theta}P_1[\xi_1
(\frac{|\xi-\v u|^2}{2\theta}(\theta-\bar{\theta})_x+\xi\cdot(\v u-{\v\bar
u})_x)M]\nonumber\\
&&~~~~~~~~~~~~~~~~~~~~~~~~+\frac{\v^2 u_1}{v}G_x-\frac{\v}{v}P_1(\xi_1G_x)+\v Q(G,G)-\v^2\bar{G}_t,
\end{eqnarray}
where
\begin{eqnarray}\label{B1.23}
\begin{cases}
\bar G=\frac1{Rv\theta}L^{-1}_M\{P_1[\xi_1(\frac{|\xi-\v u|^2}{2\theta}{\bar
\theta}_x+\xi\cdot{\v\bar u}_{x})M]\},\qquad \tilde{G}=G-{\bar G},\\[2mm]
H=\frac{4\v^3}3\frac{\mu(\theta)}{v}u_1u_{1x}+\sum_{i=2}^3\v^3\frac{\mu(\theta)}{v}
u_iu_{ix},\\[2mm]
\Theta_1^1=L_M^{-1}\Big(\frac{\v}vP_1(\xi_1{\bar G}_x)-\v Q({\bar G},{\bar G})\Big),\\[2mm]
\Theta_1^2=
L_M^{-1}\Big(\v^2G_t-\frac{\v^2u_1}vG_x+\frac{\v}{v}P_1(\xi_1\tilde{G}_x)
-\v Q(\tilde{G},\tilde{G})-2\v Q(\bar{G},\tilde{G})\Big),
\end{cases}
\end{eqnarray}
satisfying
$$
\sum_{j=1}^2\Theta^j_1=\Theta_1=
L_M^{-1}(\v^2G_t-\frac{\v^2 u_1}vG_x+\frac{\v}{v}P_1(\xi_1G_x)-\v Q(G,G)).
$$
Here,  the function $(\bar{v},\v\bar{u},\bar{\theta})(x,t)$ is the profile
to be constructed.

\

Since the velocity $\v u$ decays faster than $(v,\theta)$ in time,
the leading terms in the energy equation $\eqref{B1.21}_4$ are
\begin{eqnarray}\label{B1.28}
\v\theta_t+\v pu_{1x}=\v(\frac{\k(\theta)}v\theta_{x})_x-\v\int\frac12\xi_1|\xi|^2\Theta^1_{1x}d\xi.
\end{eqnarray}
By the definition of $\Theta_1^1$, it holds that

\begin{eqnarray}\label{B1.29}
\begin{cases}
-\v\int\frac12\xi_1|\xi|^2\Theta_{1}^1d\xi=\v^2N_1+\v^3F_1,\\
N_1=f_{11}\theta_x{\bar\theta}_x+f_{12}v_x\bar{\theta}_x
+f_{13}\bar{\theta}_x^2+f_{14}\bar{\theta}_{xx},\\
|F_1|=O(1)[(|v_x|+|\theta_x|+|\bar{\theta}_x|+\v|u_x|+\v|\bar{u}_x|)|\bar{u}_x|
+|u_x\bar\t_x|+|{\bar u}_{xx}|],
\end{cases}
\end{eqnarray}
where the coefficients $f_{1j},j=1,2,3,4$ are smooth functions of $(v,\v u,\theta)$.
By \eqref{B1.18-2}, it is expected
that the profile $(\bar{v},\v\bar{u},\bar{\theta})$ for the
Boltzmann equation satisfies $\bar\t\cong\bar{v}$. Thus,
by choosing only the leading term in \eqref{B1.28}, one obtains that
\begin{eqnarray}\label{B1.30}
\v\theta_t=\v(a(\t)\theta_x)_x+\frac{3\v^2}{5}N_{1x},
\end{eqnarray}
where $a(\t)$ is given in \eqref{B1.30-1}. Thus the leading part of \eqref{B1.30} is the nonlinear diffusion equation \eqref{B1.30-1} and
an explicit solution $\hat\t(\f{x}{\sqrt{1+t}})$ is given with the boundary conditions $\hat\t(\pm\i,t)=\t_\pm$.

To include more microscopic
effect, let the  profile $\bar{\t}\approx
\hat\t(\f{x}{\sqrt{1+t}})+\v\theta^{nf}(x,t)$, where $\theta^{nf}(x,t)$ represents the
part of the nonlinear diffusion wave coming from the non-fluid
component. Moreover, the
term $\theta^{nf}(x,t)$ in the form of
$\frac1{\sqrt{1+t}}D_1(\frac{x}{\sqrt{1+t}})$ is from $N_1$ in
\eqref{B1.30}. Note that $\theta^{nf}(x,t)$ decays faster than
$\hat\theta(x,t)$ so that it can be viewed as a perturbation around
 profile $\hat\theta(x,t)$. To construct
$\theta^{nf}(x,t)$, we linearize the equation \eqref{B1.30} around
$\hat\theta(x,t)$ and keep only the linear terms. This leads to a
linear equation for $\theta^{nf}(x,t)$ from \eqref{B1.30}
\begin{eqnarray}\label{B1.31}
\theta_t^{nf}=(a(\hat\theta)\theta_x^{nf})_x+(a'(\hat\theta)\hat\theta_x\theta^{nf})_x+\frac35\hat{N}_{1x},
\end{eqnarray}
where
$\hat{N}_1=(\hat{f}_{11}+\hat{f}_{12}+\hat{f}_{13})(\hat\theta_x)^2+
\hat{f}_{14}\hat\theta_{xx}$ with
$\hat{f}_{1j}=f_{1j}(\tilde{v},0,\hat\t),j=1,2,3,4$.
Let
$$
g_{1}(x,t)=\int_{-\infty}^x \theta^{nf}(x,t)dx,
$$
then integrating \eqref{B1.31} with respect to $x$ yields that
\begin{eqnarray}\label{B1.32}
g_{1t}=a(\hat\t)g_{1xx}+a'(\hat\t)\hat\t_x g_{1x}+\frac35{\hat
N}_1.
\end{eqnarray}
Note
that ${\hat N}_1$ takes the form of
$\frac1{1+t}D_2(\frac{x}{\sqrt{1+t}})$ and satisfies
$$
|{\hat N}_1|=O(1)\d(1+t)^{-1}e^{-\f{x^2}{4a(\t_\pm)
(1+t)}}, ~~~{\rm as}~~x\rightarrow\pm\i.
$$
We can check that there exists a
self-similar solution $g_{1}(\eta), \eta=x/\sqrt{1+t}$ for \eqref{B1.32}
with the boundary condition
$g_{1}(-\infty,t)=0,~g_{1}(+\infty,t)=\delta_1 $. Here $\delta_1$ satisfies $0<\delta_1<
\delta$. Note that even though the
function $g_1(x,t)$ depends on the constant $\delta_1$,
$\theta^{nf}(x,t)=g_{1x}(x,t)\to 0$ as
$x\rightarrow\pm\i.$ That is, the choice of the constant $\d_{1}$
has no influence on the ansantz as long as $|\d_{1}|<\d$.
From now on, we fix $\d_{1}$ so that the function $g_1(x,t)$ is
uniquely determined and its derivative $g_{1x}=\theta^{nf}$
has the property
\begin{equation*}
|\theta^{nf}|=|g_{1x}|=O(\delta)(1+t)^{-\f12}e^{-\f{x^2}{4a(\t_\pm)(1+t)}},~~~{\rm
as}~~x\rightarrow\pm\i.\label{T-nf}
\end{equation*}
Now we follow the same procedure to construct the second and third
components of the velocity profile denoted by
$\v\bar{u}_i,~i=2,3$. That is, the leading part of the equation
for $\v u_i$ coming from \eqref{B1.21} is
\begin{eqnarray}\label{B1.33}
\v^2 u_{it}=\v^2(\frac{\mu(\theta)}{\theta}u_{ix})_x-\v\int\xi_1\xi_i\Theta_{1x}^1d\xi.
\end{eqnarray}
For $i=2,3$, one gets
\begin{eqnarray}\label{B1.34}
\begin{cases}
-\v\int\xi_1\xi_i\Theta_{1}^1d\xi=\v^2N_{i}+\v^3F_{i},\\
N_{i}=f_{i1}\theta_x{\bar\theta}_x+f_{i2}v_x\bar{\theta}_x+f_{i3}\bar{\theta}_x^2
+f_{i4}\bar{\theta}_{xx},\\
|F_{i}|=O(1)((|v_x|+|\theta_x|+|\bar{\theta}_x|+\v|u_x|+\v|\bar{u}_x|)|\bar{u}_x|
+|u_x||\bar{\theta}_x|+|{\bar u}_{xx}|),
\end{cases}
\end{eqnarray}
with smooth functions $f_{ij},i=2,3,~j=1,2,3,4$. Notice that the
symbols $N_i$ and $F_i$, $i=2,3,$ used here are for the convenience
of notations.

From \eqref{B1.33} and \eqref{B1.34}, we expect that the  profile ${\bar u}_i(x,t)$
takes the form of $\f{1}{\sqrt{1+t}}h_i(\f{x}{\sqrt{1+t}})$ and
satisfies the following linear equation
\begin{eqnarray}\label{B1.35}
\v^2{\bar
u}_{it}=\v^2(\frac{\mu(\hat\t)}{\hat\t}\bar{u}_{ix})_x+\v^2\hat{N}_{ix},~i=2,3,
\end{eqnarray}
where
$\hat{N}_i=(\hat{f}_{i1}+\hat{f}_{i2}+\hat{f}_{i3})(\hat\theta_x)^2+
\hat{f}_{i4}\hat\theta_{xx}$,
$\hat{f}_{ij}=f_{ij}(\tilde{v},0,\hat\t),i=2,3,j=1,2,3,4$.

Denote
$$
g_{i}(x,t)=\int_{-\infty}^x \bar{u}_i(x,t)dx,
$$
then integrating \eqref{B1.35} with respect to $x$, one has
\begin{eqnarray}\label{B1.36}
g_{it}=\frac{\mu(\hat\t)}{\hat\t}g_{ixx}+{\hat N}_i.
\end{eqnarray}
For given
$\hat\t$, we can check that there exists a self-similar solution
 $g_i(\eta)$ with $\eta=\f{x}{\sqrt{1+t}}$ with the boundary conditions
$g_{i}(-\infty,t)=0$, $g_{i}(+\infty,t)=\delta_i$, where $\delta_i$ satisfies
$0<\delta_i<\delta$. As we explained before, the choice of the constant
$\delta_i$ is not essential. From \eqref{N-1.4}, we fix $\d_{i}$ so that the function
$g_i(x,t)$ is uniquely determined and the derivative
$g_{ix}={\bar u}_i~ (i=2,3)$ has
the following property
\begin{align*}\label{B1.37}
|\v\bar{u}_{i}|=|\v g_{ix}|=O(1)\delta\v(1+t)^{-\frac12}e^{-\frac{x^2}{4b(\theta_\pm)(1+t)}},
\quad {\rm as}~ x\to \pm\infty,
\end{align*}
where
$b(\theta_\pm)=\max\{a(\theta_\pm),\frac{\mu(\theta_\pm)}{\theta_\pm}\}$.

In summary, one can define the  profile $({\bar
v},\v\bar{u},\bar{\theta})$ for the Boltzmann equation as follows.
To satisfy the conservation of mass, one needs
$$
\v\bar{v}_t-\v\bar{u}_{1x}=0.
$$
By letting $\bar{v}=\hat\t+\v\t^{nf}$, one gets
\begin{eqnarray}\label{B1.38}
\v{\bar
u_1}=\v[a(\hat\theta)\hat\theta_x+\v a(\hat\theta)\theta_x^{nf}
+\v a'(\hat\theta)\hat\theta_x\theta^{nf}]+\frac{3\v^2}{5}\hat{N}_1.
\end{eqnarray}
However, by plugging \eqref{B1.38} into the momentum equation of \eqref{B1.21},
we have a non-conservative term containing $\v^2\hat{N}_{1t}$. To
avoid this, one defines
\begin{equation*}\label{B1.39}
\v{\bar
u_1}=\v[a(\hat\theta)\hat\theta_x+\v a(\hat\theta)\theta_x^{nf}
+\v a'(\hat\theta)\hat\theta_x\theta^{nf}].
\end{equation*}
Similarly, to avoid the non-conservative term $(|\bar{u}|^2)_t$ in
the energy equation, set
$$
\tilde{\theta}=\theta^{ns}+\v\theta^{nf}-\frac12|\v\bar{u}|^2.
$$

Therefore, the profile $({\bar v},\v\bar{u},\bar{\theta})$ is finally defined as:
\begin{eqnarray}\label{B1.40}
\begin{cases}
{\displaystyle
{\bar v}=\hat\t+\v\t^{nf},}\\
{\displaystyle
\v{\bar
u_1}=\v[a(\hat\theta)\hat\theta_x+\v a(\hat\theta)\theta_x^{nf}
+\v a'(\hat\theta)\hat\theta_x\theta^{nf}],}\\
{\displaystyle \v\bar{u}_i=\v g_{ix},~i=2,3,}\\
{\displaystyle
 {\bar
\theta}=\hat\t+\v\theta^{nf}-\frac12|\v\bar{u}|^2,}
\end{cases}
\end{eqnarray}
where $\hat\t$ is given by \eqref{B1.30-1}, $\theta^{nf}$ by \eqref{B1.31} and
$g_i$, $i=2,3$ by \eqref{B1.36}. Then a direct but tedious computation shows that
\begin{equation}\label{B1.41}
\left\{\begin{array}{llll}
{\displaystyle \v\bar{v}_{t}-\v\bar{u}_{1x}=\frac{3\v^2}{5}\hat{N}_{1x},}\\[2mm]
{\displaystyle \v^2\bar{u}_{1t}+\bar{p}_{x}=\frac{4\v^2}{3}(\frac{\mu(\bar{\theta})}{\bar{v}}\bar{u}_{1x})_{x}+\bar{R}_{1x},}\\[3mm]
{\displaystyle
\v^2\bar{u}_{it}=\v^2(\frac{\mu(\bar{\theta)}}{\bar{v}}\bar{u}_{ix})_x+\v^2\bar{N}_{ix}+\bar{R}_{ix},
i=2,3,}\\
{\displaystyle \v\bigl(\bar{e}+\frac{|\v\bar{u}|^{2}}{2}\bigr)_{t}+
(\v\bar{p}\bar{u}_1)_{x}=\v(\frac{\k(\bar{\theta})}{\bar{v}}\bar{\theta}_x)_x+\bar{H}_{x}
+\v^2\bar{N}_{1x}-\frac{2\v^2}{5}\hat{N}_{1x}+\bar{R}_{4x},}
\end{array}
\right.
\end{equation}
where
\begin{equation}\label{B1.42}
\begin{array}{ll}
\bar R_{1}&={\displaystyle
\v^2[a(\hat\t)\hat\t_{t}+(a(\hat\t)\theta^{nf})_t]+
\bar{p}-p_+-\frac43\v(\frac{\mu(\bar{\theta})}{\bar{v}}\v\bar{u}_{1x})}\\
&={\displaystyle
O(1)\delta\v^2(1+t)^{-1}e^{-\frac{x^2}{4c(\theta_\pm)(1+t)}},\quad
{\rm as}~x\to \pm\infty,}
\end{array}
\end{equation}
\begin{equation}\label{B1.43}
\begin{array}{ll}
\bar{R}_{i}&={\displaystyle
\v[\frac{\mu(\hat\t)}{\hat\t}-\frac{\mu(\bar{\theta})}{\bar{v}}]\v\bar{u}_{ix}+\v^2(\hat{N}_{i}-\bar{N}_i)}\\
&={\displaystyle
O(1)\delta\v^3(1+t)^{-3/2}e^{-\frac{x^2}{4c(\theta_\pm)(1+t)}},\quad
{\rm as}~x\to \pm\infty,~i=2,3,}
\end{array}
\end{equation}
\begin{equation}\label{B1.44-1}
\begin{array}{lll}
\bar{R}_{4}&=&{\displaystyle
\Big[\frac53\v(a(\hat\t)\hat\t_x+a(\hat\t)\theta_x^{nf}+a'(\hat\t)
\hat\t_x\theta^{nf})
-\v\frac{\k(\bar{\theta})}{\bar{v}}\bar{\theta}_x\Big]}\\
&&{\displaystyle +(\bar{p}-p_+)\v\bar{u}_{1}+\v^2(\hat{N}_{1}-\bar{N}_1)-\bar{H}}\\
&=&{\displaystyle
O(\delta)\v^3(1+t)^{-3/2}e^{-\frac{x^2}{4c(\theta_\pm)(1+t)}},\quad
{\rm as}~x\to \pm\infty,}
\end{array}
\end{equation}
\begin{equation}\label{B1.44}
\hat{N}_{i}=
O(1)\delta(1+t)^{-1}e^{-\frac{x^2}{4a(\theta_\pm)(1+t)}},\quad
{\rm as}~x\to \pm\infty,~i=1,2,3,
\end{equation}
with
$c(\theta_{\pm})=\max\{a(\theta_{\pm}),\frac12b(\theta_{\pm})\},$
$\bar{N}_{i},i=1,2,3$, and $\bar{H}$ are the corresponding
functions defined in \eqref{B1.23}, \eqref{B1.29} and \eqref{B1.34} by substituting the
variable $(v,\v u,\theta)$ by the  profile
$(\bar{v},\v\bar{u},\bar{\theta})$. Note
that the decay rate of $\bar R_i, i=2,3,4$ is of order
$\v^3(1+t)^{-3/2}$.  Furthermore, even though the decay
rate of $\bar R_1$ is still $\v^2(1+t)^{-1}$, it is sufficient to obtain the
desired a priori estimates through some subtle analysis coming from
the intrinsic dissipation mechanism in  the momentum
equations as shown in the following.

Define
\begin{equation*}\label{B1.45}
\bar{M}=\frac{\bar{v}^{-1}}{\sqrt{ (2 \pi R \bar\t)^3}}
\exp\left(-\frac{|\xi -\v\bar{u}|^2}{2R\bar\t}\right),~~\bar{G}_0=L_{\bar{M}}^{-1}\left(\f{1}{\bar v}\bar{P}_1(\xi_1\bar{M}_x)\right),
\end{equation*}
and
\begin{equation*}\label{B1.46}
\bar{f}=\bar{M}+\v\bar{G}_0.
\end{equation*}
Then it follows  from \eqref{B1.41} that
\begin{equation}\label{B1.47}
\v\bar{f}_t-\f{\v\bar{u}_1}{\bar v}\bar{f}_x+\f1{\bar v}\xi_1\bar{f}_x
=L_{\bar M}\bar{G}_0+\v Q(\bar{G}_0,\bar{G}_0)+\bar{R}_{\bar f},
\end{equation}
where
\begin{equation*}\label{B1.48}
\bar{R}_{\bar f}=\v^2\hat{B}_2(x,t,\xi)\bar{M}+\v^2\bar{G}_{0t}-\v\f{\v\bar u_1}{\bar v}\bar{G}_{0x}
+\v\bar{P}_1\Big(\f{\v}{\bar v}\bar{G}_{0x}\Big)-\v Q(\bar{G}_0, \bar{G}_0),
\end{equation*}
and $|\hat{B}_2(x,t,\xi)| =O(1)\d(1+t)^{-\f32}e^{-\frac{x^2}{4c(\theta_\pm)(1+t)}}|\xi|^3$,  as $x\to \pm\infty$.

\begin{remark}
From the definition  of $(\tilde{v}, \tilde{u}_1, \tilde\t)$ in \eqref{tilde} and the definition of $(\bar{v}, \bar{u}_1, \bar\t)$ in \eqref{B1.40}, it holds that
\begin{equation}\label{2.55}
|(\bar{v}-\tilde{v}, \bar{u}_1-\tilde{u}_1, \bar\t-\tilde\t)(x,t)|=O(1) \d\v(1+t)^{-\f12}e^{-\f{x^2}{4c(\t_\pm)(1+t)}}, 
\end{equation}
that implies that the ansantz $(\bar v, \bar u_1, \bar\t)$
well approximates   $(\tilde{v}, \tilde{u}_1,
\tilde\t)$  when $\v$ is small.
\end{remark}

\subsection{Main result}

Now we consider the system \eqref{B1.15}-\eqref{B1.16} with the initial data
\begin{equation}\label{B1.49}
(v,u,\t)|_{t=0}=(\bar{v}, \bar{u}, \bar\t)(x,0),~~~~ G(x,t)|_{t=0}=\bar{G}(x,0).
\end{equation}
Then the main result in this paper can be stated as follows.
\begin{theorem}\label{thm2.1}
Let $(\bar{v}, \bar{u}, \bar\t)(x,t)$ be the  profile defined in
\eqref{B1.40} with strength $\d=|\t_+-\t_-|$. Then there exist small
positive constants $\d_0$ and $\v_0$ and a global Maxwellian
$M_\ast=M_{[v_\ast, u_{\ast},\t_{\ast}]}$, such that when $\d\leq
\d_0$  and $\v\leq \v_0$, the Cauchy problem  \eqref{B1.15}-
\eqref{B1.16} with the initial data \eqref{B1.49} has a unique
global solution $(v,u,\t, G)$ satisfying, for any sufficiently small
but fixed positive constant $\vartheta>0$,
\begin{equation}\label{B1.50-1}
\begin{cases}
\|(v-\bar{v},\v u-\v \bar{u},\t-\bar{\t})(t)\|^2_{L^2_x}\leq C\sqrt\d\v^3(1+t)^{-1+C_0\sqrt\d},\\
\|(v-\bar{v},\v u-\v\bar{u},\t-\bar{\t})_x(t)\|^2_{L^2_x}\leq C\sqrt\d\v^2(1+t)^{-\f32+\vartheta+C_0\sqrt\d},\\
\|f_{xx}(t)\|^2_{L^2_x(L^2_{\xi}(\f1{\sqrt{M_\ast}}))}
+\|(v-\bar{v},\v u-\v\bar{u},\t-\bar{\t})_{xx}(t)\|^2_{L^2_x}\leq C\sqrt\d(1+t)^{-\f32+\vartheta+C_0\sqrt\d},\\
\|(G-\bar G)(t)\|^2_{L^2_x(L^2_{\xi}(\f1{\sqrt{M_\ast}}))}\leq C\sqrt\d(1+t)^{-\f12},\\
\|(G-\bar G)_x(t)\|^2_{L^2_x(L^2_{\xi}(\f1{\sqrt{M_\ast}}))}\leq C\sqrt\d(1+t)^{-\f32+\vartheta+C_0\sqrt\d},
\end{cases}
\end{equation}
that implies that
\begin{equation}\label{B1.50}
\begin{cases}
\|(v-\bar{v},\v u-\v\bar{u},\t-\bar{\t})(t)\|_{L^\infty_x}\leq C\d^\f14\v^{\f54}(1+t)^{-\f58+\f34\vartheta},\\
\|(v-\bar{v},\v u-\v\bar{u},\t-\bar{\t})_x(t)\|_{L^\infty_x}\leq C\d^\f14\v^{\f12}(1+t)^{-\f34+\vartheta},
\end{cases}
\end{equation}
where $C$ and $C_0$ are positive constants independent of $\v$ and $\d$.
\end{theorem}

\

The following result justifies the hydrodynamic limit of the rescaled Boltzmann equation \eqref{1.1} to the diffusion wave $(\tilde{v}, \tilde{u}_1, \tilde{\t})$ global in time.
\begin{corollary}\label{corB1.3}
Under the conditions of Theorem \ref{thm2.1}, from \eqref{2.55} and \eqref{B1.50}, it holds that
\begin{eqnarray}\label{1.55}
\begin{cases}
|(v-\tilde{v}, \t-\tilde\t)(x,t)|\leq C\v(1+t)^{-\f12}\rightarrow0,\\
|(u_1-\tilde{u}_1)(x,t)|\leq C\v^{\f14}(1+t)^{-\f12},\rightarrow0,
\end{cases}
\mbox{as}~\v\rightarrow0,
\end{eqnarray}
that is,  the fluid part $(v, u_1, \t)$ of  the solution   of the rescaled Boltzmann equation
\eqref{1.1} converges to the diffusion wave solution $(\tilde{v}, \tilde{u}_1,
\tilde\t)$ of \eqref{1.23}  in the sense of \eqref{1.55} as $\v\rightarrow0$, which reveals that $v$ and $\t$ are diffusive.
\end{corollary}

\begin{remark}
The above Corollary shows that if the zero order function  in \eqref{expansion} is not a global Maxwellian, then one has to consider the effect of diffusive wave in the diffusive limit of rescaled Boltzmann equation \eqref{1.1}. 
\end{remark}

Since the scaled velocity $\tilde{u}_1$ is actually induced by the variation of temperature $\tilde{\t}$, i.e., $\tilde{u}_1=a(\tilde{\t})\tilde{\t}_x$.  The following result shows that the scaled velocity $u_1$ is also induced by the variation of temperature $\t$ in some sense when $\v$ is small. From the definition of $\hat\t(\eta)$ with $\eta=\f{x}{\sqrt{1+t}}$ in \eqref{B1.30-1} and \eqref{N-1.4}, it can be seen that $\hat\t$ is monotonic.
To be definite and without loss of generality, let us assume that $\t_-<\t_+$, that is, $\hat\t$ is monotonically increasing.  Then there exists a positive constant $\eta_0>0$ such that
\begin{eqnarray}\label{B1.53}
\hat\t'(\eta)>c_{\eta_0}\d, \qquad \mbox{for}~~ |\eta|\leq \eta_0,
\end{eqnarray}
where $c_{\eta_0}$ depends on $\eta_0$ and $c_{\eta_0}\rightarrow0$ as $\eta_0\rightarrow+\infty$.

\begin{corollary}\label{corB1.2}
Under the conditions of Theorem \ref{thm2.1} and $\t_-<\t_+$, for any fixed $\eta_0>0$, there exists a small positive constant $\v_1=\v_1(\eta_0)\leq \v_0$, such that if $\v\leq \v_1$, then it follows from \eqref{B1.53} and \eqref{B1.50} that
\begin{equation}\label{B1.55}
\begin{cases}
0<\f{c_{\eta_0}\d}{C_{1}\sqrt{1+t}}<\f{1}{C_{1}}\hat\t_x\leq u_1(x,t)\leq C_{1}\hat\t_x,\\[3mm]
~~~~~~~~0<\f12 \hat\t_x\leq \t_x(x,t)\leq \f32\hat\t_x,
\end{cases}
\mbox{\rm for}~ |x|\leq \eta_0(1+t)^{\f12},~t\geq0,
\end{equation}
that is
\begin{equation}\label{B1.55-1}
\f{2}{3C_1}\t_x(x,t)\leq u_{1}(x,t)\leq 2C_1\t_x(x,t),~\mbox{\rm for}~ |x|\leq \eta_0(1+t)^{\f12},~t\geq 0,
\end{equation}
where $C_{1}$ is a suitably large positive constant  depending only on $\t_\pm$.
In particular, \eqref{B1.55-1} implies that variation of the  temperature induces a non-trivial flow of higher order  in the following parabolic  region
\begin{equation*}
\Big\{(x,t)~:~ |x|\leq \eta_0(1+t)^{\f12},~t\geq 0~\Big\}.
\end{equation*}
\end{corollary}

\

\section{Stability Analysis}

In this section, we will investigate the stability of the profile constructed in \eqref{B1.41}
for the Boltzmann equation \eqref{1.1}. This section is organized as follows: in Section 3.1,
the fluid type system \eqref{B1.8} is reformulated in terms of the integrated variables; Section
3.2 is devoted to the lower order estimate, while Section 3.3 is for the derivative estimate.

\subsection{Reformulated system}
We now reformulate the system by introducing a scaling for the independent variables. Set
\begin{equation}\label{B4.2}
y=\f{x}{\v},~~~\tau=\f{t}{\v^2}.
\end{equation}
In the following, we will also use the notations $(v,u,\t)(\tau,y)$ and $(\bar v,\bar u,\bar\t)(\tau,y)$, etc., in the scaled independent variables. Set the perturbation around the  profile $(\bar v,\bar u,\bar\t)(\tau,y)$ by
\begin{equation*}\label{B4.3}
\phi=v-\bar{v}, \psi=\v u-\v\bar{u},
\zeta=\theta-\bar{\theta},
\end{equation*}
and
\begin{eqnarray*}\label{B4.4}
(\Phi,\Psi,\bar{W})(y,\tau)=\int_{-\infty}^{y}\Big(\phi,~\psi,~(\t+\f{|\v u|^2}{2})-(\bar \t+\f{|\v\bar u|^2}{2})\Big)(z,\tau)dz.
\end{eqnarray*}
Then we have
$(\phi,\psi)=(\Phi,\Psi)_y$ and
$\zeta+\frac12|\Psi_y|^2+\sum_{i=1}^3\v\bar{u}_i\Psi_{iy}=\bar{W}_y$.

Subtracting \eqref{B1.41} from the equation \eqref{B1.21} and integrating the
reduced system yield
\begin{equation}\label{B4.6}
\left\{
\begin{array}{llll}
\di \Phi_{\tau}-\Psi_{1y}=-\frac2{5p_+}\v^2\hat{N}_{1},\\[0.2cm]
\di \Psi_{1\tau}+p-\bar{p}=\frac{4\v}3\Big(\frac
{\mu(\theta)}vu_{1y}
-\frac43\frac{\mu(\bar{\theta})}{\bar{v}}\bar{u}_{1y}\Big)-\v\sum_{j=1}^2\int\xi_1^2\Theta_{1}^jd\xi-\bar{R}_1,\\[0.2cm]
\di \Psi_{i\tau}=\v\Big(\frac{\mu(\theta)}{v}u_{iy}-\frac{\mu(\bar{\theta})}{\bar{v}}\bar{u}_{iy}\Big)
+\v^2(N_{i}-\bar{N}_{i})+\v^3F_{i}-\v\int\xi_1\xi_i\Theta_{1}^2d\xi-\bar{R}_i,i=2,3,\\[0.2cm]
\di \bar{W}_{\tau}+\v pu_1-\v\bar{p}\bar{u}_1=\Big(\frac{\k(\theta)}{v}\theta_y
 -\frac{\k(\bar{\theta})}{\bar{v}}\bar{\theta}_y\Big)+(H-\bar{H})+\v^2(N_{1}-\bar{N}_{1})+\v^3F_{1}
 \\[0.2cm]
\di \quad~~~~~~~~~~~~~~~~~~~~~~~~-\v\int\frac12\xi_1|\xi|^2\Theta_1^{2}d\xi-\bar{R}_4+\frac25\v^2\hat{N}_1.\\
\end{array}
\right.
\end{equation}
Since the variable $\bar{W}$ is the anti-derivative of the total
energy, not the temperature, it is more convenient to introduce
another variable
\begin{equation*}\label{B4.7}
W=\bar{W}-\v\bar{u}_1\Psi_1.
\end{equation*}
It follows that
\begin{equation*}\label{B4.8}
\zeta=W_y-Y,~\mbox{with}~Y=\frac12|\Psi_y|^2-\v\bar{u}_{1y}\Psi_1+\v\bar{u}_2\Psi_{2y}
+\v\bar{u}_3\Psi_{3y}.\end{equation*}
Using the new variable $W$ and linearizing the left hand side of
the system \eqref{B4.6} by using the formula of $H$ in \eqref{B1.23} give that
\begin{eqnarray}\label{B4.9}
\left\{
\begin{array}{llll}
\di  \Phi_{\tau}-\Psi_{1y}=-\frac3{5}\v^2\hat{N}_{1},\\[0.2cm]
\di \Psi_{1\tau}-\frac{p_+}{\bar{v}}\Phi_y+\frac{2}{3\bar{v}}W_y=
\frac43\frac{\mu(\bar{\theta})}{\bar{v}}\Psi_{1yy}+\frac43(\frac
{\mu(\theta)}v
-\frac{\mu(\bar{\theta})}{\bar{v}})\v{u}_{1y}\\[0.2cm]
\di \quad~~~~~~~~~~~~-\v\sum_{j=1}^2\int\xi_1^2\Theta_{1}^jd\xi+J_1+\frac{2}{3\bar{v}}Y-\bar R_1
\di \doteq\frac43\frac{\mu(\bar{\theta})}{\bar{v}}\Psi_{1yy}+Q_1,\\[0.2cm]
\di \Psi_{i\tau}=\frac{\mu(\bar{\theta})}{\bar{v}}\Psi_{iyy}
+\v(\frac{\mu(\theta)}{v}-\frac{\mu(\bar{\theta})}{\bar{v}}){u}_{iy}+\v^2(N_{i}-\bar{N}_{i})+\v^3 F_{i}
\\[0.2cm]
\di ~~~~~~~~~~~~~~~-\v\int\xi_1\xi_i\Theta_1^{2}d\xi-\bar R_i\doteq
\di \frac{\mu(\bar{\theta})}{\bar{v}}\Psi_{iyy}+Q_i,i=2,3,\\[0.2cm]
\di {W}_{\tau}+p_+\Psi_{1y}=\frac{\k(\bar{\theta})}{\bar{v}}W_{yy}
+(\frac{\k(\theta)}{v}-\frac{\k(\bar{\theta})}{\bar{v}}){\theta}_y
+\v^2(N_{1}-\bar{N}_{1})+\v^3F_{1}
+\frac{4\v}3\frac{\mu({\theta})}{{v}}u_{1y}\Psi_{1y}\\[0.2cm]
\di ~~~~~~~~~~~~+\v^3\sum_{i=2}^3[\frac{\mu(\theta)}{v}u_iu_{iy}-
\v\frac{\mu(\bar{\theta})}{\bar{v}}\bar{u}_i\bar{u}_{iy}]-\v\bar{u}_{1\tau}\Psi_1+J_2
-\v\int\frac12\xi_1|\xi|^2\Theta_1^{2}d\xi\\[0.2cm]
~~~~~~~~~~~~~~~~~~~~~~+\v^2\bar{u}_1
\sum_{j=1}^2\int\xi_1^2\Theta_1^jd\xi-\frac{\k(\bar{\theta})}{\bar{v}}Y_y
+\frac25\v^2\hat{N}_1+\v\bar{u}_1\bar{R}_1-\bar{R}_4\\[0.2cm]
~~~~~~~~~~~~~~~~~
\di \doteq\frac{\k(\bar{\theta})}{\bar{v}}W_{yy}+\frac25\v^2\hat{N}_1+Q_4,
\end{array}\right.
\end{eqnarray}
where
\begin{eqnarray}\label{B4.10}
\left\{
\begin{array}{llll}
\di J_1=\frac{\bar{p}-p_+}{\bar{v}}\Phi_y-[p-\bar{p}+\frac{\bar{p}}{\bar{v}}\Phi_y-\frac{2}{3\bar{v}}(\theta-\bar{\theta})]
=O(1)(\Phi_y^2+(\theta-\bar{\theta})^2+|\v\bar{u}|^4),\\[0.3cm]
\di J_2=(p_+-p)\Psi_{1y}=O(1)(\Phi_y^2+\Psi_{1y}^2+(\theta-\bar{\theta})^2+|\v\bar{u}|^4),\\[0.2cm]
\di Q_1=\frac{4\v}3(\frac
{\mu(\theta)}v
-\frac{\mu(\bar{\theta})}{\bar{v}}){u}_{1y}-\v\sum_{j=1}^2\int\xi_1^2\Theta_{1}^jd\xi+J_1
+\frac{2}{3\bar{v}}Y-\bar R_1,\\[0.2cm]
\di Q_i=\v(\frac{\mu(\theta)}{v}-\frac{\mu(\bar{\theta})}{\bar{v}}){u}_{iy}+\v^2(N_{i}-\bar{N}_{i})+\v^3 F_{i}
-\v\int\xi_1\xi_i\Theta_1^{2}d\xi-\bar R_i,~i=2,3,\\[0.2cm]
\di Q_4=(\frac{\k(\theta)}{v}-\frac{\k(\bar{\theta})}{\bar{v}}){\theta}_y
+\v^2(N_{1}-\bar{N}_{1})+\v^3F_{1}
+\frac{4\v}3\frac{\mu({\theta})}{{v}}u_{1y}\Psi_{1y}\\[0.2cm]
\di ~~~~~~~~~~~~+\v^3\sum_{i=2}^3[\frac{\mu(\theta)}{v}u_iu_{iy}-
\frac{\mu(\bar{\theta})}{\bar{v}}\bar{u}_i\bar{u}_{iy}]-\v\bar{u}_{1\tau}\Psi_1+J_2
-\v\int\frac12\xi_1|\xi|^2\Theta_1^{2}d\xi\\[0.2cm]
\di ~~~~~~~~~~~~~~~~~~~~~~+\v^2\bar{u}_1
\sum_{j=1}^2\int\xi_1^2\Theta_1^jd\xi-\frac{\k(\bar{\theta})}{\bar{v}}Y_y+\v\bar{u}_1\bar{R}_1-\bar{R}_4.
\end{array}\right.
\end{eqnarray}
The equation of  microscopic component $\tilde{G}$ given in \eqref{B1.22} in
the coordinate $(y, \tau)$ becomes
\begin{eqnarray}\label{B1.22-1}
&&v\tilde{G}_\tau-vL_M\tilde{G}=-\frac1{R\theta}P_1[\xi_1
(\frac{|\xi-\v u|^2}{2\theta}\f1\v\zeta_y+\xi\cdot\f1\v\psi_y)M]\nonumber\\
&&~~~~~~~~~~~~~~~~~~~~~~~~+\v u_1G_y-vP_1(\xi_1G_y)+\v vQ(G,G)-v\bar{G}_\tau.
\end{eqnarray}
In the scaling of \eqref{B4.2}, the equation \eqref{B1.13} reads
\begin{equation}\label{B4.11-2}
f_\tau-\v \frac{u_1}{v}f_y+\frac{\xi_1}{v}f_y= \v L_{M}G+\v^2 Q(G,G).
\end{equation}
Set
\begin{equation}\label{B4.11-3}
\tilde f\doteq f-\bar f,
\end{equation}
then from \eqref{B4.11-2} and \eqref{B1.47}, we have
\begin{eqnarray}\label{B4.11}
&&\di v\tilde{f}_\tau-\v u_1\tilde{f}_y+\xi_1\tilde{f}_y=\v vL_{M}\tilde{G}+\v\big[vL_{M}\bar{G}-\bar{v}L_{\bar{M}}\bar{G}_0\big]
+\v^2\big[vQ(G,G)-\bar{v}Q(\bar{G}_0,\bar{G}_0)\big]\nonumber\\
&&\di ~~~~~~~~~~~~~~~~~~~~~~~~~~~~~~~~~~~~~~~-\phi\bar{f}_\tau+\psi\bar{f}_y-\v v\bar{R}_{\bar f}.
\end{eqnarray}

Note that to prove the main theorem in this paper, it is sufficient to prove the following {\it a priori} estimate in the scaled independent variables based on the construction of the approximate  profile.
\begin{theorem}[A priori estimate]\label{thm3.1}
For any sufficiently small and fixed positive constant $\vartheta>0,$ there exist small positive constants $\d_2>0, \v_2>0$ and  a global Maxwellian $M_\ast=M_{[\r_\ast, u_\ast,\t_\ast]}$ such that if $\d\leq \d_2$ and $\v\leq \v_2$, then the Cauchy problem \eqref{B4.9}, \eqref{B1.22-1} and \eqref{B4.11} admits a unique smooth solution satisfying
\begin{equation}\label{B1.50-2}
\begin{cases}
\di \|(\Phi,\Psi, W)(\tau)\|^2_{L^\infty}\leq C\sqrt\d \v,~~~\|(\phi,\psi, \zeta)(\tau)\|^2_{L^2_y} \leq C\sqrt\d\v^2(1+\v^2\tau)^{-1+C_0\sqrt\d},\\[2mm]
\di \|(\phi,\psi, \zeta)_y(\tau)\|_{L^2_y}+\|(\phi,\psi, \zeta)_{yy}(\tau)\|_{L^2_y}
+\sum_{|\alpha|=2}\int_{\mathbb R}\int_{\mathbb{R}^3}\frac{|\partial^\alpha \tilde f|^2}{M_*}d\xi dy
\leq C\sqrt\d\v^3(1+\v^2\tau)^{-\f32+\vartheta+C_0\sqrt\d},\\[2mm]
\di \v\int_{\mathbb R}\int_{\mathbb{R}^3}\frac{|\tilde{G}|^2}{M_*}d\xi dy \leq C\sqrt\d (1+\v^2\tau)^{-\f12},~~
\sum_{|\a|=1}\int_{\mathbb R}\int_{\mathbb{R}^3}\frac{|\partial^\a\tilde{G}|^2}{M_*}d\xi dy\leq C\sqrt\d\v (1+\v^2\tau)^{-\f32+\vartheta+C_0\sqrt\d},
\end{cases}
\end{equation}
where $C, C_0$ are positive constants independent of $\d$ and $\v$.
\end{theorem}

In the next subsection, we will work on the reformulated system
\eqref{B4.9} and \eqref{B4.11}. Since the local existence of the solution can be proved similarly as the
discussion in \cite{Guo} and \cite{Ukai-1}, we will omit it here for
brevity. To prove the global existence, it is sufficient to close the
following a priori estimate:
\begin{equation}\label{B4.13}
\begin{array}{ll}
N(\tau)=&{ \displaystyle\sup_{0\leq s\leq \tau}\Big\{\v^{-1}\|(\Phi,\Psi, W)\|^2_{L^\infty}+\v^{-2}\|(\phi,\psi, \zeta)\|^2_{L^2}+\v^{-3}\|(\phi_y,\psi_y, \zeta_y)\|^2_{L^2}}\\[0.3cm]
&{\di +\|\int_{\mathbb{R}^3}\frac{|\tilde{G}|^2}{M_*}d\xi\|_{L^\infty}
+\int_{\mathbb R}\int_{\mathbb{R}^3}\Big(\sum_{|\alpha|=1}
\v^{-1}\frac{|\partial^\alpha\tilde{G}|^2}{M_*}
+\sum_{|\alpha|=2}\v^{-3}\frac{|\partial^\alpha\tilde{f}|^2}{M_*}\Big)d\xi
dy\Big\}\le \l_0^2},
\end{array}
\end{equation}
where $\l_0$ is positive small constant depending on the
initial data and $M_*$ is a global Maxwellian to be chosen later.

Before proving the a priori estimate \eqref{B4.13}, we list some
 lemmas based on the celebrated H-theorem for later use. The
first one is from \cite{GPS}.

\begin{lemma}\label{lem4.1}
There exists a positive constant $C>0$
such that
\begin{equation*}\label{B4.19-1}
\int_{\mathbb{R}^3}\frac{\nu(|\xi|)^{-1}Q(f,g)^2}{M}d\xi\le
C\Big\{\int_{\mathbb{R}^3}\frac{\nu(|\xi|)f^2}{M}d\xi\cdot\int_{\mathbb{R}^3}\frac{g^2}{M}+
\int_{\mathbb{R}^3}\frac{f^2}{M}d\xi\cdot\int_{\mathbb{R}^3}\frac{\nu(|\xi|)g^2}{M}\Big\},
\end{equation*}
where $M$ can be any Maxwellian so that the above integrals are
well defined.
\end{lemma}

Based on Lemma \ref{lem4.1}, the following three lemmas are from
\cite{Liu-Yang-Yu-Zhao}.

\begin{lemma}\label{lem4.2}
If $\theta/2<\theta_*<\theta$, then
there exist two positive constants
$\bar{\sigma}=\bar{\sigma}(\rho,u,\theta;\rho_*,u_*,\theta_*)>0$
and $\eta_0=\eta_0(\rho,u,\theta;\rho_*,u_*,\theta_*)>0$ such that
if $|\r-\r_*|+|\v u- u_*|+|\theta-\theta_*|<\eta_0$, we have for
$h(\xi)\in N^\bot$,
\begin{equation*}\label{B4.19}
-\int_{\mathbb{R}^3}\frac{hL_Mh}{M_*}d\xi\ge
\bar{\sigma}\int_{\mathbb{R}^3}\frac{\nu(|\xi|)h^2}{M_*}d\xi,
\end{equation*}
where $M_*=M_{[\rho_*,u_*,\theta_*]}$ and the definition of
$M=M_{[\rho,\v u,\theta]}$ can be found in \eqref{B1.4-1}.
\end{lemma}

\begin{lemma}\label{lem4.3}
Under the assumptions in Lemma \ref{lem4.2}, we
have
\begin{equation*}\label{B4.20}
\left\{
\begin{array}{l}
{ \di \int_{\mathbb{R}^3}\frac{\nu(|\xi|)}{M}|L_M^{-1}h|^2d\xi\le \bar{\sigma}^{-2}\int_{\mathbb{R}^3}\frac{\nu(|\xi|)^{-1}h^2}{M}d\xi},\\
{
\di \int_{\mathbb{R}^3}\frac{\nu(|\xi|)}{M_*}|L_M^{-1}h|^2d\xi\le
\bar{\sigma}^{-2}\int_{\mathbb{R}^3}\frac{\nu(|\xi|)^{-1}h^2}{M_*}d\xi},
\end{array}
\right.
\end{equation*}
for each $h(\xi)\in N^\bot$.
\end{lemma}

\begin{lemma}\label{lem4.4}
Under the conditions in Lemma \ref{lem4.2}, there
exists a constant $C>0$ such that for positive constants $k$ and
$\lambda$, we have
\begin{equation*}\label{B4.21}
|\int_{\mathbb{R}^3}\frac{g_1P_1(|\xi|^kg_2)}{M_*}d\xi-\int_{\mathbb{R}^3}\frac{g_1|\xi|^kg_2}{M_*}d\xi|\le
C\int_{\mathbb{R}^3}\frac{\lambda|g_1|^2+\lambda^{-1}|g_2|^2}{M_*}d\xi.
\end{equation*}
\end{lemma}

Note that \eqref{B4.13} also gives
the a priori estimates on $\|(\phi_\tau,\psi_\tau,\zeta_\tau)\|$,
$\|\partial^\alpha(\phi,\psi,\zeta)\|$ and
$\di\int\int\frac{|\partial^\alpha\tilde{G}|^2}{M_*}d\xi dx$ ($|\alpha|=2$).
In fact, from \eqref{B4.2},  \eqref{B1.14} and \eqref{B4.13}, one has
\begin{eqnarray}\label{B4.14}
&&\|(\phi_\tau,\psi_\tau,\zeta_\tau)\|^2\leq C\|(v_\tau,\v u_\tau,\t_\tau)\|^2+C\d \v^3(1+\v^2\tau)^{-\f32}\nonumber\\
&&\leq C\Big(\|(p_y-\bar p_y, \v pu_{1y}-\v\bar p\bar u_{1y})\|^2
+\|(\bar p_y, \v\bar p\bar u_{1y})\|^2
+\v^2\int\int\frac{|\tilde G_y|^2+|\bar G_y|^2}{M_*}d\xi dy\Big)+C\d \v^3(1+\v^2\tau)^{-\f32}\nonumber\\
&&\leq C\Big(\|(\phi_y,\psi_y,\zeta_y)\|^2+\d\v^2\|(\phi,\psi,\zeta)\|^2+\v^2\int\int\frac{|\tilde G_y|^2}{M_*}d\xi dy\Big)+C\d \v^3(1+\v^2\tau)^{-\f32}\nonumber\\
&&\leq C(\d+\l_0^2)\v^3,
\end{eqnarray}
where we have used the fact that
\begin{equation*}\label{B4.15}
\int\Big(\int\xi_1^2G_yd\xi\Big)^2dy\le C\int \int \frac{G_y^2}{M_*}d\xi dy.
\end{equation*}
To derive the a priori assumption on $\|\partial^\alpha(\phi,\psi,\zeta)\|$, ($|\alpha|=2$), we use the definition of $\rho$, $m=\v \rho u$ and
$\rho(\theta+\frac12|\v u|^2)$. Let $|\alpha|=2$, by \eqref{B1.4}, one can obtain
\begin{eqnarray}\label{B4.16}
&&\|\partial^\alpha(\rho, m,\rho(\theta+\frac12|\v u|^2))\|^2
\leq C\int\int\frac{|\partial^\alpha f|^2}{M_*}d\xi dy\nonumber\\
&&\leq C\int\int\frac{|\partial^\alpha\tilde{f}|^2}{M_*}d\xi+C\d\v^3(1+t)^{-\f32}\leq C(\l_0^2+\d)\v^3.
\end{eqnarray}
This yields that
\begin{equation}\label{B4.17}
\sum_{|\a|=2}\|\partial^\alpha(\phi,\psi,\zeta)\|^2\leq C\int\int\frac{|\partial^\alpha\tilde{f}|^2}{M_*}d\xi dy+\sum_{|\beta|=1}\|\partial^\beta(\phi,\psi,\zeta)\|^2+C\d\v^3(1+t)^{-\f32}\leq C(\l_0^2+\d)\v^3.
\end{equation}
Finally, one has
\begin{eqnarray}\label{B4.18}
&&\v^2\int\int\frac{|\partial^\alpha\tilde{G}|^2}{M_*}d\xi dy\leq
C\int\int\frac{|\partial^\alpha \tilde{f}|^2}{M_*}d\xi
dy+C\int\int\frac{|\partial^\alpha(M-\bar M)|^2}{M_*}d\xi dy\\
&&\leq C\int\int\frac{|\partial^\alpha \tilde{f}|^2}{M_*}d\xi dy
+C\sum_{|\a|=1,2}\|\partial^\a(\phi,\psi,\zeta)\|^2+C\d\v^4(1+t)^{-2}
\leq C(\varepsilon_0+\delta)^2\v^3,~~ |\alpha|=2.\nonumber
\end{eqnarray}

\subsection{Lower order estimate}
We are now ready to derive the lower order estimate. Multiplying
$\eqref{B4.9}_1$ by $p_+\Phi$, $\eqref{B4.9}_2$ by $\bar{v}\Psi_1$, $\eqref{B4.9}_3$ by
$\Psi_i$, $\eqref{B4.9}_4$ by $\frac{2}{3p_+}W$ with $p_+=\f23$ respectively and adding
all the  equations, one can obtain
\begin{eqnarray}\label{B4.22}
\begin{array}{ll}
\di \Big(\frac{p_+}2\Phi^2+\frac1{3p_+}W^2+\frac{\bar{v}}2\Psi_1^2+\frac12\sum_{i=2}^3\Psi_i^2\Big)_\tau
+\frac{4\mu(\bar{\theta})}{3}\Psi_{1y}^2+\sum_{i=2}^3\frac{\mu(\bar{\theta})}{\bar{v}}\Psi_{iy}^2
+\frac{2\k(\bar{\theta})}{3p_+\bar{v}}W_y^2\\[3mm]
\di =\frac25\v^2\hat{N}_1(-\Phi+\frac2{3p_+}W)+\frac12\bar{v}_\tau\Psi_1^2+{\bar v}Q_1\Psi_1 +\sum_{i=2}^3{Q}_i\Psi_i+\frac2{3p_+}WQ_4 \\[3mm]
\di ~~-(\frac{4\mu(\bar{\theta})}{3})_y\Psi_{1}\Psi_{1y}
-\sum_{i=2}^3(\frac{\mu(\bar{\theta})}{\bar{v}})_y\Psi_{i}\Psi_{iy}
-(\frac{2\lambda(\bar{\theta})}{3p_+\bar{v}})_yWW_y  +(\cdots)_y.
\end{array}
\end{eqnarray}
Here and in the sequel the notation $(\cdots)_y$ represents
the term in the conservative form so that it vanishes after
integration. Since it has no effect on the energy estimates, we do
not write them out in detail.

Note that the term $Q_1\Psi_1$ contains $(1+t)^{-1}\Psi_1$ which
can not be controlled by the dissipation from the viscosity and
heat conductivity. So is the term
$\tilde{N}_1(-\Phi+\frac2{3p_+}W)$. As we will see later, an
intrinsic dissipation associated with the  profile is
derived by the diagonal method and weighted energy estimate to
control the above two terms. Let us consider the equations for the
conservation of the mass, the first component of velocity and
energy by defining
\begin{equation*}\label{B4.23-1}
V=(\Phi,\Psi_1,W)^t,
\end{equation*}
where $(\cdot,\cdot,\cdot)^t$ means the transpose of the vector $(\cdot,\cdot,\cdot)$. Then from \eqref{B4.9}, we have
\begin{equation}\label{B4.23}
V_\tau+A_1V_y=A_2V_{yy}+A_3,
\end{equation}
where
\begin{equation*}\label{B4.24}
A_1=\left(\begin{array}{cccc}&0&-1&0\\&-\frac{p_+}{\bar{v}}&0&\frac{2}{3\bar{v}}\\
&0&p_+&0\end{array}\right), \quad
A_2=\left(\begin{array}{cccc}
&0&0&0\\&0&\frac{4\mu(\bar{\theta})}{3\bar{v}}&0\\
&0&0&\frac{\lambda(\bar{\theta})}{\bar{v}}\end{array}\right),
\end{equation*}
\begin{equation*}\label{B4.25}
A_3=(-\frac3{5}\v^2\hat{N}_1,~~Q_1,~~Q_4+\frac25\v^2\hat{N}_1)^t.
\end{equation*}
Direct computation shows that the eigenvalues of the matrix $A_1$
are $\lambda_1,0,\lambda_3$. Here
$\lambda_3=-\lambda_1=\sqrt{\frac{5p_+}{3\bar{v}}}$. The
corresponding normalized left and right eigenvectors can be chosen
as
\begin{eqnarray*}
l_1=\sqrt{3/10}(-1,-\frac{5}{3\lambda_3},\frac2{3p_+}),~
l_2=\sqrt{2/5}(1,0,\frac1{p_+}),~
l_3=\sqrt{3/10}(-1,\frac{5}{3\lambda_3},\frac2{3p_+}),\label{B4.26}\\r_1=\sqrt{3/10}(-1,-\lambda_3,p_+)^t,~
r_2=\sqrt{2/5}(1,0,\frac32{p_+})^t,~r_3=\sqrt{3/10}(-1,\lambda_3,p_+)^t,\label{B4.27}
\end{eqnarray*}
such that
\begin{eqnarray*}\label{B4.29}
l_ir_j=\delta_{ij},~i,j=1,2,3,\quad LA_1R=\Lambda=\left(\begin{array}{cccc}&\lambda_1&0&0\\&0&0&0\\
&0&0&\lambda_3\end{array}\right),
\end{eqnarray*}
with
$$L=(l_1,l_2,l_3)^t,~ R=(r_1,r_2,r_3).$$
Let
\begin{eqnarray*}\label{B4.29-1}
B=LV=(b_1,b_2,b_3),
\end{eqnarray*}
 then multiplying the
equations \eqref{B4.23} by the matrix $L$ yields that
\begin{eqnarray}\label{B4.30}
B_\tau+\Lambda B_y=LA_2RB_{yy}+2LA_2R_yB_y+[(L_\tau+\Lambda
L_y)R+LA_2R_{yy}]B+LA_3.
\end{eqnarray}
A direct computation shows that $LA_2R=A_4 \mbox{ is a non-negative matrix}$.
From \eqref{B4.30}, we will apply weighted energy method to derive an
intrinsic dissipation. Since we have assumed that
$\hat\t_{y}> 0$. Let $v_1=\f{\hat\t}{\theta_+}$, then $|v_1-1|\le
C\delta$. Multiplying \eqref{B4.30} by
$\bar{B}=(v_1^nb_1,b_2,v_1^{-n}b_3)$ with a large positive integer
$n$ which will be chosen later, we have
\begin{eqnarray}\label{B4.32}
\begin{array}{ll} \Big(\frac{1}2v_1^nb_1^2+\frac12b_2^2+\frac{1}{2}v_1^{-n}b_3^2\Big)_\tau
-(\frac{v_1^n}2)_\tau b_1^2-(\frac{v_1^{-n}}2)_\tau b_3^2
+\bar{B}_yA_4B_y+\bar{B}A_{4y}B_y\\[0.2cm]
-\frac{1}2v_1^{n-1}(n\lambda_1v_{1y}+v_1\lambda_{1y})b_1^2
+\frac{1}2v_1^{-n-1}(n\lambda_3v_{1y}-v_1\lambda_{3y})b_3^2\\[0.2cm]
=2\bar{B}LA_2R_yB_y+\bar{B}[L_tR+LA_2R_{yy}]B+\bar{B}\Lambda
L_xRB+\bar{B}LA_3+(\cdots)_x.
\end{array}
\end{eqnarray}
Let
\begin{eqnarray*}
&&E_1=\int\Big(\frac{p_+}2\Phi^2+\frac1{3p_+}W^2+\frac{\bar{v}}2\Psi_1^2+\frac12\sum_{i=2}^3\Psi_i^2\Big)dy
+\int(\frac{v_1^n}2b_1^2
+\frac12b_2^2+\frac{v_1^{-n}}{2}b_3^2)dy,\label{B4.33}\\
&&K_1=\int(\frac{4\mu(\bar{\theta})}{3}\Psi_{1y}^2+\sum_{i=2}^3\frac{\mu(\bar{\theta})}{\bar{v}}\Psi_{iy}^2
+\frac{2\lambda(\bar{\theta})}{3p_+\bar{v}}W_y^2+{B}_yA_4B_y)dy.\label{B4.34}
\end{eqnarray*}
Note that
\begin{eqnarray}\label{B4.35}
&&|\int(\bar{B}-B)_yA_4B_ydy\leq C\delta\int|B_y|^2dy+C\delta^{-1}\int|\hat\theta_y|^2|B|^2dy\nonumber\\
&&~~~~~~~~~~~~~~~~~~~~~~~~\leq  C\v^2\delta(1+t)^{-1}E_1+C\delta
K_1+C\delta\int|\Phi_y|^2dy.
\end{eqnarray}

Similarly, the terms in the last second line of \eqref{B4.32}, $\bar{B}A_{4y}B_y$, $\bar{B}LA_2R_yB_y$ and
$\bar{B}[L_\tau R+LA_2R_{yy}]B$ satisfy the same estimate. For
$\bar{B}\Lambda L_yRB$ and $\bar{B}LA_3$, we need to use the
explicit presentation. By the choice of the characteristic matrix
$L$ and $R$, we have
$$
\Lambda L_yR=\frac12\lambda_{3y}\left(\begin{array}{cccc}&1&0&-1\\&0&0&0\\
&1&0&-1\end{array}\right),~~~~~
LA_3=\left(\begin{array}{c}{ \sqrt{\frac{2}{15}}\frac1{p_+}(\v^2\hat{N}_1+Q_4)-\sqrt{\frac{5}{6}}\frac{Q_1}{\lambda_3}}\\
{ \sqrt{\frac{2}{5}}\frac{Q_4}{p_+} }\\
{ \sqrt{\frac{2}{15}}\frac1{p_+}(\v^2\hat{N}_1+Q_4)+\sqrt{\frac{5}{6}}\frac{Q_1}{\lambda_3}}\end{array}\right).
$$
Thus
\begin{eqnarray}
\bar{B}\Lambda
L_yRB=\frac12\lambda_{3y}(v_1^nb_1^2
+v_1^{-n}b_1b_3-v_1^nb_1b_3-v_1^{-n}b_3^2),\nonumber\\
\bar{B}LA_3=\sqrt{\frac{2}{15}}\frac1{p_+}\v^2\hat{N}_1
(v_1^nb_1+v_1^{-n}b_3)+q_1v_1^nb_1+q_2b_2+q_3v_1^{-n}b_3,\label{B4.37}
\end{eqnarray}
where
$$
q_1=\sqrt{\frac{2}{15}}\frac1{p_+}Q_4-\sqrt{\frac{5}{6}}\frac{Q_1}{\lambda_3},~q_2=\sqrt{\frac{2}{5}}\frac{Q_4}{p_+},
~q_3=\sqrt{\frac{2}{15}}\frac1{p_+}Q_4+\sqrt{\frac{5}{6}}\frac{Q_1}{\lambda_3}.
$$
 Combine \eqref{B4.22}, \eqref{B4.32}, \eqref{B4.35}-\eqref{B4.37}, we
have by choosing $n$ sufficiently large,
\begin{eqnarray}\label{B4.38}
 E_{1\tau}+\frac12
K_1+2\int|\hat\t_{y}|(b_1^2+b_3^2)|dy\leq C\v^2\delta(1+t)^{-1}(E_1+1)
+C\delta\int\Phi_y^2dy+I_{nf},
\end{eqnarray}
where
\begin{eqnarray}\label{B4.39}
I_{nf}=\int\bar{v}Q_1\Psi_1dy+\int\sum_{i=2}^3{Q}_i\Psi_idy
+\int\frac2{3p_+}WQ_4dy+\int
(q_1v_1^nb_1+q_2b_2+q_3v_1^{-n}b_3)dy.
\end{eqnarray}
Here we have used the fact that
\begin{equation}\label{B4.40}
-\Phi+\frac2{3p_+}W=\sqrt{5/6}(b_1+b_3),
\end{equation}
and
\begin{equation*}\label{B4.41}
\v^2\int|\hat{N}_1|(|b_1|+|b_3|)dy\leq
C\delta\int|\hat\t_{y}|(b_1^2+b_3^2)|dy+C\v^2\delta(1+t)^{-1},
\end{equation*}
and for $n$ large enough,
\begin{eqnarray*}\label{B2.34}
-\frac12v_1^{n-1}(n\lambda_1v_{1y}+2v_1\lambda_{1y})b_1^2+\frac12v_1^{-n-1}
(n\lambda_3v_{1y}-2v_1\lambda_{3y})b_3^2-\tilde{B}\Lambda L_yRB\geq 3|\hat\theta_{y}|(b_1^2+b_3^2).
\end{eqnarray*}

Even though $Q_1$ contains the term $R_1$ with the decay rate
 $\frac{\v^2}{1+t}$,  the terms in \eqref{B4.39} involving $Q_1$ have factor
$b_1$ or $b_3$ because
\begin{equation}\label{B4.42}
\Psi_1=\sqrt{3/10}\lambda_3(b_3-b_1).
\end{equation}
Thus the
terms $\bar{v}Q_1\Psi_1$, $q_1v_1^nb_1$ and $q_3v_1^{-n}b_3$ can
be controlled by the intrinsic dissipation on $b_1$ and $b_3$ as
shown later.  The estimates on the other terms involving $Q_i ~(i=2,3,4)$ are straightforward because from \eqref{B1.42}-\eqref{B1.44} and \eqref{B4.10}, they decay at least in the order of
$\v^3(1+t)^{-3/2}$.
For brevity, we only estimate $\di \int
\bar{v}Q_1\Psi_1dy$ and $\di \int q_2b_2dy$ as follows for
illustration.

\

\noindent {\bf \underline{Estimate on $\di \int \bar{v}Q_1\Psi_1dy$}:}

\

From \eqref{B4.42}, we have
\begin{equation}\label{B4.43}
\int
\bar{v}Q_1\Psi_1dy=\sqrt{\frac{3}{10}}\int\bar{v}Q_1\lambda_3(b_3-b_1)dy.
\end{equation}
Here we only consider the integral
\begin{equation*}\label{B4.44}
I_1=\int\bar{v}Q_1\lambda_3b_1dy,
\end{equation*}
and the other term in \eqref{B4.43}
can be estimated similarly. By the definition of $Q_1$ in \eqref{B4.10},
we have
\begin{eqnarray*}\label{B4.44-1}
&&I_1=\int
\bar{v}\lambda_3b_1[\frac{4\v}3(\frac{\mu(\theta)}{v}-\frac{\mu(\bar{\theta})}{\bar{v}})u_{1y}
+J_1+\frac2{3\bar{v}}Y]dy-\int\bar{v}\lambda_3b_1\bar{R}_1dy-\v\int\bar{v}\lambda_3b_1\sum_{j=1}^2\int\xi_1^2\Theta_1^jd\xi
dy\nonumber\\
&&~~~={  I^1_1+I_1^2+I_1^3}.
\end{eqnarray*}
Since
\begin{eqnarray*}\label{N-2.36}
&&\int\big|\f{4\v}3\big(\f{\mu(\t)}{v}-\f{\mu(\bar\t)}{\bar v}\big)u_{1y}\big|\cdot|b_1|dy\nonumber\\
&&\leq C(\d+\lambda_0) (K_1+\|\Phi_y\|^2)+C\d\v^2(1+t)^{-1} E_1+C(\d+\lambda_0)\|\psi_{1y}\|^2+C\d\v^5(1+t)^{-\f52},
\end{eqnarray*}
and
\begin{eqnarray*}\label{N-2.37}
&&\int\big(|J_1|+|\f{Y}{\bar v}|\big)\cdot|b_1|dy
\leq C(\d+\l_0) (K_1+\|\Phi_y\|^2)+C\d\v^2(1+t)^{-1}E_1+C\d\v^5(1+t)^{-\f52},
\end{eqnarray*}
we obtain
\begin{eqnarray}\label{B4.45}
I_1^1&\leq& C(\d+\l_0) (K_1+\|\Phi_y\|_{L^2}^2)+C(\d+\l_0) \|\psi_{1y}\|_{L^2}^2\nonumber\\
&+&C\d\v^2(1+t)^{-1}E_1+C\d\v^2(1+t)^{-1}.
\end{eqnarray}
On the other hand, from \eqref{B1.42}, we have
$$
\bar R_1=O(1)\delta\v^2(1+t)^{-1}e^{-\frac{x^2}{4c(\theta_\pm)(1+t)}},\quad
{\rm as}~x\to \pm\infty.
$$
From \eqref{N-1.4}, $\hat\t_y$ satisfies
$$
|\hat\t_y|=O(\delta)\v(1+t)^{-\f12}e^{-\frac{x^2}{4a(\theta_\pm)(1+t)}},\quad
{\rm as}~x\to \pm\infty.
$$
Thus, by \eqref{B1.12} and the assumption on the  profile, we
have
\begin{equation}\label{B4.46}
\k(\theta_\pm)=\f52\mu(\t_\pm)>\frac{5}{4}\mu(\theta_\pm).
\end{equation}
Since $a(\theta_\pm)=\f{3\k(\theta_\pm)}{5\theta_\pm},$ $b(\theta_\pm)=\max\{a(\theta_\pm),\f{\mu(\theta_\pm)}{\theta_\pm}\}$ and $c(\t_\pm)=\max\{a(\t_\pm),\f12b(\t_\pm)\}$ , it follows from \eqref{B4.46} that
$a(\theta_\pm)>\frac23c(\theta_\pm)$, which leads to
\begin{equation}\label{B4.47}
|I_1^2|\le
\f1{16}\int|\hat\t_y|b_1^2dy+C\delta\v^2(1+t)^{-1}.
\end{equation}

\

We now estimate the integral $I_1^3$. Let $M_*$ be a global
Maxwellian with the state $(\rho_*,u_*,\theta_*)$ satisfying
$\frac12\theta<\theta_*<\theta$ and
$|\r-\r_*|+|\v u-u_*|+|\theta-\theta_*|\le \eta_0$ so that Lemma \ref{lem4.2} holds. Note that,
\begin{equation}\label{B4.48}
I_1^3=-\v\int\bar{v}\lambda_3b_1\int\xi_1^2\Theta_1^1d\xi
dy-\v\int\bar{v}\lambda_3b_1\int\xi_1^2\Theta_1^2d\xi
dy=:I_1^{31}+I_1^{32}.
\end{equation}
The estimation on $I_1^{31}$ is straightforward by using the
intrinsic dissipation on $b_1$ and \eqref{B1.23}.
\begin{eqnarray}\label{B4.49}
&&|I_1^{31}|=|-\v\int\bar{v}\lambda_3b_1\int\xi_1^2L_M^{-1}[\frac1vP_1(\xi_1\bar{G}_y)
-\v Q(\bar{G},\bar{G})]d\xi dy|\nonumber\\
&&\leq C\int
|b_1|\Big(|(\v\bar{u}_y,\bar{\theta}_y)|^2+|(\v\bar{u}_{yy},\bar{\theta}_{yy}|
+|(\v\bar{u}_y,\bar{\theta}_y)||(v_y,\v u_y,\theta_y)|\Big)dy\nonumber\\
&&\leq C\delta\int|\hat\t_y|b_1^2dy+
C\delta\v^2(1+t)^{-1}+C\delta\|(\phi_y,\psi_y,\zeta_y)\|^2.
\end{eqnarray}
The estimation on $I_1^{32}$ is more complicated and it will be
 divided into five parts as follows. From \eqref{B1.23}, it holds that
\begin{eqnarray}\label{B4.51}
&&I_1^{32}=-\v\int\bar{v}\lambda_3b_1\int\xi_1^2L_M^{-1}(G_\tau)d\xi
dy+\v\int\bar{v}\lambda_3b_1\int\xi_1^2\frac{\v u_1}{v}L_M^{-1}(G_y)d\xi
dy\nonumber\\
&&~~~~~~~~~
 -\v\int\bar{v}\lambda_3b_1\frac1v\int\xi_1^2L_M^{-1}[P_1(\xi_1\tilde{G}_y)]d\xi
dy+\v^2\int\bar{v}\lambda_3b_1\int\xi_1^2L_M^{-1}[Q(\tilde{G},\tilde{G})]d\xi
dy\nonumber\\
&&~~~~~~~~~
+2\v^2\int\bar{v}\lambda_3b_1\int\xi_1^2L_M^{-1}[Q(\tilde{G},\bar{G})]d\xi
dy=:\sum_{i=1}^5I_1^{32i}.
\end{eqnarray}
For the integral $I_1^{321}$, one has
\begin{eqnarray}\label{B4.52}
I_1^{321}=-\v\int\bar{v}\lambda_3b_1\int\xi_1^2L_M^{-1}(\tilde{G}_\tau)d\xi
dy-\v\int\bar{v}\lambda_3b_1\int\xi_1^2L_M^{-1}(\bar{G}_\tau)d\xi
dy=:I_1^{3211}+I_1^{3212}.
\end{eqnarray}
Note that the linearized operator $L_M^{-1}$ satisfies that, for any
$h\in N^\bot$,
\begin{eqnarray}\label{B4.53}
\begin{array}{l}
(L_M^{-1}h)_\tau=L_M^{-1}(h_\tau)-2L_M^{-1}\{Q(L_M^{-1}h,M_\tau)\},\\[2mm]
(L_M^{-1}h)_y=L_M^{-1}(h_y)-2L_M^{-1}\{Q(L_M^{-1}h,M_y)\}.
\end{array}
\end{eqnarray}
Then it follows that
\begin{eqnarray}\label{B4.54}
\begin{array}{ll}
\di I_1^{3211}=-\v\int\bar{v}\lambda_3b_1\int\xi_1^2\big(L_M^{-1}\tilde{G}\big)_\tau d\xi
dy-2\v\int\bar{v}\lambda_3b_1\int\xi_1^2L_M^{-1}\{Q(L_M^{-1}\tilde{G},M_\tau)\}d\xi
dy\\[3mm]
\di ~~~~~~~=-\big(\v\int\bar{v}\lambda_3b_1\int\xi_1^2L_M^{-1}\tilde{G}d\xi
dy\big)_\tau+\v\int(\bar{v}\lambda_3b_1)_\tau\int\xi_1^2L_M^{-1}\tilde{G}d\xi
dy\\[3mm]
\di ~~~~~~~~~~~~~-2\v\int\bar{v}\lambda_3b_1\int\xi_1^2L_M^{-1}\{Q(L_M^{-1}\tilde{G},M_\tau)\}d\xi
dy.
\end{array}
\end{eqnarray}
The H\"{o}lder inequality and Lemma \ref{lem4.3} yield that
\begin{equation*}\label{B4.55}
|\int\xi_1^2L_M^{-1}\tilde{G}d\xi|^2\le
C\int\xi_1^4\nu(|\xi|)^{-1}M_*d\xi\cdot\int\frac{\nu(|\xi|)}{M_*}|L_M^{-1}\tilde{G}|^2d\xi\le
C\int\frac{\nu(|\xi|)}{M_*}|\tilde{G}|^2d\xi.
\end{equation*}
Moreover, from Lemmas \ref{lem4.1}-\ref{lem4.3}, one has
\begin{equation}\label{B4.56}
\begin{array}{ll}
{
\di |\int\xi_1^2L_M^{-1}\{Q(L_M^{-1}\tilde{G},M_\tau)\}d\xi|^2\le
C\int\frac{\nu(|\xi|)}{M_*}|L_M^{-1}\{Q(L_M^{-1}\tilde{G},M_\tau)\}|^2d\xi}
\\[3mm]
\quad{
\di \leq C\int\frac{\nu(|\xi|)^{-1}}{M_*}|Q(L_M^{-1}\tilde{G},M_\tau)|^2d\xi\le
C\int\frac{\nu(|\xi|)}{M_*}|L_M^{-1}\tilde{G}|^2d\xi\cdot\int\frac{\nu(|\xi|)}{M_*}|M_\tau|^2d\xi
}\\[3mm]
\quad {\di  \le
C(v_\tau^2+\v^2u_\tau^2+\theta_\tau^2)\int\frac{\nu(|\xi|)^{-1}}{M_*}|\tilde{G}|^2d\xi.
}\end{array}
\end{equation}
Combining \eqref{B4.54}-\eqref{B4.56} gives that
\begin{equation}\label{B4.57}
\begin{array}{ll}
\di I_1^{3211}\leq
-\Big(\v\int\bar{v}\lambda_3b_1\int\xi_1^2L_M^{-1}\tilde{G}d\xi
dy\Big)_\tau+C\b\|(\Phi_\tau,\Psi_\tau,W_\tau)\|^2+C\delta\v^2(1+t)^{-1}E_1\\[2mm]
\quad
 \di ~~~~~~~~~~~~~~~~~~~~~~+C_{\b}~\v^2\int\int\frac{\nu(|\xi|)}{M_*}|\tilde{G}|^2d\xi
dy+C\l_0^2|(\phi_\tau,\psi_\tau,\zeta_\tau)\|^2,
\end{array}
\end{equation}
where and in the sequel $\b$ is small positive constant to be
chosen later and $C_\b$ is a positive constant depending on $\b$. By
the definition of $\bar{G}$ in \eqref{B1.23}, similar to the estimate in
\eqref{B4.49}, one has
\begin{equation}\label{B4.58}
\begin{array}{l}
{\di |I_1^{3212}|=|\v\int\bar{v}\lambda_3b_1\int\xi_1^2L_M^{-1}(\bar{G}_\tau)d\xi
dy|}\\[2mm]
\di \quad { \leq
C\delta\v^2(1+t)^{-1}E_1+C\delta\v^3(1+t)^{-\frac32}+C\delta\|(\phi_\tau,\psi_\tau,\zeta_\tau)\|^2.}
\end{array}
\end{equation}
Substituting \eqref{B4.57} and \eqref{B4.58} into \eqref{B4.52}  implies that
\begin{equation}\label{B4.59}
\begin{array}{ll}
\di I_1^{321}\leq
-\Big(\v\int\bar{v}\lambda_3b_1\int\xi_1^2L_M^{-1}\tilde{G}d\xi
dy\Big)_\tau+C\b\|(\Phi_\tau,\Psi_\tau,W_\tau)\|^2+C\delta\v^2(1+t)^{-1}E_1\\[2mm]
\di \quad~~~~~~~+C_{\b}~\v^2\int\int\frac{\nu(|\xi|)}{M_*}|\tilde{G}|^2d\xi
dy+C(\d+\l_0)\|(\phi_\tau,\psi_\tau,\zeta_\tau)\|^2+C\delta\v^3(1+t)^{-\frac32}.
\end{array}
\end{equation}
The estimation on  $I_1^{32i} ~(i=2,4,5)$ is straightforward by
using the Cauchy inequality and Lemmas \ref{lem4.1}-\ref{lem4.3}. First, it holds that
\begin{eqnarray}\label{B4.60}
&&|I_1^{322}|\leq
C\delta\v^2(1+t)^{-1}E_1+C\l_0K_1+C\v^2\int\int\frac{\nu(|\xi|)}{M_*}|\tilde{G}_y|^2d\xi
dy\nonumber\\
&&~~~~~~~~~~~~~~~~+C\d\v^3(1+t)^{-\f32}+C\d\v^2\|(\phi,\psi,\zeta)_y\|^2.
\end{eqnarray}
Since
\begin{eqnarray*}\label{B4.61}
\begin{array}{l}
{\di
|\int\xi_1^2L_M^{-1}\{Q(\tilde{G},\bar{G})\}d\xi|^2\le
C\int\frac{\nu(|\xi|)}{M_*}|L_M^{-1}\{Q(\tilde{G},\bar{G})\}|^2d\xi}
\\[2mm]
\quad{\di
\le
C\int\frac{\nu(|\xi|)^{-1}}{M_*}|Q(\tilde{G},\bar{G})|^2d\xi\le
C\int\frac{\nu(|\xi|)}{M_*}|L_M^{-1}\tilde{G}|^2d\xi\cdot\int\frac{\nu(|\xi|)}{M_*}|\bar{G}|^2d\xi
}\\[2mm]
\di \quad { \le
C|(\v\bar{u}_x,\bar{\theta}_x)|^2\int\frac{\nu(|\xi|)}{M_*}|\tilde{G}|^2d\xi,
}\end{array}
\end{eqnarray*}
and
\begin{eqnarray*}\label{B4.62}
\begin{array}{l}
{\di
|\int\xi_1^2L_M^{-1}\{Q(\tilde{G},\tilde{G})\}d\xi|\le
C(\int\frac{\nu(|\xi|)}{M_*}|L_M^{-1}\{Q(\tilde{G},\tilde{G})\}|^2d\xi)^{\frac12}}
\\[2mm]
\quad{\di
\le
C(\int\frac{\nu(|\xi|)^{-1}}{M_*}|Q(\tilde{G},\tilde{G})|^2d\xi)^{\frac12}\le
C\int\frac{\nu(|\xi|)}{M_*}|\tilde{G}|^2d\xi,
}\end{array}
\end{eqnarray*}
it follows that
\begin{eqnarray}\label{B4.63}
|I_1^{324}|+|I_1^{325}|\le
C(\delta+\l_0)\v^2\int\int\frac{\nu(|\xi|)}{M_*}|\tilde{G}|^2d\xi
dy+C\delta\v^2(1+t)^{-1}E_1.
\end{eqnarray}
The estimate on $I_1^{323}$ is similar to the one for $I_1^{321}$.
First, notice that
\begin{equation*}\label{B4.64}
P_1(\xi_1\tilde{G}_y)=\{P_1(\xi_1\tilde{G})\}_y
+\sum_{j=0}^4\langle\xi_1\tilde{G},\chi_j\rangle P_1(\chi_{jy}).
\end{equation*}
From \eqref{B4.53} and Lemmas \ref{lem4.1}-\ref{lem4.3}, we have
\begin{equation}\label{B4.65}
\begin{array}{l}
\di I_1^{323}=\v\int(\frac{\bar{v}}{v}\lambda_3b_1)_y\int\xi_1^2L_M^{-1}[P_1(\xi_1\tilde{G})]d\xi
dy-\v\int\frac{\bar{v}}{v}\lambda_3b_1\int\xi_1^2L_M^{-1}[\sum_{j=0}^4\langle\xi_1\tilde{G},\chi_j \rangle P_1(\chi_{jy})]d\xi dy\\[2mm]
\di \quad ~~~~~~~~~~~~~~~~~~~~-2\v\int\frac{\bar{v}}{v}\lambda_3b_1\int\xi_1^2L_M^{-1}\{Q(L_M^{-1}[P_1(\xi_1\tilde{G})],M_y)\}d\xi
dy\\[2mm]
\di \leq C_{\b}~\v^2\int\int\frac{\nu(|\xi|)}{M_*}|\tilde{G}|^2d\xi
dy +C\delta\v^2(1+t)^{-1}E_1+C(\l_0+\b)(K_1+\|\Phi_y\|^2)+C\l_0\|(\phi_y,\psi_y,\zeta_y)\|^2,
\end{array}
\end{equation}
where we have used the fact that
$$
|\langle\xi_1\tilde{G},\chi_j\rangle|^2\le
C\int\frac{\nu(|\xi|)\tilde{G}^2}{M_*}d\xi.
$$
Substituting \eqref{B4.59}, \eqref{B4.60}, \eqref{B4.63} and \eqref{B4.65} into \eqref{B4.51} gives that
\begin{equation*}\label{B4.66}
\begin{array}{l}
\di I_1^{32}\leq
-\Big(\v\int\bar{v}\lambda_3b_1\int\xi_1^2L_M^{-1}\tilde{G}d\xi
dy\Big)_\tau+C\delta\v^2(1+t)^{-1}E_1+C(\l_0+\b)(K_1+\|\Phi_y\|^2)\\[2mm]
\di ~~~~~~~~~~+C\b\|(\Phi,\Psi,W)_\tau\|^2
+C_{\b}~\v^2\int\int\frac{\nu(|\xi|)}{M_*}|\tilde{G}|^2d\xi
dy+C\v^2\int\int\frac{\nu(|\xi|)}{M_*}|\tilde{G}_y|^2d\xi
dy\\[2mm]
\di \quad~~~~~~~ +C(\delta+\l_0)\sum_{|\alpha|=1}\|\partial^\alpha(\phi,\psi,\zeta)\|^2
+C\delta\v^{3}(1+t)^{-\frac32},
\end{array}
\end{equation*}
which implies by \eqref{B4.48} and \eqref{B4.49} that
\begin{equation}\label{B4.67}
\begin{array}{l}
\di I_1^{3}\leq
-\Big(\v\int\bar{v}\lambda_3b_1\int\xi_1^2L_M^{-1}\tilde{G}d\xi
dy\Big)_\tau+C\delta\v^2(1+t)^{-1}E_1+C(\l_0+\b)(K_1+\|\Phi_y\|^2)\\[2mm]
\di ~~~~~~~~~~+C\b\|(\Phi,\Psi,W)_\tau\|^2
+C_{\b}~\v^2\int\int\frac{\nu(|\xi|)}{M_*}|\tilde{G}|^2d\xi
dy+C\v^2\int\int\frac{\nu(|\xi|)}{M_*}|\tilde{G}_y|^2d\xi
dy\\[2mm]
\di \quad~~~~~~~ +C(\delta+\l_0)\sum_{|\alpha|=1}\|\partial^\alpha(\phi,\psi,\zeta)\|^2
+C\delta\v^{2}(1+t)^{-1}.
\end{array}
\end{equation}
And finally, \eqref{B4.45}, \eqref{B4.47} and \eqref{B4.67} yield the estimate
on $I_1$ as follows.
\begin{equation}\label{B4.68}
\begin{array}{l}
\di I_1\leq
-\Big(\v\int\bar{v}\lambda_3b_1\int\xi_1^2L_M^{-1}\tilde{G}d\xi
dy\Big)_\tau+C\delta\v^2(1+t)^{-1}E_1+C(\l_0+\b)(K_1+\|\Phi_y\|^2)\\[2mm]
\di ~~~~~~~~~~+C\b \|(\Phi,\Psi,W)_\tau\|^2
+C_{\b}~\v^2\int\int\frac{\nu(|\xi|)}{M_*}|\tilde{G}|^2d\xi
dy+C\v^2\int\int\frac{\nu(|\xi|)}{M_*}|\tilde{G}_y|^2d\xi
dy\\[2mm]
\di \quad~~~~~~~ +C(\delta+\l_0)\sum_{|\alpha|=1}\|\partial^\alpha(\phi,\psi,\zeta)\|^2
+C\delta\v^{2}(1+t)^{-1}+\f1{16}\int|\hat\t_y|b_1^2dy,
\end{array}
\end{equation}
which completes the estimate on the term
$\di \int\bar{v}Q_1\Psi_1dx$.

\

\noindent{\bf\underline{ Estimate on $\di \int q_2b_2dy$}:}

\

Notice that the profile has no intrinsic dissipation on $b_2$. Fortunately, it holds that
$q_2=\sqrt{\frac25}\frac{Q_4}{p_+}$ and $Q_4$ has the decay rate as $\v^3(1+t)^{-\frac32}$. Thus the estimation on $\di \int q_2b_2dy$
can be directly obtained even though there is no intrinsic
dissipation on $b_2$. For
example,
\begin{equation*}\label{B4.69}
\begin{array}{l}
\di |\int\v\bar{u}_1\bar R_1b_2dy|\le
C\delta\v^2(1+t)^{-1}E_1+C\delta\v^3(1+t)^{-\frac32},\\[2mm]
\di |\int\int \v^2\bar{u}_1b_2\xi_1^2\Theta_1d\xi
dy|\leq C\delta\v^2(1+t)^{-1}E_1
+C(\delta+\l_0)\sum_{|\alpha|=1}\|\partial^\alpha(\phi,\psi,\zeta)\|^2\\[2mm]
\di \qquad~~~~~~~~~~~~~~~ +C\d\v^2\int\int\frac{\nu(|\xi|)\tilde{G}^2}{M_*}d\xi
dy +C\int\int\frac{\nu(|\xi|)}{M_*}(\tilde G_\tau^2+\tilde G_y^2)d\xi dy+C\delta\v^3(1+t)^{-\frac32}.
\end{array}
\end{equation*}
And the term $\di \v\int\int\xi_1|\xi|^2\Theta_1^2
b_2d\xi dy$ can be estimated similarly as  for
$I_1^{32}$ where the intrinsic dissipation on $b_1, b_3$ is not needed. Notice also that all the other terms in $q_2$ are of
higher order. Therefore, one has
\begin{align}\label{B4.70}
& I_2=\int q_2b_2d\xi dy\leq (\v\int\int
\hat{A}(\xi,b_2)L_M^{-1}\tilde{G}d\xi
dy)_\tau+C\delta\v^2(1+t)^{-1}E_1+C(\l_0+\b)(K_1+\|\Phi_y\|^2)\nonumber\\
& ~~~~~~~~~~+C\b\|(\Phi,\Psi,W)_\tau\|^2
+C_{\b}~\v^2\int\int\frac{\nu(|\xi|)}{M_*}|\tilde{G}|^2d\xi
dy+C\v^2\int\int\frac{\nu(|\xi|)}{M_*}|\tilde{G}_y|^2d\xi
dy\nonumber\\
& \quad~~~~~~~ +C(\delta+\l_0)\sum_{|\alpha|=1}\|\partial^\alpha(\phi,\psi,\zeta)\|^2
+C\delta\v^{2}(1+t)^{-1},
\end{align}
where $\hat{A}(\xi, b_2)$ is a linear function of $b_2$ and a polynomial function of
$\xi$. Using \eqref{B4.68}, \eqref{B4.70} and \eqref{B4.38}, we get
\begin{equation}\label{B4.71}
\begin{array}{l}
\di E_{1\tau}+\Big(\int\int\v\hat{A}_1(\xi,B)L_M^{-1}\tilde{G}d\xi
dy\Big)_\tau+\frac14
K_1+\int|\hat\t_{y}|(b_1^2+b_3^2)dy\\[3mm]
\di ~~~~~\leq C\delta\v^2(1+t)^{-1}E_1+C(\l_0+\b)(K_1+\|\Phi_y\|^2)+C\b\|(\Phi,\Psi,W)_\tau\|^2\\[2mm]
\di ~~~~~~~~~~~~+C_{\b}\v^2\int\int\frac{\nu(|\xi|)}{M_*}|\tilde{G}|^2d\xi
dy+C\v^2\int\int\frac{\nu(|\xi|)}{M_*}|\tilde{G}_y|^2d\xi
dy\\[2mm]
\di \quad~~~~~~~~~~ +C(\delta+\l_0)\sum_{|\alpha|=1}\|\partial^\alpha(\phi,\psi,\zeta)\|^2
+C\delta\v^{2}(1+t)^{-1},
\end{array}
\end{equation}
where we have used the smallness of $\delta$ and $\varepsilon_0$.
Here $\hat{A}_1$ is a linear function of $B=(b_1,b_2,b_3)^t$ and a polynomial function of
$\xi$.

\

Note that $K_1$ does not contain the norm $\|\Phi_y\|^2$. To
complete the lower order inequality, we have to estimate $\Phi_y$.
From $\eqref{B4.9}_2$, we have
\begin{equation}\label{B4.72}
\frac{4\mu(\bar{\theta})}{3\bar{v}}\Phi_{y\tau}-\Psi_{1\tau}+\frac{p_+}{\bar{v}}\Phi_y=\frac2{3\bar{v}}W_y-
\frac{8\mu(\bar{\theta})}{15p_+\bar{v}}\v^2\hat{N}_{1y}-Q_1.
\end{equation}
Multiplying \eqref{B4.72} by $\Phi_y$ yields
\begin{equation*}\label{B4.73} (\frac{2\mu(\bar{\theta})}{3\bar{v}}\Phi_y^2)_\tau-(\frac{2\mu(\bar{\theta})}{3\bar{v}})_\tau\Phi_y^2
-\Phi_y\Psi_{1\tau}+\frac{p_+}{\bar{v}}\Phi_y^2
=\Big(\frac2{3\bar{v}}W_y-\frac{8\mu(\bar{\theta})}{15p_+\bar{v}}\v^2\hat{N}_{1y}-Q_1\Big)\Phi_y.
\end{equation*}
Since
\begin{equation*}\label{B4.74}
\Phi_y\Psi_{1\tau}=(\Phi_y\Psi_1)_\tau-(\Phi_\tau\Psi_1)_y+\Psi_{1y}^2-\frac2{5p_+}\v^2\hat{N}_1\Psi_{1y},
\end{equation*}
we can obtain
\begin{equation}\label{B4.75}
\Big(\int\frac{2\mu(\bar{\theta})}{3\bar{v}}\Phi_y^2-\Phi_y\Psi_1dy\Big)_\tau
+\int\frac{p_+}{2\bar{v}}\Phi_y^2dy\\[2mm]
\leq C\|(\Psi_{1y},W_y)\|^2+C\delta\v^3(1+t)^{-3/2}+\int Q_1^2dy.
\end{equation}
The formula \eqref{B4.10} for $Q_1$
 and the Cauchy inequality directly yield
\begin{equation}\label{B4.76}
\begin{array}{l}
\di \int Q_1^2dy\leq
C(\d+\l_0)(K_1+\|\Phi_y\|^2)+
C\l_0\sum_{|\alpha|=1}\|\partial^\alpha(\phi,\psi,\zeta)\|^2\\[2mm]
\di ~~~~~~~~~~~~~~~~~~+C\delta\v^3(1+t)^{-3/2}+C\int|\int\xi_1^2\Theta_1d\xi|^2 dy.
\end{array}
\end{equation}
And using Lemmas \ref{lem4.1}-\ref{lem4.3} implies
\begin{align}\label{B4.77}
&\int|\int\xi_1^2\Theta_1d\xi|^2d\xi dy\leq
C\v^2\int\int\frac{\nu(|\xi|)}{M_*}(|\tilde{G}_y|^2+|\tilde{G}_\tau|^2)d\xi dy+C(\delta+\l_0)\v^4\int\int\frac{\nu(|\xi|)|\tilde{G}|^2}{M_*}d\xi dy\nonumber\\
&~~~~~~~~~~~~~~~~~~~~~~~~~~~~~+C(\delta+\l_0)\sum_{|\alpha|=1}\|\partial^\alpha(\phi,\psi,\zeta)\|^2
+C\delta\v^3(1+t)^{-3/2}.
\end{align}
Substituting  \eqref{B4.76} and \eqref{B4.77} into \eqref{B4.75} yields
\begin{equation}\label{B4.78}
\begin{array}{l}
\di \Big(\int\frac{2\mu(\bar{\theta})}{3\bar{v}}\Phi_y^2-\Phi_y\Psi_1dy\Big)_\tau
+\int\frac{p_+}{4\bar{v}}\Phi_y^2dy\\
\di \leq
C_2K_1+C_2\v^2\int\int\frac{\nu(|\xi|)}{M_*}(|\tilde{G}_y|^2+|\tilde{G}_\tau|^2)d\xi dy+C_2(\delta+\l_0)\v^4\int\int\frac{\nu(|\xi|)|\tilde{G}|^2}{M_*}d\xi dy\\[2mm]
\di ~~~~~~~~~~~~~~~~+C_2(\delta+\l_0)\sum_{|\alpha|=1}\|\partial^\alpha(\phi,\psi,\zeta)\|^2
+C_2\delta\v^3(1+t)^{-3/2}.
\end{array}
\end{equation}
Multiplying \eqref{B1.22-1} by $\v^2\frac{\tilde{G}}{M_*}$, one can obtain
\begin{eqnarray}\label{B4.80}
\begin{array}{ll}
\di (\v^2\frac{v\tilde{G}^2}{2M_*})_\tau-\v^2\frac{v\tilde{G}}{M_*}L_M\tilde{G}=\Big\{-\frac1{R\theta}P_1[\xi_1
(\frac{|\xi-\v u|^2}{2\theta}\f1\v\zeta_y+\xi\cdot\f1\v\psi_y)M]\\[2mm]
\di ~~~~~~~~~~~~~~~~~~~~~~~~+\v u_1G_y-P_1(\xi_1G_y)+\v vQ(G,G)-v\bar{G}_\tau\Big\}\cdot
\v^2\frac{\tilde{G}}{M_*}.
\end{array}
\end{eqnarray}
Integrating \eqref{B4.80} with respect to $\xi$ and $y$ and using the
Cauchy inequality and Lemmas \ref{lem4.1}-\ref{lem4.3}, one has
\begin{eqnarray}\label{B4.81}
\begin{array}{l}
\di \Big(\v^2\int\int\frac{1}{2M_*}|\tilde{G}|^2d\xi
dy\Big)_\tau+\frac{3\bar{\sigma}}{4}\v^2\int\int\frac{\nu(|\xi|)}{M_*}|\tilde{G}|^2d\xi dy\\[3mm]
\di \leq
C_3\delta\v^3(1+t)^{-3/2}+C_3\sum_{|\alpha|=1}\|\partial^\alpha(\phi,\psi,\zeta)\|^2
+C_3~\v^2\int\int\frac{\nu(|\xi|)}{M_*}|\tilde{G}_y|^2d\xi dy.
 \end{array}
\end{eqnarray}
On the other hand, since $(\Phi,\Psi,W)_\tau$ can be represented by
$(\Phi,\Psi,W)_y$ and $(\Phi,\Psi,W)_{yy}$ from the equation
\eqref{B4.9}, we can get an estimate for $(\Phi,\Psi,W)_\tau$ as
follows.
\begin{align}\label{B4.82}
& \|(\Phi,\Psi,W)_\tau\|^2\leq
C_4(K_1+\|\Phi_y\|^2)+C_4\sum_{|\alpha|=1}\|\partial^\alpha(\phi,\psi,\zeta)\|^2
+C_4\delta\v^3(1+t)^{-\frac32}\nonumber\\
&\quad~~~~~~~ +C_4\v^2\int\int\frac{\nu(|\xi|)}{M_*}(|\tilde{G}_y|^2+|\tilde{G}_\tau|^2)d\xi dy
+C_4(\d+\l_0)\v^4\int\int\frac{\nu(|\xi|)|\tilde{G}|^2}{M_*}d\xi dy.
\end{align}
Now we can complete the lower order estimate. Since $\hat{A}_1$ is a
linear function of the vector $B$ and a polynomial of $\xi$, we
get
\begin{equation*}\label{B4.83}
|\v\int\int\hat{A}_1(\xi,B)L_M^{-1}\tilde{G}d\xi dy|\le
\frac14E_1+C\v^2\int\int\frac{|\tilde{G}|^2}{M_*}d\xi dy.
\end{equation*}
We choose large constants $\bar{C}_1>1$, $\bar{C}_2>1$, $\bar{C}_3>1$ and
small constant $\b$ such that
\begin{equation*}\label{B4.83-1}
\begin{array}{l}
\di \bar{C}_1E_1+\bar{C}_1\v\int\int\hat{A}_1L_M^{-1}\tilde{G}d\xi
dy+\bar{C}_2\int\big(\frac{2\mu(\bar{\theta})}{3\bar{v}}\Phi_y^2-\Phi_y\Psi_1\big)dy
+\bar{C}_3\v^2\int\int\frac{|\tilde{G}|^2}{2M_*}d\xi dy\\[3mm]
\di ~~~~\geq\frac12\bar{C}_1E_1+\bar{C}_2\int\frac{\mu(\bar{\theta})}{3\bar{v}}\Phi_y^2dy
+\frac{\bar{C}_3}{4}\v^2\int\int\frac{\tilde{G}^2}{M_*}d\xi dy,
\end{array}
\end{equation*}
\begin{equation*}\label{B4.84}
\begin{array}{l}
\di \Big(\frac{\bar{C}_1}{4}-C_2\bar{C}_2-\bar{C}_1C_1C_4\Big)K_1
+\int\Big[\bar{C}_2\frac{p_+}{4\bar{v}}-(\d+\b+\l_0)\bar{C}_1C_1(1+C_4)\Big]\Phi_y^2dy\\[3mm]
\di \geq  \frac{\bar{C}_1}{8}K_1+
\bar{C}_2\int\frac{p_+}{8\bar{v}}\Phi_y^2dy,
\end{array}
\end{equation*}
and
\begin{equation*}\label{B4.85}
\frac{\bar{\sigma}}{2}\bar{C}_3-\bar{C}_1C_1C_4(\d+\b+\l_0)
-C_{\b}\bar{C}_1-C_2\bar{C}_2\v^2(\d+\l_0)\geq
\frac{\bar{\sigma}}{4}\bar{C}_3.
\end{equation*}
Hence, by multiplying \eqref{B4.71} by $\bar{C}_1$, \eqref{B4.78} by $\bar{C}_2$, \eqref{B4.81} by
$\bar{C}_3$, \eqref{B4.82} by
$C_1(\delta+\varepsilon_0+\varepsilon_1)\bar{C}_1$ and adding all
these inequalities together, we have

\begin{eqnarray}\label{B4.87}
&&E_{2\tau}+K_2+\int|\hat\theta_y|(b_1^2+b_3^2)dy\leq C_5\delta\v^2(1+t)^{-1}(E_2+1)\nonumber\\
&&~~~~~~~~~~+C_5\v^2\int\int\frac{\nu(|\xi|)}{M_*}(|\tilde{G}_y|^2+|\tilde{G}_\tau|^2)d\xi dy +C_5\sum_{|\alpha|=1}\|\partial^\alpha(\phi,\psi,\zeta)\|^2,
\end{eqnarray}
where
\begin{eqnarray*}
&&E_2=\bar{C}_1E_1+\bar{C}_1\int\int\v\hat{A}_1L_M^{-1}\tilde{G}d\xi
dy+\bar{C}_2\int\big(\frac{2\mu(\bar{\theta})}{3\bar{v}}\Phi_y^2-\Phi_y\Psi_1\big)dy+
\bar{C}_3\v^2\int\int\frac{|\tilde{G}|^2}{2M_*}d\xi dy,\label{B4.88}\\
&&K_2=\frac{\bar{C}_1}{16}K_1+\bar{C}_2
\int\frac{p_+}{16\bar{v}}\Phi_y^2dy+\|(\Phi,\Psi,W)_\tau\|^2
+\frac{\bar{\sigma}}{8}\bar{C}_3\v^2\int\int\frac{\nu(|\xi|)|\tilde{G}|^2}{M_*}d\xi
dy.\label{B4.89}
\end{eqnarray*}

\subsection{Derivative estimate}

In this subsection, we derive the higher order estimate for $(\Phi,\Psi,W)$. First, by the definition of $\Theta_1^1$ in \eqref{B1.23}, it holds that
\begin{eqnarray}\label{B4.91}
\begin{cases}
\di -\v\int\xi_1^2\Theta_{1}^1d\xi=\v^2N_{4}+\v^3F_{4},\\
\di N_{4}=f_{41}\theta_x{\bar\theta}_x+f_{42}v_x\bar{\theta}_x+f_{43}\bar{\theta}_x^2
+f_{44}\bar{\theta}_{xx},\\
\di |F_{4}|=O(1)\big[(|v_x|+|\theta_x|+|\bar{\theta}_x|+|\v u_x|+|\v \bar{u}_x|)|\bar{u}_x|
+|u_x||\bar{\theta}_x|+|{\bar u}_{xx}|\big].
\end{cases}
\end{eqnarray}
From \eqref{B1.15}
and \eqref{B1.41}, one has
\begin{eqnarray}\label{B4.90}
\left\{
\begin{array}{l}
\di \phi_\tau-\psi_{1y}=-\frac{3}{5}\v^2\hat{N}_{1y},\\[2mm]
\di \psi_{1\tau}+(p-\bar{p})_y=
\frac{4\v}3(\frac{\mu(\theta)}{v}u_{1y}-\frac{\mu(\bar{\theta})}{\bar{v}}\bar{u}_{1y})_y+Q_5,\\[2mm]
\di \psi_{i\tau}=\v(\frac{\mu(\theta)}{v}u_{iy}-\frac{\mu(\bar{\theta})}{\bar{v}}\bar{u}_{iy})_y
+Q_{4+i},i=2,3,\\[2mm]
\di \zeta_\tau+\v pu_{1y}-\v\bar{p}\bar{u}_{1y}=(\frac{\k(\theta)}{v}\theta_{y}
-\frac{\k(\bar{\theta})}{\bar{v}}\bar{\theta}_{y})_y+Q_8,
\end{array}
\right.
\end{eqnarray}
where
\begin{equation*}\label{B4.92}
\begin{cases}
\di Q_5=\v^2(N_{4y}-\bar{N}_{4y})+\v^3F_{4y}-\v\int\xi_1^2\Theta_{1y}^2d\xi+(\v^2\bar{N}_{4y}-\bar{R}_{1y}),\\
\di Q_{4+i}=\v^2(N_{iy}-\bar{N}_{iy})+\v^3F_{iy}-\v\int\xi_1\xi_{i}\Theta_{1y}^2d\xi-\bar{R}_{iy},i=2,3,\\
\di Q_8=\frac2{5}\v^2\hat{N}_{1y}
+\v^2(N_{1y}-\bar{N}_{1y})+\v^3F_{1y}+\frac43\frac{\mu(\theta)}{v}\v^2u_{1y}^2
+\sum_{i=2}^3\frac{\mu(\theta)}{v}\v^2u_{iy}^2\\
\di ~~~~~~~-\frac12\v\int\xi_1|\xi|^2\Theta_{1y}^2d\xi
+\sum_{i=1}^3\v^2u_i\int\xi_1\xi_i\Theta_{1y}d\xi
-\bar{H}_{y}-\bar{R}_{4y}+\frac12(|\v\bar{u}|^2)_\tau+\v\bar{p}_y\bar{u}_1,
\end{cases}
\end{equation*}
and $N_i, F_i ~(i=1,2,3,4)$ are defined in \eqref{B1.29}, \eqref{B1.34} and \eqref{B4.91} respectively and $\bar{N}_{i}~(i=1,2,3,4)$ is the corresponding function of $N_i~(i=1,2,3,4)$ by substituting the
variable $(v,u,\theta)$ by the  profile
$(\bar{v},\bar{u},\bar{\theta})$.

We will use the convex entropy for the fluid system to obtain the first-order derivative estimates of
$(\Phi_y,\Psi_y,W_y)$. Multiplying $\eqref{B4.90}_2$ by $\psi_1$ and $\eqref{B4.90}_3$ by $\psi_i$, one has
\begin{align*}\label{4.57}
(\f12\sum_{i=1}^3\psi_i^2)_\tau-(p-\bar{p})\psi_{1y}+\frac{4\v}3\Big(\frac{\mu(\theta)}{v}u_{1y}
-\frac{\mu(\bar{\theta})}{\bar{v}}\bar{u}_{1y}\Big)\psi_{1y}
+\v(\frac{\mu(\theta)}{v}u_{iy}-\frac{\mu(\bar{\theta})}{\bar{v}}\bar{u}_{iy})\psi_{iy}=\sum_{i=1}^3Q_{4+i}\psi_i+(\cdots)_y.
\end{align*}
Since $p-\bar{p}=\f23\bar{\theta}(\frac1v-\frac1{\bar
v})+\frac{2\zeta}{3v}$,  we obtain
\begin{eqnarray} \label{4.58}
&&(\frac12\sum_{i=1}^3\psi_i^2)_\tau-\f23\bar{\theta}(\frac{1}{v}-\frac{1}{\bar{v}})\phi_\tau
-\frac{2}{3v}\zeta\psi_{1y}+\frac{4}{3}\frac{\mu(\theta)}{v}\psi_{1y}^2+\sum_{i=2}^3\frac{\mu(\theta)}{v}\psi_{iy}^2
+\frac{4\v}3(\frac{\mu(\theta)}{v}-\frac{\mu(\bar{\theta})}{\bar{v}})
\bar{u}_{1y}\psi_{1y}\nonumber\\
&&+\v\sum_{i=2}^3(\frac{\mu(\theta)}{v}-\frac{\mu(\bar{\theta})}{\bar{v}})\bar{u}_{iy}\psi_{iy}
=\sum_{i=1}^3\psi_i Q_{i+4}+\frac{2\v^2}{5}\bar{\theta}(\frac{1}{v}-\frac{1}{\bar{v}})\hat{N}_{1y}
+(\cdots)_x.
\end{eqnarray}
Let
\begin{eqnarray*}\label{4.59}
\hat{\Phi}(s)=s-1-\ln s,
\end{eqnarray*}
then it holds that
\begin{eqnarray}\label{4.60}
&&\{\f23\bar{\theta}\hat{\Phi}(\frac{v}{\bar{v}})\}_\tau=\f23\bar{\theta}_\tau\hat{\Phi}(\frac{v}{\bar{v}})
+\f23\bar{\theta}(-\frac{1}{v}+\frac{1}{\bar{v}})\phi_\tau
 +\f23\bar{\theta}(-\frac{v}{\bar{v}^2}+\frac{1}{ \bar{v}}) \bar{v}_t
+\f23\bar{\theta}(-\frac{1}{v}+\frac{1}{\bar{v}})\bar{v}_\tau\nonumber\\
&&=\f23\bar{\theta}(-\frac{1}{v}+\frac{1}{\bar{v}})\phi_\tau-\bar{p}\hat{\Psi}(\frac{v}{\bar{v}})\bar{v}_\tau
+\bar{v}\bar{p}_\tau\hat{\Phi}(\frac{v}{\bar{v}}),
\end{eqnarray}
where
\begin{eqnarray*}\label{4.61}
\hat{\Psi}(s)=s^{-1}-1+\ln{s}.
\end{eqnarray*}
It is easy to check that $\hat{\Phi}(1)=\hat{\Phi}'(1)=\hat{\Psi}(1)=\hat{\Psi}'(1)=0$ and $\hat{\Phi}(s)$, $\hat{\Psi}(s)$ are
strictly convex around $s=1$.
Substituting \eqref{4.60} into \eqref{4.58} yields that
\begin{eqnarray} \label{4.62}
&&\Big(\frac12\sum_{i=1}^3\psi_i^2+\f23\bar{\theta}\hat{\Phi}(\frac{v}{\bar{v}})\Big)_\tau
-\frac{2}{3v}\zeta\psi_{1y}+\frac{4}{3}\frac{\mu(\theta)}{v}\psi_{1y}^2+\sum_{i=2}^3\frac{\mu(\theta)}{v}\psi_{iy}^2
+\frac{4\v}3(\frac{\mu(\theta)}{v}-\frac{\mu(\bar{\theta})}{\bar{v}})
\bar{u}_{1y}\psi_{1y}\nonumber\\
&&+\v\sum_{i=2}^3(\frac{\mu(\theta)}{v}-\frac{\mu(\bar{\theta})}{\bar{v}})\bar{u}_{iy}\psi_{iy}
=\sum_{i=1}^3\psi_i Q_{i+4}+\frac{2\v^2}{5}\bar{\theta}(\frac{1}{v}-\frac{1}{\bar{v}})\hat{N}_{1y}
-\bar{p}\hat{\Psi}(\frac{v}{\bar{v}})\bar{v}_\tau
+\bar{v}\bar{p}_\tau\hat{\Phi}(\frac{v}{\bar{v}}).
\end{eqnarray}
On the other hand, multiplying $\eqref{B4.90}_4$ by $\f{\zeta}{\t}$, it holds that
\begin{eqnarray}\label{N2.63}
\f{\zeta}{\t}\zeta_\tau+\v(pu_{1y}-\bar{p}\bar{u}_{1y})\f{\zeta}{\t}
=\Big(\f{\k(\t)}{v}\t_y-\f{\k(\bar\t)}{\bar v}\bar{\t}_y\Big)_y\f{\zeta}{\t}+Q_8\f{\zeta}{\t}.
\end{eqnarray}
One can compute that
\begin{equation}\label{N2.64}
\f{\zeta}{\t}\zeta_\tau=\Big(\bar{\t}\hat{\Phi}(\frac{\t}{\bar{\t}})\Big)_{\tau}
+\bar{\t}_\tau\hat\Psi(\f{\t}{\bar\t})=\Big(\bar{\t}\hat{\Phi}(\frac{\t}{\bar{\t}})\Big)_{\tau}
+O(1)\d\v^2(1+t)^{-1}|\zeta|^2,
\end{equation}
\begin{equation}\label{N2.65}
\v(pu_{1y}-\bar{p}\bar{u}_{1y})\f{\zeta}{\t}=\frac{2\zeta}{3v}\psi_{1y}+\v(p-\bar p)\bar{u}_{1y}\f{\zeta}{\t}=\frac{2\zeta}{3v}\psi_{1y}+O(1)\d\v^2(1+t)^{-1}|(\psi,\zeta)|^2,
\end{equation}
and
\begin{equation}\label{N2.66}
\begin{array}{ll}
\di \Big(\f{\k(\t)}{v}\t_y-\f{\k(\bar\t)}{\bar v}\bar{\t}_y\Big)_y\f{\zeta}{\t}&\di =(\cdots)_y-\f{\bar\t\k(\t)}{v\t^2}\zeta_y^2-\f{\k(\t)\bar\t_y\z\z_y}{v\t^2}-\f{\bar\t\bar\t_y\z_y-|\bar\t_y|^2\z}{\t^2}\big(\f{\k(\t)}{v}-\f{\k(\bar\t)}{\bar v}\big)\\[3mm]
&\di \leq (\cdots)_y-\f34\f{\bar\t\k(\t)}{v\t^2}\zeta_y^2+C\d\v^2(1+t)^{-1}|(\phi,\zeta)|^2.
\end{array}
\end{equation}
Substituting \eqref{N2.64}, \eqref{N2.65} and \eqref{N2.66} into \eqref{N2.63} yields that
\begin{eqnarray}\label{N2.67}
\Big(\bar{\t}\hat{\Phi}(\frac{\t}{\bar{\t}})\Big)_{\tau}+\frac{2\zeta}{3v}\psi_{1y}+\f34\f{\bar\t\k(\t)}{v\t^2}\zeta_y^2
\leq  (\cdots)_y+C\d\v^2(1+t)^{-1}|(\phi,\zeta)|^2+|Q_8\f{\zeta}{\t}|.
\end{eqnarray}
Combining \eqref{N2.67} and \eqref{4.62} and using Cauchy inequality, one has
\begin{equation}\label{B4.93}
\begin{array}{ll}
\di E_{3\tau}+\frac34K_3\leq &\di
C\delta\v^2(1+t)^{-1}E_3+C\d\v^4(1+t)^{-2}+C\v^2(1+t)^{-1}\int|\hat\t_y|(b_1^2+b_3^2)dy\\
&\di
+\left|\int\f{2}{5}\v^2\hat{N}_{1y}\Big(-\f{\bar\t}{v\bar{v}}\phi+\f{\zeta}{\t}\Big)dy\right|
+\sum_{i=1}^4|I_{i+2}|,
\end{array}
\end{equation}
where
\begin{eqnarray}\label{B4.95}
\begin{cases}
E_3=\int\Big(\frac12\sum_{i=1}^3\psi_i^2+R\bar{\theta}\hat{\Phi}(\frac{v}{\bar{v}})+\bar{\theta}
\hat{\Phi}(\frac{\theta}{\bar{\theta}})\Big)dy, \\[3mm]
K_3=\int\Big(\frac{4}{3}\frac{\mu(\theta)}{v}\psi_{1y}^2
+\sum_{i=2}^3\frac{\mu(\theta)}{v}\psi_{iy}^2+\f34\frac{\bar\t\k(\theta)}{v\theta^2}\zeta_y^2\Big)dy,
\end{cases}
\end{eqnarray}
and
\begin{eqnarray*}
&& I_3=\int\v^2(N_{4y}-\bar{N}_{4y})\psi_1 dy+\int\v^2F_{4y}\psi_1 dy-\v\int\int\xi_1^2\Theta_{1y}^2\psi_1 d\xi dy\label{B4.96}\\
&& I_{2+i}=\int\v^2(N_{iy}-\bar{N}_{iy})\psi_i dy+\int\v^3F_{iy}\psi_i dy-\v\int\int\xi_1\xi_i\Theta_{1y}^2\psi_id\xi dy,~~i=2,3,\label{B4.97}\\
&& I_6=\int\v^2(N_{1y}-\bar{N}_{1y})\f{\zeta}{\t}dy+\int\v^3F_{1y}\f{\zeta}{\t}dy
-\frac12\v\int\int\xi_1|\xi|^2\Theta_{1y}^2\f{\zeta}{\t}d\xi dy\nonumber\\
&&~~~~~~~~~~~~~~~~~~~~~~+\sum_{i=1}^3\v^2 u_i\int\int\xi_1\xi_i\Theta_{1y}\f{\zeta}{\t}d\xi dy.\label{B4.98}
\end{eqnarray*}
In the estimate of \eqref{B4.93}, we have used the estimate like
\begin{eqnarray*}\label{B4.99-1}
&&|\int (\v^2\bar{N}_{4y}-\bar{R}_{1y})\psi_{1}dy|=|\int (\v^2\bar{N}_{4y}-\bar{R}_{1y})\Psi_{1y}dy|
\leq C\int |(\v^2\bar{N}_{4yy}-\bar{R}_{1yy})|(|b_1|+|b_3|)dy\nonumber\\
&&\leq C\d\v^4(1+t)^{-2}+C\v^2(1+t)^{-1}\int|\hat\t_y|(b_1^2+b_3^2)dy.
\end{eqnarray*}
Now, we calculate the terms on the right hand side of \eqref{B4.93}. Firstly,  a direct calculation yields
\begin{equation*}\label{B4.99}
-\f{\bar\t}{v\bar{v}}\phi+\f{\zeta}{\t}=\f{1}{\bar\t}(\f{2}{3p_+}W_y-\Phi_y)+O(1)\Big[|(\phi,\psi,\zeta)|^2+|\v\bar{u}_{1y}\Psi_1|+\d\v(1+t)^{-1}\Big],
\end{equation*}
and thus by combining \eqref{B4.40} and \eqref{B4.42}, it holds that
\begin{align}\label{B4.100}
&|\int\f{2}{5}\v^2\hat{N}_{1y}\Big(-\f{\bar\t}{v\bar{v}}\phi+\f{\zeta}{\t}\Big)dy|\nonumber\\
&\leq C\delta\v^2(1+t)^{-1}E_3+C\d\v^4(1+t)^{-2}
+C\int\v^2[|(\hat{N}_{1y}\f{1}{\bar\t})_y|+|\hat{N}_{1y}\v\bar{u}_{1y}|]\cdot|(b_1,b_3)|dy \nonumber\\
&\leq C\delta\v^2(1+t)^{-1}E_3+C\d\v^4(1+t)^{-2}+C\v^2(1+t)^{-1}\int|\hat\t_y|(b_1^2+b_3^2)dy.
\end{align}
By using the definition of $N_4$ and $F_4$ in \eqref{B4.91}, one can obtain
\begin{eqnarray}\label{B4.101}
&&|\int\v^2(N_{4y}-\bar{N}_{4y})\psi_1 dy+\int\v^3F_{4y}\psi_1 dy|\leq \f1{32}K_3+\int\v^4(N_{4}-\bar{N}_{4})^2 dy+\v^6F_{4}^2dy\nonumber\\
&&\leq \f1{32}K_3+C\d\|\phi_y\|^2+C\delta\v^2(1+t)^{-1}E_3+C\d\v^4(1+t)^{-2}.
\end{eqnarray}
And by using Lemma \ref{lem4.1}-\ref{lem4.3}, one has that
\begin{eqnarray}\label{B4.102}
&&|\v\int\int\xi_1^2\Theta_{1y}^2\psi_1 d\xi dy|\leq \f1{32}K_3+C\v^2\int|\int\xi_1^2\Theta_{1}^2d\xi|^2dy\nonumber\\
&&\leq \f1{32}K_3+C\v^2\int\int\frac{\nu(|\xi|)}{M_*}|(\tilde{G}_y,\tilde{G}_\tau)|^2d\xi dy + C\d\sum_{|\alpha|=1}\|\partial^\alpha(\phi,\psi,\zeta)\|^2
+C\delta\v^2(1+t)^{-1}E_3\nonumber\\
&&~~~+C\d\v^4(1+t)^{-2} +C\v^4\Big[\d(1+t)^{-1}+\|\int\frac{|\tilde{G}|^2}{M_*}d\xi\|_{L^\infty}\Big]\int\int\frac{\nu(|\xi|)|\tilde{G}|^2}{M_*}d\xi dy.
\end{eqnarray}
Combining \eqref{B4.101} and \eqref{B4.102} yields that
\begin{eqnarray}\label{B4.103}
&&I_3\leq \f1{32}K_3+C\v^2\int\int\frac{\nu(|\xi|)}{M_*}|(\tilde{G}_y,\tilde{G}_\tau)|^2d\xi dy + C\d\sum_{|\alpha|=1}\|\partial^\alpha(\phi,\psi,\zeta)\|^2
+C\delta\v^2(1+t)^{-1}E_3\nonumber\\
&&~~~~~~~+C\d\v^4(1+t)^{-2} +C\v^4\Big[\d(1+t)^{-1}+\|\int\frac{|\tilde{G}|^2}{M_*}d\xi\|_{L^\infty}\Big]\int\int\frac{\nu(|\xi|)|\tilde{G}|^2}{M_*}d\xi dy.
\end{eqnarray}
Similarly, $I_4, I_5, I_6$ can be controlled by the right hand side of \eqref{B4.103}. Substituting \eqref{B4.100} and \eqref{B4.103} into  \eqref{B4.93} gives that
\begin{eqnarray}\label{B4.93-1}
&&E_{3\tau}+\frac12K_3\leq
 C_6\v^2\int\int\frac{\nu(|\xi|)}{M_*}|(\tilde{G}_y,\tilde{G}_\tau)|^2d\xi dy+C_6\d\sum_{|\a|=1}\|\partial_\a(\phi,\psi,\zeta)\|^2\nonumber\\
&&~~~~~~~+C_6\delta\v^2(1+t)^{-1}E_3+C_6\d\v^4(1+t)^{-2}+C_6\v^2(1+t)^{-1}\int|\hat\t_y|(b_1^2+b_3^2)dy\nonumber\\
&&
~~~~~~~~~~~
+C_6\v^4\Big[\d(1+t)^{-1}+\|\int\frac{|\tilde{G}|^2}{M_*}d\xi\|_{L^\infty}\Big]\int\int\frac{\nu(|\xi|)|\tilde{G}|^2}{M_*}d\xi dy.
\end{eqnarray}

\

Note that the norm $\|\phi_y\|$ is not included in $K_3$ (see
\eqref{B4.95}). To complete the first-order derivative estimate, we follow the
same way as to estimate $\Phi_y$ in the previous section.  By using the $\eqref{B4.90}_1$, we can rewrite
the equation $\eqref{B4.90}_2$ as
\begin{eqnarray}\label{B4.105}
&&\frac43\frac{\mu(\bar{\theta})}{\bar{v}}\phi_{y\tau}-\psi_{1\tau}-(p-\bar{p})_y=-\frac{4\v^2}{5}
\frac{\mu({\bar{\theta}})}{\bar{v}}\hat{N}_{1yy}-\frac{4}{3}
(\frac{\mu({\bar{\theta}})}{\bar{v}})_y\psi_{1y}\nonumber\\
&&~~~~~~~~~~~~~~~~~~~~~~~~~~~~~~~~~~~~~~~-\frac{4\v}3\Big[(\frac{\mu(\theta)}{v}u_{1y}
-\frac{\mu(\bar{\theta})}{\bar{v}}\bar u_{1y})\Big]_y
+\v\int\xi_1^2\Theta_{1y}d\xi+\bar{R}_{1y}.
\end{eqnarray}
 Multiplying the equation \eqref{B4.105} by $\phi_y$, one has
\begin{align}\label{B4.105-1}
&\frac23(\frac{\mu(\bar{\theta})}{\bar{v}}\phi_{y}^2)_\tau-
\frac23(\frac{\mu(\bar{\theta})}{\bar{v}})_\tau\phi_{y}^2-\psi_{1\tau}\phi_y-(p-\bar{p})_y\phi_y\\
& =\Big\{-\frac{4\v^2}{5}
\frac{\mu({\bar{\theta}})}{\bar{v}}\hat{N}_{1yy}-\frac43(\frac{\mu({\bar{\theta}})}{\bar{v}})_y\psi_{1y}
-\frac{4\v}3[(\frac{\mu({\theta})}{v}-\frac{\mu({\bar{\theta}})}{\bar{v}})u_{1y}]_y
+\v\int\xi_1^2\Theta_{1y}d\xi+\bar{R}_{1y}\Big\}\phi_y.\nonumber
\end{align}
Since
\begin{equation*}\label{B4.106}
-(p-\bar{p})_y=\frac{\bar{p}}{\bar{v}}\phi_y-\frac{2}{3\bar{v}}\zeta_y+(\frac{p}{v}
-\frac{\bar{p}}{\bar{v}})v_y-\f23(\frac{1}{v}-\frac{1}{\bar{v}})\theta_y,
\end{equation*}
and
\begin{equation*}\label{B4.107}
\phi_y\psi_{1\tau}=(\phi_y\psi_1)_\tau-(\phi_\tau\psi_1)_y+\psi_{1y}^2-\frac{3\v^2}{5}
\hat{N}_{1y}\psi_{1y},
\end{equation*}
integrating \eqref{B4.105-1} with respect to $y$ and using the Cauchy
inequality yield
\begin{eqnarray}\label{B4.109}
&&\Big(\int\frac{2\mu(\bar{\theta})}{3\bar{v}}\phi_y^2-\phi_y\psi_1dy\Big)_\tau+\int
\frac{\bar{p}}{2\bar{v}}\phi_y^2dy\leq
C_7K_3+C_7\delta\v^2(1+t)^{-1}E_3+C_7\delta\v^5(1+t)^{-\frac52}
\nonumber\\
&&~~~+C_7(\d+\l_0)\v\sum_{|\alpha|=2}\|\partial^\alpha(\phi,\psi,\zeta)\|^2+C_7\d\|\partial_\tau(\phi,\psi,\zeta)\|^2+ C\v^2\int |\int\xi_1^2\Theta_{1y}d\xi|^2dy,
\end{eqnarray}
where we have used the fact that
\begin{equation*}\label{B4.110}
\int|(\frac{p}{v}
-\frac{\bar{p}}{\bar{v}})v_y-\f23(\frac{1}{v}-\frac{1}{\bar{v}})\theta_y||\phi_y|dy\leq
\f18\|\phi_y\|^2+C\delta\v^2(1+t)^{-1}E_3+CK_3.
\end{equation*}
It follows from  \eqref{B1.23} and Lemmas \ref{lem4.1}-\ref{lem4.2} that
\begin{align}\label{B4.111}
&\v^2\int |\int\xi_1^2\Theta_{1y}d\xi|^2dy\leq C\delta\v^5(1+t)^{-\f52}+C\d\v^4(1+t)^{-1}\sum_{|\alpha|=1}\|\partial^\alpha(\phi,\psi,\zeta)\|^2 \nonumber\\
&~~~~~~~~~~~~~~~~~~~~~~~~~~~~~~~~~+C\d\v^2\sum_{|\alpha|=2}\|\partial^\alpha(\phi,\psi,\zeta)\|^2+CJ_3,
\end{align}
where
\begin{align}\label{B4.111-1}
& J_3\doteq\Bigg[\v^2\sum_{|\alpha|=2}\int\int\frac{\nu(|\xi|)}{M_*}|\partial^\alpha\tilde{G}|^2d\xi dy+\v^4(\d+\l_0)\sum_{|\alpha|=1}
\int\int\frac{\nu(|\xi|)}{M_*}|\partial^\alpha\tilde{G}|^2d\xi dy\nonumber\\
& ~~~~~~~~+\v^4\Big(\d\v(1+t)^{-1}+\int\int\frac{|\tilde{G}_y|^2}{M_*}d\xi dy+\int\int\frac{|\tilde{G}_{yy}|^2}{M_*}d\xi dy\Big) \int\int\frac{\nu(|\xi|)|\tilde{G}|^2}{M_*}d\xi dy\Bigg].
\end{align}
To estimate $(\phi,\psi,\zeta)_\tau$, we use \eqref{B4.90} to obtain
\begin{eqnarray}\label{B4.112}
&&\|\partial_\tau(\phi,\psi,\zeta)\|^2\leq C_8(K_3+\|\phi_y\|^2)+C_8\delta\v^2(1+t)^{-1}E_3+C_8\delta\v^5(1+t)^{-\frac52}\nonumber\\
&&~~~~~~~~~~~~~~~~~~~~~~~~~~~~~
+C_8\sum_{|\alpha|=2}\|\partial^\alpha(\phi,\psi,\zeta)\|^2+C_8J_3.
\end{eqnarray}
Thus we choose large constants $\bar{C}_4$ and $\bar{C}_5$
so that
\begin{eqnarray*}\label{B4.113}
\bar{C}_4E_3+\bar{C}_5\int\big(\frac{2\mu(\bar{\theta})}{3\bar{v}}
\phi_y^2-\phi_y\psi_1\big)dy
\ge\frac{\bar{C}_4}{2}E_3+\bar{C}_5\int\frac{\mu(\bar{\theta})}{3\bar{v}}\phi_y^2dy,
\end{eqnarray*}
and
\begin{eqnarray*}\label{B4.114}
\frac12\bar{C_4}-\bar{C}_5C_7-C_8\geq
\frac18\bar{C_4},\quad \bar{C}_5\int\frac{\bar{p}}{2\bar{v}}\phi_y^2dy
-C_8\|\phi_y\|^2\geq
\frac{\bar{C}_5}{4}\int\frac{\bar{p}}{\bar{v}}\phi_y^2dy.
\end{eqnarray*}
Let
\begin{eqnarray*}\label{B4.115}
E_4=\bar{C}_4\v^{-2}E_3+\bar{C}_5\v^{-2}\int\big(\frac{2\mu(\bar{\theta})}{3\bar{v}}
\phi_y^2-\phi_y\psi_1\big)dy,
\end{eqnarray*}
\begin{eqnarray*}\label{B4.116}
K_4=\frac{1}8\bar{C}_4\v^{-2}K_3+\frac{\bar{C}_5}{4}\v^{-2}\int\frac{\bar{p}}{\bar{v}}
\phi_y^2dy+\v^{-2}\|(\phi_\tau,\psi_\tau,\zeta_\tau)\|^2.
\end{eqnarray*}
Then from \eqref{B4.93-1},
\eqref{B4.109}, \eqref{B4.111} and \eqref{B4.112}, we have the following estimate on the $(\phi,\psi,\zeta)$
\begin{eqnarray}\label{B4.117}
&&E_{4\tau}+K_4\leq C_9\delta\v^2(1+t)^{-1}E_4+C_9\d\v^2(1+t)^{-2}+C_9\v^{-2}\sum_{|\alpha|=2}\|\partial^\alpha(\phi,\psi,\zeta)\|^2 \nonumber\\
&&~~~~~~~~~~~~~~~~~~~~~~~~
+C_9(1+t)^{-1}\int|\hat\t_y|(b_1^2+b_3^2)dy+C_9\v^{-2}J_3,
\end{eqnarray}
where $J_3$ is defined in \eqref{B4.111-1}.

\

Define
\begin{eqnarray}\label{B4.115-1}
E_5=E_4+\v\int\int\frac{|\tilde{G}|^2}{2M_*}d\xi dy,~~~
K_5=K_4+\frac{\bar{\sigma}}{4}
\v\int\int\frac{\nu(|\xi|)}{M_*}|\tilde{G}|^2d\xi dy.
\end{eqnarray}
Then from \eqref{B4.117} and \eqref{B4.81}, one has
\begin{equation*}\label{B4.119}
\begin{array}{ll}
\di E_{5\tau}+K_5\leq C_{10}\delta\v^2(1+t)^{-1}E_5+C_{10}\d\v^2(1+t)^{-\f32}
+C_{10}\v^{-2}\sum_{|\alpha|=2}\|\partial^\alpha(\phi,\psi,\zeta)\|^2 \\
\di ~~~~~+C_{10}(\d+\l_0)\v^2\sum_{|\alpha|=1}
\int\int\frac{\nu(|\xi|)}{M_*}|\partial^\alpha\tilde{G}|^2d\xi dy+C_{10}\sum_{|\alpha|=2}\int\int\frac{\nu(|\xi|)}{M_*}|\partial^\alpha\tilde{G}|^2d\xi dy\\
\di ~~~~~+C_{10}(1+t)^{-1}\int|\hat\t_y|(b_1^2+b_3^2)dy.
\end{array}
\end{equation*}
Next we derive the higher order derivative estimate. Applying
$\partial_y$ to \eqref{B4.90} yields that
\begin{eqnarray}\label{B4.90-1}
\left\{
\begin{array}{l}
\di \phi_{y\tau}-\psi_{1yy}=-\frac{3}{5}\v^2\hat{N}_{1yy},\\[2mm]
\di \psi_{1y\tau}+\f{p_+}{\bar v}\zeta_{yy}-\f{p_+}{\bar v}\phi_{yy}=
\frac{4\v}3(\frac{\mu(\theta)}{v}u_{1y}-\frac{\mu(\bar{\theta})}{\bar{v}}\bar{u}_{1y})_{yy}+Q_9,\\[2mm] \di \psi_{iy\tau}=\v(\frac{\mu(\theta)}{v}u_{iy}-\frac{\mu(\bar{\theta})}{\bar{v}}\bar{u}_{iy})_{yy}
+Q_{8+i},i=2,3,\\[2mm]
\di \zeta_{y\tau}+p_+\psi_{1yy}=(\frac{\k(\theta)}{v}\theta_{y}
-\frac{\k(\bar{\theta})}{\bar{v}}\bar{\theta}_{y})_y+\frac2{5}\v^2\hat{N}_{1y}+Q_{12},
\end{array}
\right.
\end{eqnarray}
where
\begin{align*}
& Q_9=\frac{p-\bar{p}_+}{v}\phi_{yy}+(\f{p}{v}-\frac{\bar{p}}{v})\phi_{yy}
+O(1)(|\bar{v}_{yy}|\cdot|(\phi,\zeta)|+|\phi\zeta_{yy}|)\nonumber\\
&~~~~~~~-\f43\Big(\f{1}{v^2}v_y\t_y-\f{1}{\bar{v}^2}\bar{v}_y\bar{\t}_y\Big)
+\f43\Big(\f{\t}{v^3}v_y^2-\f{1}{{\bar\t}^3}\bar{v}_y^2\Big)
-\v\int\xi_1^2\Theta_{1yy}d\xi-\bar R_{1yy},\nonumber\\
& Q_{i+8}=-\v\int\xi_1\xi_{i}\Theta_{1yy}d\xi-\bar{R}_{iyy},i=2,3,\\
& Q_{12}=-\v\bar{u}_{1yy}(p-\bar p)-\v(p_yu_y-\bar{p}_y\bar{u}_y)+(p_+-\bar p)\psi_{1yy}+Q_{13y}\nonumber\\
&~~~~~~~~~~~~~~~~~~~~~~~~~~-\frac12\v\int\xi_1|\xi|\Theta_{1yy}^2d\xi
+\sum_{i=1}^3\v^2 (u_i\int\xi_1\xi_i\Theta_{1y}d\xi)_y,\\
& Q_{13}=\frac43\frac{\mu(\theta)}{v}\v^2u_{1y}^2+\v^2\sum_{i=2}^3\frac{\mu(\theta)}{v}u_{iy}^2
-\bar{H}_{1y}-\bar{R}_{4y}+\frac12(|\v\bar{u}|^2)_\tau+\v\bar{p}_y\bar{u}_1.
\end{align*}
Multiplying $\eqref{B4.90-1}_1$ by $p_+\phi_y$, $\eqref{B4.90-1}_2$ by $\bar v\psi_{1y}$,  $\eqref{B4.90-1}_3$ by $\psi_{iy}$, $\eqref{B4.90-1}_4$ by $\zeta_{y}$, we have
\begin{align}\label{B4.124}
&\Big[\int\big(\frac{p_+}{2}\phi_y^2+\frac{\bar{v}}{2}\psi_{1y}^2+\sum_{i=2}^3\psi_{iy}^2
+\frac{1}{2}\zeta_y^2\big)dy\Big]_\tau
+\f34\int\Big[\frac{4\mu(\theta)}{3v}\psi_{1yy}^2+\sum_{i=2}^3\frac{\mu(\theta)}{v}\psi_{iyy}^2
+\frac{\k(\theta)}{v}\zeta_{yy}^2\Big]dy\nonumber\\
& \leq C\d\v(1+t)^{-\f12}\|(\phi_y,\psi_y,\zeta_y)\|^2+C\d\v^3(1+t)^{-\f32}\|(\phi,\psi)\|^2
+C\|(\phi,\psi)\|\|(\phi_y,\psi_y,\zeta_y)\|^{3}\nonumber\\
&~~~~~~~~~~+C\|(\phi_y,\psi_y,\zeta_y)\|^{\f{10}{3}}+C\d\v^6(1+t)^{-3}
+C\v^2\int\Big|\int|\xi|^3\Theta_{1y}d\xi\Big|^2 dy\nonumber\\
& \leq C\d\v(1+t)^{-\f12}\sum_{|\alpha|=1}\|\partial^\alpha(\phi,\psi,\zeta)\|^2+C\d\v^3(1+t)^{-\f32}\|(\phi,\psi)\|^2
+C\|(\phi,\psi)\|\|(\phi_y,\psi_y,\zeta_y)\|^{3}\nonumber\\
& ~~~~~~~~~~+C\|(\phi_y,\psi_y,\zeta_y)\|^{\f{10}{3}}+C\delta\v^5(1+t)^{-\f52}
+C\d\v^2\sum_{|\alpha|=2}\|\partial^\alpha(\phi,\psi,\zeta)\|^2+CJ_3,
\end{align}
where we have used  \eqref{B4.111} in the last inequality.

Let
\begin{equation*}\label{B4.125}
E_6=\int\Big[\frac{p_+}{2}\phi_y^2+\frac{\bar{v}}{2}\psi_{1y}^2+\sum_{i=2}^3\psi_{iy}^2
+\frac{1}{2}\zeta_y^2\Big]dy,~~~K_6=\int\Big[\frac{4\mu(\theta)}{3v}\psi_{1yy}^2+\sum_{i=2}^3\frac{\mu(\theta)}{v}\psi_{iyy}^2
+\frac{\k(\theta)}{v}\zeta_{yy}^2\Big]dy,
\end{equation*}
then \eqref{B4.124} implies
\begin{eqnarray}\label{B4.126}
&&E_{6\tau}+\f12K_6\leq C_{11}\d\v(1+t)^{-\f12}\sum_{|\alpha|=1}\|\partial^\alpha(\phi,\psi,\zeta)\|^2
+C_{11}\d\v^3(1+t)^{-\f32}\|(\phi,\psi,\zeta)\|^2+C_{11}\delta\v^5(1+t)^{-\f52}\nonumber\\
&&~~+C_{11}\|(\phi,\psi)\|\|(\phi_y,\psi_y,\zeta_y)\|^{3}
+C_{11}\|(\phi_y,\psi_y,\zeta_y)\|^{\f{10}{3}}
+C_{11}\d\v^2\sum_{|\alpha|=2}\|\partial^\alpha(\phi,\psi,\zeta)\|^2+C_{11}J_3.
\end{eqnarray}

To get the estimate on $\phi_{yy}$, we use the momentum
equation $\eqref{B1.14}_2$. Applying $\partial_y$ on $\eqref{B1.14}_2$, it holds that
\begin{equation}\label{B4.127}
\psi_{1y\tau}+(p-\bar{p})_{yy}+\v\bar{u}_{1y\tau}+\bar{p}_{yy}=-\v\int\xi_1^2G_{yy}d\xi.
\end{equation}
Note that
\begin{equation*}\label{B4.128}
(p-\bar{p})_{yy}=-\frac{p}{v}\phi_{yy}+\frac{2}{3v}\zeta_{yy}-\frac1v(p-\bar{p})\bar{v}_{yy}
-\frac{\phi}{v}\bar{p}_{yy}-\frac{2v_y}{v}(p-\bar{p})_y-\frac{2\bar{p}_y}{v}\phi_y,
\end{equation*}
then multiplying \eqref{B4.127} by $-\phi_{yy}$ and integrating the
reduced equation with respect to $y$  give that
\begin{eqnarray}\label{B4.129}
&&(-\int\psi_{1y}\phi_{yy}dy)_\tau+\int\frac{p}{2v}\phi_{yy}^2dy\nonumber\\
&&\leq C_{12}K_6+C_{12}\d\v(1+t)^{-\f12}\|(\phi,\psi,\zeta)_y\|^2
+C_{12}\d\v^3(1+t)^{-\f32}\|(\phi,\psi,\zeta)\|^2+C_{12}\delta\v^5(1+t)^{-\f52}\nonumber\\
&&~~~~~
+C_{12}\|(\phi_y,\psi_y,\zeta_y)\|^{\f{10}{3}}
+C_{12}\v^2\int\int\f{\nu(|\xi|)}{M_\ast}|\tilde{G}_{yy}|^2d\xi dy.
\end{eqnarray}
To estimate $(\phi,\psi,\zeta)_{y\tau}$ and $(\phi,\psi,\zeta)_{\tau\tau}$,
we also use the original fluid-type equation \eqref{B1.14}. Here we only
consider the term $\di \int\psi_{1y\tau}^2dy$ because the other terms can
be estimated similarly. It follows from $\eqref{B1.14}_2$ that
\begin{equation}\label{B4.130}
\psi_{1y\tau}=-(p-\bar{p})_{yy}-\v \bar{u}_{1y\tau}-\bar{p}_{yy}-\v\int\xi_1^2G_{yy}d\xi.
\end{equation}
By \eqref{B4.130} and using the Cauchy inequality, it holds that
\begin{eqnarray}\label{B4.131}
&&\|\psi_{1y\tau}\|^2\leq C_{13}(K_6+\|\phi_{yy}\|^2)+C_{13}\d\v(1+t)^{-\f12}\|(\phi,\psi,\zeta)_y\|^2
+C_{13}\d\v^3(1+t)^{-\f32}\|(\phi,\psi,\zeta)\|^2\nonumber\\
&&~~~~~~~~~~~~+C_{13}\delta\v^5(1+t)^{-\f52}
+C_{13}\|(\phi_y,\psi_y,\zeta_y)\|^{\f{10}{3}}
+C_{13}\v^2\int\int\f{\nu(|\xi|)}{M_\ast}|\tilde{G}_{yy}|^2d\xi dy.
\end{eqnarray}
Let $\bar{C}_6$ and $\bar{C}_7$ be suitably large constants, then it follows from \eqref{B4.126}, \eqref{B4.129} and \eqref{B4.131} that
\begin{eqnarray}\label{B4.132}
&&\bar{C}_7\Big(\bar{C}_6E_6-\int\psi_{1y}\phi_{yy}dy\Big)_\tau
+\sum_{|\a|=2}\|\partial^\a(\phi,\psi,\zeta)\|^2\nonumber\\
&&\leq C_{14}\d\v(1+t)^{-\f12}\sum_{|\alpha|=1}\|\partial^\alpha(\phi,\psi,\zeta)\|^2
+C_{14}\d\v^3(1+t)^{-\f32}\|(\phi,\psi,\zeta)\|^2+C_{14}\delta\v^5(1+t)^{-\f52}\nonumber\\
&&~~~~~~+C_{14}\|(\phi,\psi)\|\|(\phi_y,\psi_y,\zeta_y)\|^{3}
+C_{14}\|(\phi_y,\psi_y,\zeta_y)\|^{\f{10}{3}}+C_{14}J_3,
\end{eqnarray}
where $J_3$ is defined in \eqref{B4.111-1}.

\

To close the a priori argument, we need to estimate the non-fluid component $\partial^\a\tilde{G}, |\a|=1,2$. Applying $\partial_y$ on \eqref{B1.22-1}, we have
\begin{eqnarray}\label{B4.134-1}
&&v\tilde{G}_{y\tau}-vL_M\tilde{G}_{y}=-v_y\tilde{G}_{\tau}
+v_yL_M\tilde{G}+2Q(M_y,\tilde{G})-\Big\{\frac1{R\theta}P_1[\xi_1
(\frac{|\xi-\v u|^2}{2\theta}\f1\v\zeta_y+\xi\cdot\f1\v\psi_y)M]\Big\}_y\nonumber\\
&&~~~~~~~~~~~~~~~~~~~~~~~~+\Big\{\v u_1G_y-P_1(\xi_1G_y)+\v vQ(G,G)-v\bar{G}_\tau\Big\}_y.
\end{eqnarray}
Multiplying \eqref{B4.134-1} by $\frac{\tilde{G}_y}{M_*}$, then integrating the reduced equation with respect to $\xi$ and $y$ and using the
Cauchy inequality and Lemmas \ref{lem4.1}-\ref{lem4.3}, we have
\begin{eqnarray*}\label{B4.133-1}
&&\Big(\int\int\frac{v|\tilde{G}_y|^2}{2M_*}d\xi
dy\Big)_\tau+\frac{3\bar{\sigma}}{4}\int\int\frac{\nu(|\xi|)}{M_*}|\tilde{G}_y|^2d\xi dy\nonumber\\
&&\leq C\sum_{|\a|=2}\int\int\frac{\nu(|\xi|)}{M_*}|\partial^\a\tilde{G}|^2d\xi dy+
C_3\delta\v^3(1+t)^{-5/2}+C_3\v^{-2}\sum_{|\alpha|=2}\|\partial^\alpha(\phi,\psi,\zeta)\|^2\nonumber\\
&&~~~~~~~~+C_3\d(1+t)^{-1}\sum_{|\alpha|=1}\|\partial^\alpha(\phi,\psi,\zeta)\|^2
+C_3\v^{-2}\sum_{|\alpha|=1}\|\partial^\alpha(\phi,\psi,\zeta)\|^6\nonumber\\
&&~~~~~~~~+C_3\sum_{|\alpha|=1}\|\partial^\alpha(v,u,\t)\|^2_{L^\infty}
\int\int\frac{\nu(|\xi|)}{M_*}|\tilde{G}|^2d\xi dy.
\end{eqnarray*}
Similarly, we can obtain the  estimate for $\tilde{G}_{\tau}$. Hence, one obtains that
\begin{eqnarray}\label{B4.133}
&&\Big(\sum_{|\a|=1}\int\int\frac{v|\partial^\a\tilde{G}|^2}{2M_*}d\xi
dy\Big)_\tau+\frac{3\bar{\sigma}}{4}\sum_{|\a|=1}\int\int\frac{\nu(|\xi|)}{M_*}|\partial^\a\tilde{G}|^2d\xi dy\nonumber\\
&&\leq C\sum_{|\a|=2}\int\int\frac{\nu(|\xi|)}{M_*}|\partial^\a\tilde{G}|^2d\xi dy+
C_3\delta\v^3(1+t)^{-5/2}+C_3\v^{-2}\sum_{|\alpha|=2}\|\partial^\alpha(\phi,\psi,\zeta)\|^2\nonumber\\
&&~~~~~~~~+C_3\d(1+t)^{-1}\sum_{|\alpha|=1}\|\partial^\alpha(\phi,\psi,\zeta)\|^2
+C_3\v^{-2}\sum_{|\alpha|=1}\|\partial^\alpha(\phi,\psi,\zeta)\|^6\nonumber\\
&&~~~~~~~~+C_3\sum_{|\alpha|=1}\|\partial^\alpha(v,u,\t)\|^2_{L^\infty}
\int\int\frac{\nu(|\xi|)}{M_*}|\tilde{G}|^2d\xi dy.
\end{eqnarray}

Finally, we need the highest order estimate  to control $\di \sum_{|\a|=2}\int\int\frac{\nu(|\xi|)}{M_*}|\partial^\a\tilde{G}|^2d\xi dy$ and $\di \int\psi_{1y}\phi_{yy}dy$ in \eqref{B4.132}. To estimate $\di \int\psi_{1y}\phi_{yy}dy$, it is sufficient to study the a priori estimate for $\di \sum_{|\a|=2}\int\int\frac{v|\partial^\a\tilde{f}|^2}{2M_*}d\xi dy$ due to \eqref{B4.16} and \eqref{B4.17}. Applying $\partial^\a,~|\a|=2$ to \eqref{B4.11}, one obtains that
\begin{eqnarray}\label{B4.11-1}
&&v\partial^\a\tilde{f}_\tau-\v vL_{M}\partial^\a\tilde{G}-\v u_1\partial^\a\tilde{f}_y+\xi_1\partial^\a\tilde{f}_y
=-\partial^\a v\tilde{f}_\tau+\v\partial^\a{u}_1\tilde{f}_y
-\sum_{|\beta|=1}\Big[\partial^{\a-\beta}v\partial^\beta\tilde{f}_\tau-\v\partial^{\a-\beta}u_1\partial^\beta \tilde{f}_y\Big]\nonumber\\
&&~~~~+\v\partial^\a\big[vL_{M}\bar{G}-\bar{v}L_{\bar{M}}\bar{G}_0\big]
+\v^2\partial^\a\big[vQ(G,G)-\bar{v}Q(\bar{G}_0,\bar{G}_0)\big]
+\partial^\a\Big[-\phi\bar{f}_\tau+\psi\bar{f}_y-\v v\bar{R}_{\bar f}\Big].
\end{eqnarray}
Multiplying \eqref{B4.11-1} by $\di \frac{\partial^\a\tilde{f}}{M_*}$, integrating the reduced equation with respect to $\xi$ and $y$
 and using the Cauchy inequality and Lemmas \ref{lem4.1}-\ref{lem4.4}, similar
to the argument used in \cite{Huang-Wang-Yang}, one gets that
\begin{eqnarray}\label{B4.134}
&&\Big(\sum_{|\a|=2}\int\int\frac{v|\partial^\a\tilde{f}|^2}{2M_*}d\xi
dy\Big)_\tau+\frac{3\bar{\sigma}}{4}\sum_{|\a|=2}\v^2\int\int\frac{\nu(|\xi|)}{M_*}|\partial^\a\tilde{G}|^2d\xi dy\nonumber\\
&&\leq
C_3\delta\v^5(1+t)^{-5/2}+C_3(\d+\eta_0+\l_0^{\f14})\sum_{|\alpha|=2}\|\partial^\alpha(\phi,\psi,\zeta)\|^2
+C\sum_{|\alpha|=1}\|\partial^\alpha(\phi,\psi,\zeta)\|^{\f{10}{3}}\nonumber\\
&&~~~~~~+C_3\d\v(1+t)^{-\f12}\sum_{|\alpha|=1}\|\partial^\alpha(\phi,\psi,\zeta)\|^2
+C_3\d\v^3(1+t)^{-\f32}\|(\phi,\psi,\zeta)\|^2
\nonumber\\
&&~~~~~~~~+\frac{C_3}{\l_0}\Big[\d^2\v^4(1+t)^{-2}+\sum_{|\alpha|=1}^2\|\partial^\alpha(\phi,\psi,\zeta)\|^4\Big]
\int\int\frac{\nu(|\xi|)}{M_*}|\tilde{G}|^2d\xi dy\nonumber\\
&&~~~~~~~~~~+C_3(\d+\eta_0+\l_0^{\f14})\v^2\int\int\frac{\nu(|\xi|)}{M_*}|\tilde{G}_y|^2d\xi dy.
\end{eqnarray}
Choose large constants $\bar{C}_8>1$ and $\bar{C}_9>1$
such that
\begin{eqnarray*}\label{B4.135}
&&E_7=\f{\bar{C}_8\bar{C}_7}{\v^3}\Big(\bar{C}_6E_6-\int\psi_{1y}\phi_{yy}dy\Big)
+\f1{\v}\sum_{|\a|=1}\int\int\frac{v|\partial^\a\tilde{G}|^2}{2M_*}d\xi dy
+\f{\bar{C}_9}{\v^3}\sum_{|\a|=2}\int\int\frac{v|\partial^\a\tilde{f}|^2}{2M_*}d\xi dy\nonumber\\
&&\geq \frac{c_1}{\v^3}\Big(\|(\phi,\psi,\zeta)_y\|^2+\sum_{|\a|=2}\|\partial^\a(\phi,\psi,\zeta)\|^2
+\sum_{|\a|=2}\int\int\frac{|\partial^\a \tilde f|^2}{M_*}d\xi dy\Big)\nonumber\\
&&~~~~~~~~~~~~~~~~~~~~~~~~~~~~~~~+\f{c_1}{\v}\sum_{|\a|=1}\int\int\frac{|\partial^\a\tilde{G}|^2}{2M_*}d\xi dy-C\d(1+t)^{-\f32}.
\end{eqnarray*}
Let
\begin{equation*}\label{B4.136}
K_7=\f{\bar{C}_8}{4\v^3}\sum_{|\a|=2}\|\partial^\a(\phi,\psi,\zeta)\|^2
+\f{\bar\sigma}{4\v}\sum_{1\leq |\a|\leq2}\int\int\frac{\nu(|\xi|)}{M_*}|\partial^\a\tilde{G}|^2d\xi dy.
\end{equation*}
Then from \eqref{B4.132}, \eqref{B4.133} and \eqref{B4.134}, one obtains  that
\begin{eqnarray}\label{B4.137}
&&E_{7\tau}+K_7\leq
C\delta\v^2(1+t)^{-5/2}+C\f{1}{\v^3}\sum_{|\alpha|=1}\|\partial^\alpha(\phi,\psi,\zeta)\|^{\f{10}{3}}
+C\d(1+t)^{-\f32}\|(\phi,\psi,\zeta)\|^2\nonumber\\
&&~+C\Big[\d(1+t)^{-\f12}+\f{1}{\v}\|(\phi,\psi)\|\cdot\|(\phi_y,\psi_y,\zeta_y)\|\Big]\sum_{|\alpha|=1}\f{1}{\v^2}\|\partial^\alpha(\phi,\psi,\zeta)\|^2
\\
&&~+C\Big[\d\v(1+t)^{-1}+\sum_{|\alpha|=1}^2\big(\f1{\v}\|\partial^\alpha(\phi,\psi,\zeta)\|^2
+\v\int\int\frac{|\partial^\a\tilde{G}|^2}{M_*}d\xi dy\big)\Big]
\int\int\frac{\nu(|\xi|)}{M_*}|\tilde{G}|^2d\xi dy.\nonumber
\end{eqnarray}

From \eqref{B4.117}, \eqref{B4.115-1} and \eqref{B4.137} and using the smallness of $\d,\l_0$ and $\v$, we have
\begin{eqnarray}\label{B4.138}
&&(E_4+E_7)_\tau+\f12(K_4+K_7)\leq C\d\v^2(1+t)^{-1}E_4+C\d\v^2(1+t)^{-2}+C(1+t)^{-1}\int|\hat\t_y|(b_1^2+b_3^2)dy\nonumber\\
&&~~+C\Big[\d\v(1+t)^{-1}+\sum_{|\alpha|=1}^2\big(\f1{\v}\|\partial^\alpha(\phi,\psi,\zeta)\|^2
+\v\int\int\frac{|\partial^\a\tilde{G}|^2}{M_*}d\xi dy\big)\Big]
\int\int\frac{\nu(|\xi|)}{M_*}|\tilde{G}|^2d\xi dy,
\end{eqnarray}
and
\begin{eqnarray}\label{B4.139}
&&(E_5+E_7)_\tau+\f12(K_5+K_7)\nonumber\\
&&\leq C\d\v^2(1+t)^{-1}E_5+C\d\v^2(1+t)^{-\f32}+C(1+t)^{-1}\int|\hat\t_y|(b_1^2+b_3^2)dy.
\end{eqnarray}

\section{The Proof of Main Result}

For a suitable large constant $\bar{C}_9$, by combining \eqref{B4.87} and \eqref{B4.139} and using  the smallness of $\d,\l_0$ and $\v$,  we have
\begin{eqnarray*}\label{B4.141}
E_{8\tau}+K_8\leq C\d\v^2(1+\v^2\tau)^{-1}E_8+C\d\v^2(1+\v^2\tau)^{-1},
\end{eqnarray*}
where
\begin{eqnarray*}\label{B4.142}
E_8=\bar{C}_9E_2+E_5+E_7,~~~~K_8=\f14(K_2+K_5+K_7)+\int|\hat\t_y|(b_1^2+b_3^2)dy.
\end{eqnarray*}
Note that
\begin{eqnarray*}\label{B4.143}
&&E_8\geq \|(\Phi,\Psi,W)\|^2+\Big\{\f{c_2}{\v^2}\|(\phi,\psi,\zeta)\|^2
+\v\int\int\frac{v|\tilde{G}|^2}{M_*}d\xi dy\Big\}\nonumber\\
&&~~~~~~~+\bigg\{\frac{c_1}{\v^3}\Big(\|(\phi,\psi,\zeta)_y\|^2+\sum_{|\a|=2}\|\partial^\a(\phi,\psi,\zeta)\|^2
+\sum_{|\a|=2}\int\int\frac{|\partial^\a f|^2}{M_*}d\xi dy\Big)\nonumber\\
&&~~~~~~~~~~~~~~~~~~~~~~~~~~~~~~~+\f{c_1}{\v}\sum_{|\a|=1}\int\int\frac{|\partial^\a\tilde{G}|^2}{2M_*}d\xi dy-C\d(1+\v^2\tau)^{-\f32}\bigg\},
\end{eqnarray*}
\begin{eqnarray*}\label{B4.144}
&&K_8\geq \sum_{|\beta|=1}\|\partial^\beta(\Phi,\Psi,W)\|^2
+c_2\bigg\{\frac{1}{\v^2}\sum_{|\a|=1}\|\partial^\a(\phi,\psi,\zeta)\|^2
+\v\int\int\frac{\nu(|\xi|)}{2M_*}|\tilde{G}|^2d\xi dy\bigg\}\nonumber\\
&&~~~~~~~+\bigg\{\frac{c_1}{\v^3}\sum_{|\a|=2}\|\partial^\a(\phi,\psi,\zeta)\|^2
+\f1{\v}\sum_{|\a|=1, 2}\int\int\frac{\nu(|\xi|)}{M_*}|\partial^\a\tilde{G}|^2d\xi dy\bigg\}+\f12\int|\hat\t_y|(b_1^2+b_3^2)dy,
\end{eqnarray*}
and
\begin{eqnarray}\label{B4.149-1}
\v^2E_7\leq C\d\v^2(1+\v^2\tau)^{-\f32}+C(K_4+K_7),~~\mbox{and} ~\v^2(E_5+E_7)\leq C\d\v^2(1+\v^2\tau)^{-\f32}+CK_8.
\end{eqnarray}

Then the Gronwall inequality yields that
\begin{eqnarray}\label{B4.145}
E_8\leq C\sqrt{\delta}(1+\v^2\tau)^{C_0\sqrt{\delta}},~~~~~~~~
\int_0^\tau K_8ds\leq
C\sqrt{\delta}(1+\v^2\tau)^{C_0\sqrt{\delta}}.
\end{eqnarray}
Hence, it holds that
\begin{eqnarray}\label{B4.146}
\|(\Phi,\Psi,W)\|^2\le
C\sqrt\delta(1+\v^2\tau)^{C_0\sqrt{\delta}}.
\end{eqnarray}
Multiplying \eqref{B4.139} by $(1+\v^2\tau)$ gives
\begin{eqnarray}\label{B4.147}
&&[(1+\v^2\tau)(E_5+E_7)]_\tau+\f12(1+\v^2\tau)(K_5+K_7)\nonumber\\
&&\leq C\v^2(E_5+E_7)+C\d\v^2(1+\v^2\tau)^{-\f12}+C_9\int|\hat\t_y|(b_1^2+b_3^2)dy.
\end{eqnarray}
Integrating \eqref{B4.147} with respect to $\tau$ and using \eqref{B4.145} and \eqref{B4.149-1}, one has that
\begin{eqnarray*}\label{B4.148}
&&(1+\v^2\tau)(E_5+E_7)+\int_{0}^{\tau}\f12(1+\v^2s)(K_5+K_7)ds\nonumber\\
&&\leq C\v^2\int_{0}^{\tau}(E_5+E_7)ds+C\sqrt{\d}(1+\v^2\tau)^{\f12}\nonumber\\
&&\leq C\sqrt{\d}(1+\v^2\tau)^{\f12}+C\int_{0}^{\tau}K_8ds\leq  C\sqrt{\d}(1+\v^2\tau)^{\f12},
\end{eqnarray*}
which yields
\begin{equation}\label{B4.148-1}
(E_5+E_7)\leq C\sqrt{\d}(1+\v^2\tau)^{-\f12}.
\end{equation}
In particular, one has
\begin{equation}\label{B4.148-2}
\v\int\int\frac{|\tilde{G}|^2}{M_*}d\xi dy\leq C\sqrt{\d}(1+\v^2\tau)^{-\f12}.
\end{equation}
On the other hand, multiplying \eqref{B4.139} by  $(1+\v^2\tau)^{\f12}$, it holds
\begin{eqnarray}\label{B4.149}
\int_{0}^{\tau}(1+\v^2s)^{\f12}(K_5+K_7)ds\leq C\sqrt{\d}(1+\v^2\tau)^{C_0\sqrt{\delta}}.
\end{eqnarray}

Multiplying \eqref{B4.138} by  $(1+\v^2\tau)$ and using \eqref{B4.148-1} and \eqref{B4.149-1}, one can obtain
\begin{eqnarray}\label{B4.150}
&&[(1+\v^2\tau)(E_4+E_7)]_\tau+\f12(1+\v^2\tau)(K_4+K_7)\leq C\v^2(E_4+E_7)+C\d\v^2(1+\v^2\tau)^{-1}\nonumber\\
&&~~~~~~~~~~~~~+C_9\int|\hat\t_y|(b_1^2+b_3^2)dy
+C\v(1+\v^2\tau)^{\f12}\int\int\frac{\nu(|\xi|)}{M_*}|\tilde{G}|^2d\xi dy\nonumber\\
&&~~~~~~~~~~~~~\leq C\d\v^2(1+\v^2\tau)^{-1}+CK_8+C(1+\v^2\tau)^{\f12}K_5.
\end{eqnarray}
Integrating \eqref{B4.150} with respect to $\tau$ and using \eqref{B4.145} and \eqref{B4.149}, one has
\begin{eqnarray}\label{B4.152}
(E_4+E_7)\leq C\sqrt{\d}(1+\v^2\tau)^{-1+C_0\sqrt{\delta}},
~~~~\int_{0}^{\tau}(1+\v^2s)(K_4+K_7)ds\leq C\sqrt{\d}(1+\v^2\tau)^{C_0\sqrt{\delta}}.
\end{eqnarray}
Therefore, it holds
\begin{eqnarray}\label{B4.152-1}
\|(\phi,\psi,\zeta)(\tau)\|^2\leq C\sqrt{\d}\v^2(1+\v^2\tau)^{-1+C_0\sqrt{\delta}}.
\end{eqnarray}

Multiplying \eqref{B4.137} by $(1+\v^2\tau)^{\f32-\vartheta}$ with $\vartheta>0$ in Theorem \ref{thm3.1} and using \eqref{B4.18},  \eqref{B4.149-1}, \eqref{B4.149}, \eqref{B4.152}  and the smallness of $\d$, one has
\begin{eqnarray}\label{B4.153}
&&[(1+\v^2\tau)^{\f32-\vartheta}E_{7}]_{\tau}=(\f32-\vartheta)(1+\v^2\tau)^{\f12-\vartheta}\v^2E_{7}
+(1+\v^2\tau)^{\f32-\vartheta}E_{7\tau}\nonumber\\
&&\leq C\d K_8+C(1+\v^2\tau)^{\f12-\vartheta}(K_4+K_7)+C\d(1+\v^2\tau)\sum_{|\a|=1}\f{1}{\v^2}\|\partial^\a(\phi,\psi,\zeta)\|^2
\nonumber\\
&&~~~~~+\frac{C}{\v^3}(1+\v^2\tau)^{\f32-\vartheta}\sum_{|\a|=1}\|\partial^\a(\phi,\psi,\zeta)\|^{\f{10}{3}}
+\frac{C}{\v^3}(1+\v^2\tau)^{\f32-\vartheta}\|(\phi,\psi,\zeta)\|\sum_{|\a|=1}\|\partial^\a(\phi,\psi,\zeta)\|^{3}
\nonumber\\
&&~~~~~+C\v(1+\v^2\tau)^{\f12-\vartheta+C_0\sqrt\d}\int\int\frac{\nu(|\xi|)}{M_*}|\tilde{G}|^2d\xi dy+C\d\v^2(1+\v^2\tau)^{-1-\vartheta}\nonumber\\
&&\leq C\d K_8+C\d(1+\v^2\tau)(K_4+K_7) +C(1+\v^2\tau)^{\f12-\vartheta+C_0\sqrt\d}K_5
+C\d\v^2(1+\v^2\tau)^{-1-\vartheta}\\[2mm]
&&~~ +\frac{C}{\v^3}(1+\v^2\tau)^{\f32-\vartheta}\sum_{|\a|=1}\|\partial^\a(\phi,\psi,\zeta)\|^{\f{10}{3}}
+\frac{C}{\v^3}(1+\v^2\tau)^{\f32-\vartheta}\|(\phi,\psi,\zeta)\|\sum_{|\a|=1}\|\partial^\a(\phi,\psi,\zeta)\|^{3}.\nonumber
\end{eqnarray}
By using \eqref{B4.14} and \eqref{B4.152}, one can get
\begin{eqnarray}\label{B4.154}
&&\frac{1}{\v^3}\int_{0}^{\tau}(1+\v^2s)^{\f32-\vartheta}\sum_{|\a|=1}\|\partial^\a(\phi,\psi,\zeta)\|^{\f{10}{3}}ds
+\frac{1}{\v^3}\int_{0}^{\tau}(1+\v^2s)^{\f32-\vartheta}\|(\phi,\psi,\zeta)\|\sum_{|\a|=1}\|\partial^\a(\phi,\psi,\zeta)\|^{3}ds\nonumber\\
&&\leq \frac{1}{\v} \int_{0}^{\tau}(1+\v^2s)^{\f32-\vartheta}\Big[(1+\v^2s)^{-\f23+\f23C_0\sqrt\d}
+(1+\v^2s)^{-1+C_0\sqrt\d}\Big]\sum_{|\a|=1}\|\partial^\a(\phi,\psi,\zeta)\|^2ds\nonumber\\
&&\leq \frac{1}{\v} \int_{0}^{\tau}(1+\v^2s)^\f12\sum_{|\a|=1}\|\partial^\a(\phi,\psi,\zeta)\|^2ds
\leq C\v\int_{0}^{\tau}(1+\v^2s)^\f12K_4ds\leq C\sqrt\d\v(1+\v^2\tau)^{C_0\sqrt\d},
\end{eqnarray}
provided that $C_0\sqrt\d\leq \vartheta$.
Thus  integrating \eqref{B4.153} over $[0,\tau]$ and using \eqref{B4.145}, \eqref{B4.149}, \eqref{B4.152} and \eqref{B4.154} yield that
\begin{eqnarray*}\label{B4.155}
E_7\leq C \sqrt\d(1+\v^2\tau)^{-\f32+\vartheta+C_0\sqrt\d},
\end{eqnarray*}
which immediately implies
\begin{eqnarray}\label{B4.155-1}
&&\frac{1}{\v^3}\Big(\|(\phi,\psi,\zeta)_y\|^2+\sum_{|\a|=2}\|\partial^\a(\phi,\psi,\zeta)\|^2
+\sum_{|\a|=2}\int\int\frac{|\partial^\a \tilde f|^2}{M_*}d\xi dy\Big)\nonumber\\
&&~~~~~~~~~~~~~~~~~~~+\f{1}{\v}\sum_{|\a|=1}\int\int\frac{|\partial^\a\tilde{G}|^2}{2M_*}d\xi dy
\leq C \sqrt\d(1+\v^2\tau)^{-\f32+\vartheta+C_0\sqrt\d}.
\end{eqnarray}

\noindent\textbf{Proof of Theorem \ref{thm3.1}}: Combining \eqref{B4.146}, \eqref{B4.148-2}, \eqref{B4.155-1} and \eqref{B4.152-1} and using the Sobolev inequality, it holds that
\begin{eqnarray}\label{B4.156-1}
&&\|(\Phi,\Psi,W)\|^2_{L^\infty}\leq C\|(\Phi,\Psi,W)\|\Big(\|(\phi,\psi,\zeta)\|+\d\v(1+\v^2\tau)^{-\f12}\Big)\leq
C\sqrt\d\v,
\end{eqnarray}
and
\begin{eqnarray}\label{B4.156-2}
\|\int\frac{|\tilde{G}|^2}{M_*}d\xi\|_{L^\infty} \leq
\Big(\int\int\frac{|\tilde{G}|^2}{M_*}d\xi dy\Big)^{\f12}\Big(\int\int\frac{|\tilde{G}_y|^2}{M_*}d\xi dy\Big)^{\f12}\leq  C\sqrt{\d}(1+\v^2\tau)^{-\f12}.
\end{eqnarray}
Therefore, \eqref{B4.148-2}, \eqref{B4.152-1}, \eqref{B4.155-1}, \eqref{B4.156-1} and \eqref{B4.156-2} verify the a priori assumption \eqref{B4.13} if we choose $\l_0=\d^{\f18}$.
Hence,  the proof of Theorem \ref{thm3.1} is completed. $\hfill\Box$

\

\noindent\textbf{Proof of Theorem \ref{thm2.1}}:
The proof of \eqref{B1.50-1} can be obtained directly from \eqref{B1.50-2} by using the transformation \eqref{B4.2} of the scaled variables $(y,\tau)$ and the original variables $(x,t)$. By combining \eqref{B1.50-1} and Sobolev inequality, \eqref{B1.50} can be derived immediately. Thus the proof of Theorem \ref{thm2.1} is completed. $\hfill\Box$

\vspace{1cm}

{\bf Acknowledgments.}
The authors would like to thank Renjun Duan for insightful discussion
during the third author's visit in the Chinese University of Hong Kong.
Feimin Huang is partially supported by National Basic Research Program
of China (973 Program) under Grant No. 2011CB808002 and by National Center for Mathematics
and Interdisciplinary Sciences, AMSS, CAS and the CAS Program for Cross $\&$ Cooperative
Team of the Science $\&$ Technology Innovation. Yi Wang is supported by National
Natural Sciences Foundation of China No. 10801128 and No. 11171326. Yong Wang is partially supported by National Natural
Sciences Foundation of China No. 11371064 and 11401565. Tong Yang is supported by the
General Research Fund of Hong Kong, CityU 103412.


\begin{thebibliography}{99}
\bibitem{Atkinson-Peletier}
\newblock  F. V. Atkinson and  L. A. Peletier,
\newblock  \emph{Similarity solutions of the nonlinear diffusion equation},
\newblock  Arch. Rat. Mech. Anal., \textbf{54} (1974), 373--392.

\bibitem{BGL}
\newblock  C. Bardos, F. Golse and D. Levermore,
\textit{Fluid dynamic limits of kinetic equations, I. Formal
derivations,}  \newblock J. Statis. Phys., \textbf{63}, (1991), 323-344; II.
\textit{Convergence proofs for the Boltzmann equation},
\newblock Comm. Pure Appl.
Math., \textbf{46}, (1993), 667-753, .

\bibitem{BGLY}
\newblock  C. Bardos, C. Levermore,  S. Ukai, and  T. Yang,
\textit{Kinetic equations: fluid dynamical limits and viscous
heating,} \newblock  Bull. Inst. Math. Acad. Sin. (N.S.) \textbf{3} (2008), 1-49.

\bibitem{BU}
\newblock C. Bardos and  S. Ukai,
\textit{The classical incompressible Navier-Stokes limit of the
Boltzmann equation},
\newblock Math. Models Methods Appl. Sci., \textbf{1} (1991),
235-257.


\bibitem{Boltzmann}
\newblock   L. Boltzmann, (translated by Stephen G. Brush),
\newblock ``Lectures on Gas Theory,''
\newblock Dover Publications, Inc. New York, 1964.

\bibitem{Caflisch}
\newblock   R. E. Caflisch,
   \emph{The fluid dynamical limit of the nonlinear Boltzmann
equation},
\newblock   Comm. Pure Appl. Math., \textbf{33} (1980), 491-508.



\bibitem{Cercignani-Illner-Pulvirenti}
\newblock   C. Cercignani, R. Illner and M. Pulvirenti,
 ``\textit{The Mathematical Theory of Dilute Gases},''
\newblock Springer-Verlag, Berlin, 1994.

\bibitem{CC}
\newblock S. Chapman and T. G. Cowling,
 ``\textit{The Mathematical Theory of Non-Uniform Gases},''
\newblock 3rd edition, Cambridge University Press, 1990.


\bibitem{Esposito1}
\newblock A. De Masi,   R. Esposito and J. L. Lebowitz,
\emph{ Incompressible Navier-Stokes and Euler limits of
the Boltzmann equation},
\newblock Comm. Pure Appl. Math. \textbf{42} (1989), no. 8, 1189-1214.

\bibitem{Duyn-Peletier}
\newblock C. T. Duyn and L. A. Peletier,
\newblock \emph{A class of similarity
solution of the nonlinear diffusion equation},
\newblock Nonlinear Analysis, T.M.A., \textbf{1} (1977), 223-233.

\bibitem{Diperna-Lions}
\newblock R. J. DiPerna and P. L. Lions,
\newblock \emph{On the Cauchy Problem for Boltzmann Equations: Global Existence and Weak Stability},
\newblock Annals of Mathematics, \textbf{130}  (1989),  no. 2, 321-366.



\bibitem{Esposito}
\newblock R. Esposito and M. Pulvirenti,
\newblock \emph{From particle to
fluids},
\newblock in ``Handbook of Mathematical Fluid Dynamics,'' Vol. \textbf{III},
North-Holland, Amsterdam, (2004), 1-82.

\bibitem{Guo2}
\newblock  R. Esposito, Y. Guo, C. Kim and R. Marra,
  \emph{Stationary solutions to the Boltzmann equation in the Hydrodynamic limit}.
\newblock Preprint,  arXiv:1502.05324.

\bibitem{Golse}
\newblock F. Golse,
\emph{The Boltzmann equation and its hydrodynamic
limits}, Evolutionary equations. Vol. II, 159-301, Handb. Differ.
Equ., Elsevier/North-Holland, Amsterdam, 2005.

\bibitem{GL}
\newblock F. Golse and  D. Levermore,
\emph{ Stokes-Fourier and acoustic limits for the Boltzmann equation:
convergence proofs},
\newblock Comm. Pure Appl. Math. \textbf{55} (2002), no. 3, 336-393.


\bibitem{GPS}
\newblock F. Golse, B. Perthame and C. Sulem,
 \emph{On a boundary layer
problem for the nonlinear Boltzmann equation},
\newblock Arch. Ration. Mech.
Anal., \textbf{103} (1986), 81-96.

\bibitem{GS1}
\newblock F. Golse and L. Saint-Raymond,
\emph{The Navier-Stokes limit of the Boltzmann equation for bounded
collision kernels},
\newblock Invent. Math. \textbf{155} (2004), no. 1, 81-161.


\bibitem{GS}
\newblock F. Golse and L. Saint-Raymond,
\emph{The incompressible Navier-Stokes
limit of the Boltzmann equation for hard cutoff potentials},
\newblock J. Math. Pures Appl., {\bf 91} (2009), 508-552.

\bibitem{Grad}
\newblock H. Grad,
 ``Asymptotic Theory of the
  Boltzmann Equation II,"
\newblock in ``Rarefied Gas Dynamics'' (J. A. Laurmann, ed.), Vol. \textbf{1}, Academic Press, New York, (1963), 26--59.

\bibitem{Guo}
\newblock  Y. Guo,
  \emph{The Boltzmann equation in the
whole space},
\newblock Indiana Univ. Math. J., \textbf{53} (2004), 1081-1094.

\bibitem{Guo1}
\newblock  Y. Guo,
  \emph{ Boltzmann diffusive limit beyond the Navier-Stokes approximation}.
\newblock Comm. Pure Appl. Math. \textbf{59} (2006), no. 5, 626-687.




\bibitem{Huang-Matsumura-Xin}
\newblock F. M. Huang, A. Matsumura and Z. P. Xin,
\newblock \emph{ Stability of Contact Discontinuities for the 1-D
Compressible Navier-Stokes Equations},
\newblock Arch. Rat.
Mech. Anal., \textbf{179}(2005), 55-77.

\bibitem{Huang-Wang-Wang-Yang}
\newblock F. M. Huang, Y. Wang, Y. Wang and T. Yang,
 \emph{ The  Limit of the Boltzmann
Equation to the  Euler Equations for Riemann Problems},
\newblock SIAM J. Math. Anal., \textbf{45} (2013), no. 3, 1741-1811.


\bibitem{Huang-Wang-Yang}
\newblock F. M. Huang, Y. Wang and T. Yang,
 \emph{Hydrodynamic limit of the Boltzmann equation with contact
discontinuities},
\newblock Comm. Math. Phy., \textbf{295} (2010), 293-326.


\bibitem{Huang-Wang-Yang-2}
\newblock F. M. Huang, Y. Wang and T. Yang,
\newblock \emph{Fluid dynamic limit to the Riemann solutions of  Euler equations:
 I. Superposition of rarefaction waves and contact discontinuity},
\newblock Kinet. Relat. Models, \textbf{3} (2010), 685-728.


\bibitem{Huang-Xin-Yang}
\newblock F. M. Huang, Z. P. Xin and T. Yang,
\newblock \emph{Contact discontinuity with general perturbations for gas motions},
\newblock Adv. Math., \textbf{219} (4) (2008), 1246-1297.

\bibitem{Huang-Yang} F. M. Huang and  T. Yang,\emph{ Stability of contact
discontinuity for the Boltzmann equation}, J. Differential Equations,
\textbf{229} (2006), 698-742.

\bibitem{Jang}
\newblock J. Jang,
\emph{ Vlasov-Maxwell-Boltzmann diffusive limit},
\newblock Arch. Ration. Mech. Anal. \textbf{194} (2009), no. 2, 531-584.

\bibitem{Xiong-Jiang}
\newblock N. Jiang  and L. J. Xiong,
\emph{Diffusive limit of the Boltzmann equation with fluid initial layer in the periodic domain}, to appear in SIAM, 2015.

\bibitem{Levermore}
\newblock D. Levermore and N. Masmoudi,
\emph{From the Boltzmann equation to an incompressible Navier-
Stokes-Fourier system},
\newblock Arch. Ration. Mech. Anal. \textbf{196} (2010), no. 3, 753-809.


\bibitem{L-M}
\newblock P.L. Lions and N. Masmoudi,
\emph{ From the Boltzmann equations to the equations of incompressible
fluid mechanics. I}
\newblock   Arch. Ration. Mech. Anal. \textbf{158} (2001), no. 3, 173-193.


\bibitem{L-M1}
\newblock P.L. Lions and N. Masmoudi,
\emph{From the Boltzmann equations to the equations of incompressible
fluid mechanics. II},
\newblock Arch. Ration. Mech. Anal. \textbf{158} (2001), no. 3, 195-211.

\bibitem{Liu-Zhao}
\newblock S. Q. Liu and H. J. Zhao,
\emph{ Diffusive expansion for solutions of the Boltzmann equation in the whole space},
\newblock J. Differential Equations, \textbf{250} (2011), no. 2, 623-674.

\bibitem{Liu-Yang-Yu}
\newblock T. Liu, T. Yang, and S. H.  Yu,
 \emph{Energy method for the Boltzmann equation},
\newblock Physica D, \textbf{188} (2004), 178--192.

\bibitem{Liu-Yang-Yu-Zhao}
\newblock T. Liu, T. Yang, S. H.  Yu and H. J. Zhao,
\newblock \emph{Nonlinear stability of rarefaction waves
for the Boltzmann equation},
\newblock Arch. Rat. Mech. Anal., \textbf{181} (2006), 333--371.

\bibitem{Liu-Yu}
\newblock T. Liu and S. H. Yu,
  \emph{Boltzmann  equation:
  Micro-macro decompositions  and positivity of shock profiles},
\newblock Commun. Math. Phys., \textbf{246} (2004), 133-179.

\bibitem{Maxwell} J. C. Maxwell, \emph{On the dynamical theory of gases}, Phil. Trans. Roy. Soc. London,
157(1866), 49-88.

\bibitem{Masmoudi}
\newblock  N. Masmoudi,
\emph{Examples of singular limits in hydrodynamics}, Handbook of differential equations: evolutionary equations. Vol. III, 195-275, Handb. Differ. Equ., Elsevier/North-Holland, Amsterdam, 2007.

\bibitem{Masmoudi-S}
\newblock N. Masmoudi and  L. Saint-Raymond,
\emph{From the Boltzmann equation to the Stokes-Fourier system
in a bounded domain},
\newblock  Comm. Pure Appl. Math. \textbf{56} (2003), no. 9, 1263-1293.

\bibitem{Saint}
\newblock L. Saint-Raymond, \emph{From the BGK model to the Navier-Stokes equations},  Ann. Sci. École Norm.
Sup.   \textbf{36} (2003), no. 2, 271-317.




\bibitem{Sone-2} Y. Sone, Molecular Gas Dynamics, Theory, Techniques, and
Applications, Birkh${\rm\ddot{a}}$user, Boston, 2006.


\bibitem{Ukai-1}
\newblock S. Ukai,
\emph{Solutions of the Boltzmann equation, Pattern and Waves} - Qualitative
Analysis of Nonlinear Differential Equations (eds. M.Mimura and T.Nishida), Studies of
Mathematics and Its Applications 18, 37-96, Kinokuniya-North-Holland, Tokyo, 1986.



\bibitem{Ukai-3} S. Ukai and K. Asano, \emph{ The Euler limit and initial layer of the Boltzmann equation},
Hokkaido Math. J, \textbf{12}(1983), 311-332.


\bibitem{Xin-Zeng}
\newblock Z. P. Xin and H. H.  Zeng,
\newblock \emph{Convergence to the rarefaction
waves for the nonlinear Boltzmann equation and compressible
Navier-Stokes equations},
\newblock J. Diff. Eqs., \textbf{249} (2010), 827--871.


\bibitem{Yu}
\newblock S. H. Yu,
\newblock \emph{Hydrodynamic limits with shock
waves of the Boltzmann equations},
\newblock Commun. Pure Appl. Math, \textbf{58} (2005), 409--443.

\end{thebibliography}
\end{document}